\documentclass{article}
\usepackage{amssymb}
\usepackage{amsmath}
\usepackage{latexsym}
\usepackage[mathscr]{euscript}
\usepackage{graphicx}
\usepackage{amscd}
\usepackage{amsthm}
\usepackage{fullpage}
\usepackage{subfigure}
\usepackage{wrapfig}

\usepackage{latexsym}

\newcommand{\la}{\to}
\newcommand{\hk}{\hookrightarrow}

\newcommand{\br}{{\mathbb R}}
\newcommand{\bz}{{\mathbb Z}}

\newcommand{\bp}{{\mathbb P}}

\newcommand{\cc}{\mathcal{C}}

\newcommand{\cn}{\mathcal{N}}
\newcommand{\cm}{\mathcal{M}}
\newcommand{\cp}{\mathcal{P}}
\newcommand{\cs}{\mathcal{S}}

\newcommand{\msa}{{\mathcal M}^{\sigma}_\Gamma (M;\vec{a})}

\newcommand{\p}{\partial}

\newcommand{\eps}{\epsilon}

\newcommand{\med}{\medskip}
\newcommand{\bfl}{\begin{flushleft}}
\newcommand{\efl}{\end{flushleft}}
\newcommand{\G}{\Gamma}

\newcommand{\ccg}{\cc_{\Gamma}}
\newcommand{\ag}{Aut (\Gamma)}
\newcommand{\aog}{Aut_0(\Gamma)}
\newcommand{\ago}{Aut(\Gamma_1)}
\newcommand{\agt}{Aut(\Gamma_2)}
\newcommand{\ccgm}{\cm_{\G} (M)}
\newcommand{\ccngm}{\cm^N_{\G}( M)}
\newcommand{\ccalgm}{\cm^\alpha_{\G}( M)}
\newcommand{\tccgm}{\tilde\cm_{\G}( M)}

\newcommand{\xr}{\xrightarrow}
\newcommand{\g}{\Gamma}
\newcommand{\ct}{\mathcal{T}}
\newcommand{\cv}{\mathcal{V}}
\newcommand{\cw}{\mathcal{W}}
\newcommand{\cgm}{\cm_\g}
\newcommand{\gn}{\tilde{\Gamma}}
\newcommand{\coker}{\rm coker\ }
\newcommand{\msan}{{\mathcal M}_{\g}^N(X;\vec{a})}
\newtheorem{theorem}{Theorem}
\newtheorem{definition}[theorem]{Definition}
\newtheorem{lemma}[theorem]{Lemma}
\newtheorem{corollary}[theorem]{Corollary}
\newtheorem{proposition}[theorem]{Proposition}

\begin{document}

\title{Morse field theory}
\author{Ralph Cohen  \thanks{The first author was partially supported  
by a grant from the NSF} \\ Department of Mathematics\\ Stanford  
University\\
Stanford, CA 94305  \and
Paul Norbury \thanks{The second author was partially supported by a  
grant from the Australian Academy of Science} \\ Department of  
Mathematics  and Statistics\\
University of Melbourne\\ Melbourne, Australia 3010}

\maketitle
\maketitle
\begin{abstract}
In this paper we define and study the moduli space of metric-graph- 
flows in a manifold $M$.
This is a space of smooth maps from a finite graph to $M$, which,  
when restricted to each edge,
is a gradient flow line of a smooth (and generically Morse) function  
on $M$.  Using the model of Gromov-Witten theory, with this moduli  
space replacing the space of stable holomorphic curves in a  
symplectic manifold, we obtain invariants, which
are (co)homology operations in $M$.  The invariants obtained in this  
setting are classical cohomology operations
such as cup product, Steenrod squares, and Stiefel-Whitney classes.   
We show that these operations satisfy invariance and gluing  
properties that  fit together to give the structure of a topological  
quantum field theory. By   considering equivariant operations with  
respect to the action of the automorphism group of the graph, the  
field theory has more structure. It is analogous to a homological  
conformal field theory.  In particular we show that classical  
relations such as the Adem relations and Cartan formulae are   
consequences of these field theoretic properties.
These operations are defined and studied using two different  
methods.  First, we use algebraic topological techniques to define  
appropriate virtual fundamental classes of these moduli spaces. This  
allows us to define the operations via the corresponding intersection  
numbers of the moduli space.    Secondly, we use geometric and  
analytic techniques to study the smoothness and compactness  
properties of these moduli spaces.  This will
allow us to define these   operations on the level of Morse-Smale  
chain complexes, by appropriately counting metric-graph-flows with   
particular boundary conditions.
\end{abstract}

\tableofcontents

\section*{Introduction}

In this paper we construct a moduli space of graphs $|\ccg|/Aut\G$
associated to a fixed oriented graph $\g$.  It is   built from a  
category
$\ccg$ in which the objects are graphs and morphisms are
homotopy equivalences.  We use this moduli space to study families  of
maps of graphs into a manifold,  which allows us to probe the
topology of the manifold.   The moduli space  is described in detail in
section 1.  For the moment it is best understood by its
following properties.  To each element of $|\ccg|$ is associated an
oriented, compact metric graph---where edges are given lengths---and an
orientation preserving homotopy equivalence from the metric graph to
the given graph $\g$ that collapses edges and vertices.  The space
$|\ccg|$ is contractible, and admits a free $\ag$ action, and hence  
the quotient is a model for the classifying space,
\[|\ccg|/Aut\g \simeq BAut\Gamma.\]
In particular when $\g$ has non-trivial automorphisms $|\ccg|/Aut\g$ has
non-trivial homology.

Given a fixed closed manifold $M$, we then thicken this moduli space  
by defining a space $\cs_\g $  whose points
are pairs, $(x, \mu)$, where $x \in |\ccg|$, and $\mu$ is a labeling  
of the edges of the graph by smooth functions on $M$.  We call $\cs_\G 
$ the space of metric-Morse structures on $M$, and define the moduli  
space of such structures
to be the quotient space,
$$
\cm_\G = \cs_\G/\ag.
$$
It will be easily seen that in thickening the moduli space $|\ccg|/Aut 
\G$ to define the moduli space of structures, $\cm_\G$,  we did not  
change the homotopy type, so that $\cm_\G \simeq B \ag$.   It is for  
this reason in our notation we suppress the dependence on $M$ of the  
moduli space $\cm_\G$.

   We can then define a moduli space $\cm_\G(M)$ of metric graph  
flows in $M$.  This space  consists of isomorphism classes of pairs, $ 
(\sigma, \gamma)$, where $\sigma \in \cs_\G$ is a metric-Morse  
structure on $\G$, and $\gamma $ is a continuous map from the graph  
to $M$, which, when restricted to a given edge, is a gradient flow  
line of the smooth function labeling that edge with respect to the  
parameterization of the edge coming from the orientation and  
metric.   Since $\cm_\G \simeq B \ag$,
   we can take as a representative of a homology class of the  
automorphism group $\ag$, a
family of metric-Morse structures on the graph $\G$.  When the  
structures in the moduli space of metric-graph-flows are restricted  
to vary in a family representing a fixed  homology class  of the  
automorphism group, we will have a finite dimensional moduli space.    
By studying the topology of this moduli space by two different  
methods (one using algebraic topology, to define Pontrjagin-Thom  
constructions and induced ``umkehr maps" in homology, and the other  
using geometry and analysis to understand the smoothness,  
transversality, and compactness properties of these moduli spaces),  
we obtain Gromov-Witten type invariants of  $M$.  For example,  the  
ring structure (cup product) in the cohomology of the manifold arises  
as such an invariant when the graph is a tree with three edges, and  
the family of structures is a single point.
    Further classical invariants such as Steenrod squares and Stiefel- 
Whitney classes of the manifold arise when we take higher dimensional  
families of structures representing nontrivial elements of the  
homology of the automorphism group.

The approach in this paper is designed specifically to deal with {\em
families} of metric-graphs mapping to manifolds.  Graphs are the  
essential objects
here.  Functions on the manifold are quite peripheral and do not even  
need to be
Morse.  The title ``Morse field theory'' primarily refers to integral  
flow lines
of  gradient vector fields on a manifold as well as the
Morse complex and cohomology operations defined on the Morse complex.

\med
The original goal of this project was to understand the Gromov-Witten  
formalism in the setting of Morse theory, where
the analysis is considerably easier.  In this model, the role of  
oriented, metric graphs fills the role of oriented surfaces with a  
conformal class of metric.   Maps from these graphs to manifolds that  
satisfy gradient flow equations
fill the role of $J$-holomorphic maps from a Riemann surface to a  
symplectic manifold.  This project took its original form in the work  
of M. Betz in his Ph.D thesis \cite{BetCat} written under the  
direction of the first author, and in the research announcement \cite 
{BCoGra}.  Similar constructions were discovered by Fukaya \cite 
{FukMor} in which he described his beautiful ideas on the $A_\infty$- 
structure of   Morse homotopy.  In particular those ideas have been  
used in the   work of Fukaya and Oh regarding deformations of $J$- 
holomorphic disks in the cotangent bundle of a manifold  \cite{fukoh}.

This present paper contains new ideas involving families of metric- 
Morse structures, as well as constructions of virtual fundamental  
classes of these moduli spaces, that allow  us to  the define    
equivariant
invariants, investigate their properties,  plus provide proofs of old  
ideas on non-equivariant invariants
\cite{BetCat,BCoGra,FukMor}.  As mentioned above, we show how to deal  
with families both
by using algebraic topological methods, and by using geometric and  
analytical techniques.  The algebraic topological techniques allow us  
to define generalized Pontrjagin-Thom constructions and resulting  
umkehr maps, which in turn allow the definition of virtual  
fundamental classes.  These techniques are based on the generalized  
Pontrjagin-Thom constructions defined by the first author and J.  
Klein in \cite{cohenumkehr}.    In particular these techniques allow  
us to avoid transversality (smoothness) and compactness issues that  
arise from the geometric viewpoint.  However, because the geometric  
viewpoint is quite important in its own right, in the second half of  
the paper we study these analytic issues, and prove the appropriate  
transversality and compactness results.  This allows a second  
definition of the invariants  that are defined  on  the level of   
Morse-Smale chain complexes, by counting the number of graph flows in  
a manifold that satisfy appropriate boundary conditions.

The moduli space $\cm_\G$  is somewhat  analogous to the moduli space  
$  \cm_g$
of Riemann surfaces homeomorphic to a given surface, and more
generally to $\cm_{g,n}$, the space of Riemann surfaces with $n$
marked points, when the graphs come equipped with marked, univalent  
vertices.  A
point in the Teichmuller space $\ct(\Sigma)$ of a topological surface
$\Sigma$
(with $n$ labeled points)
is a pair $(\Sigma',h)$ where
$\Sigma'$ is a complete hyperbolic surface and $h:\Sigma'\to\Sigma$ is
a homeomorphism well-defined up to isotopy.  Teichmuller space is
contractible and admits an action of the group of isotopy classes of
homeomorphisms of $\Sigma$, known as the mapping class group of
$\Sigma$.  The quotient of $\ct(\Sigma)$ by the mapping class group is
the moduli space of hyperbolic structures on $\Sigma$, which appears
in algebraic geometry as the moduli space $  \cm_g$ of Riemann  
surfaces.  In our setup, the
contractible space $|\ccg|$  plays the role of Teichmuller space,  
$Aut \, \g$
plays the role of the mapping class group, although unlike the  
mapping class group  it acts freely, and the
metric graph and homotopy equivalence $h:\g'\to\g$ is analogous to the
isotopy class of homeomorphism   $h:\Sigma'\to\Sigma$.

A further analogy between $\cm_\G$ and $\cm_{g,n}$ comes
from the fact that $\cm_{g,n}$ is homotopy equivalent to the moduli
space of metric ribbon graphs---finite graphs whose vertices are at  
least trivalent, and come equipped with a cyclic ordering of
(half-)edges at each vertex and lengths on edges---divided by
automorphisms \cite{penner}.   This analogy will be pursued further  
by the first author
in order to describe a Morse theoretic interpretation of string  
topology, and the relation between
string topology operations and $J$-holomorphic curves in the  
cotangent bundle \cite{cohenstring}.  A description of these ideas  
was given in \cite{cohenmontreal}.

A labeling $\mu$ of the edges of a graph in $|\ccg|$ by functions on  
$M$ is, in some sense, analogous
to choosing a compatible almost complex structure $J$  on a  
symplectic manifold.  In both cases the space
of choices of these structures is contractible, and each choice  
allows the definition of the relevant differential equations used to  
define a point in the moduli space (a $J$-holomorphic curve in the  
Gromov-Witten setting, and a gradient graph flow in our setting).

\med
Aside from the study of these moduli spaces of graphs and graph  
flows, and the resulting definition of the graph invariants  
(operations), the main  result of this paper is that these invariants  
fit together to define an appropriate field theory.
Recall that an $n$-dimensional  {\em topological quantum field  
theory} (TQFT) over a ring $R$   assigns to every closed $n-1$- 
dimensional manifold $N$, an $R$-module $Z(N)$ and to every cobordism  
$W$ from $N_1$ to $N_2$, (i.e $W$ is an $n$-manifold with boundary $ 
\p W = N_1 \sqcup N_2$), an operation
$$
\mu_W : Z(N_1) \to Z(N_2),
$$
which is a map of $R$-modules.  This structure is supposed to satisfy  
certain properties, the most important of which is gluing:  If $W_1$  
is a cobordism from $N_1$ to $N_2$, and $W_2$ is a cobordism from $N_2 
$ to $N_3$, $W = W_1  \cup_{N_2} W_2$ is the ``glued cobordism"  from  
$N_1$ to $N_3$ obtained by identifying the boundary components of $W_1 
$ and $W_2$ corresponding to $N_2$,  then we require
$$
\mu_{W_1  \cup_{N_2} W_2} = \mu_{W_2} \circ \mu_{W_1} : Z(N_1) \to Z 
(N_2) \to Z(N_3).
$$
These operations only depend on the diffeomorphism classes of the  
cobordisms. See \cite{atiyah} for details.

In the simplest case when $n=1$, we choose to relax
   the manifold condition, and think of graphs with univalent  
vertices as defining generalized cobordisms between
   zero dimensional manifolds.    These univalent vertices can be  
thought to have signs attached to them, according
   to whether the edge they lie on is oriented via an arrow that  
points toward or away from the vertex.  Alternatively we can think of  
these univalent vertices as ``incoming" or ``outgoing".

For a given manifold $M$, the Morse field theory functor $Z_M$ assigns
to each oriented point, $Z_M({\rm point})=H_*M$.
   Given a graph $\G$ with $p$ incoming and $q$ outgoing univalent  
vertices (i.e a generalized cobordism between $p$ points and $q$  
points), as well as a   homology class $\alpha \in H_*(\cgm) = H_*(B 
\ag)$, the graph invariants described above can be viewed as a  
homology operation
$$
q^\alpha_\G : H_*(M)^{\otimes p} \to H_*(M)^{\otimes q}.
$$
We prove that these operations satisfy gluing and a certain  
invariance properties.  This is the ``Morse field theory" of the title.
We remark that it is a well known folk theorem that a $2$-dimensional  
quantum field theory is equivalent to a
Frobenius algebra $A$.  That is, $A$ is an algebra over a field $k$,   
equipped with a ``trace map"  $\theta :A \to k$,
such that the pairing
$$
A \otimes A \xr{multiply} A \xr{\theta} k
$$
is nonsingular.   A well known example of a Frobenius algebra is the  
homology of a connected, closed, oriented manifold, $H_*(M)$, where  
the product is the intersection product, and the trace map is the  
projection onto the $H_0$ summand.  The resulting nondegeneracy is a  
manifestation of Poincare duality.
As we will see, the basic  Frobenius algebra  of $H_*(M)$ is realized  
by our Morse field theory, when the
homology classes $\alpha$ are simply taken to be the generator $ 
\alpha = 1 \in H_0(B(\ag))$.  In other words,
the basic Frobenius algebra structure is the nonequivariant version  
of our field theory, achieved by choosing a fixed metric-Morse  
structure on the graph.  It is interesting that by choosing families  
of these structures  we obtain
operations $$q_\G : H_*(B(\ag)) \otimes H_*(M)^{\otimes p} \to H_*(M)^ 
{\otimes q}$$ that satisfy the appropriate
gluing and invariance properties.  Thus we get an extended Frobenius  
algebra structure on $H_*(M)$,  whose operations we prove
encompass such classical operations in algebraic topology as Steenrod  
squares and  Stiefel-Whitney classes.
This structure is analogous to the structure in $2$-dimensional field  
theory, where given a connected genus $g$-cobordism between $p$  
circles and $q$ circles, one has an operation,
$$\mu : H_*(\cm_{g, p+q}) \otimes Z(S^1)^{\otimes p} \to Z(S^1)^ 
{\otimes q}
$$
which satisfy gluing laws.  Such a field theory is called a  
homological conformal field theory \cite{manin}.

We will also prove that the field theoretic properties (invariance  
and gluing) of our operations force the classical
relations among cohomology operations such as the Adem relations and  
the Cartan formulae.

\med
The organization of this paper is as follows.  In sections 1 and 2 we  
define the moduli spaces
of metric graph structures, $\cm_\G$, as well as the moduli space of  
graph flows in a manifold, $\cm_\G(M)$.
These are described in algebraic topological terms, using categories  
of graphs, following ideas of Culler-Vogtmann \cite{cullervogtmann},  
Igusa \cite{igusa}, and Godin \cite{godin}.   We then describe a  
generalized Pontrjagin-Thom construction that allows us to define  
fundamental classes of these moduli spaces, without having to study
smoothness or compactness issues.  We then define the invariants (the  
graph operations) in section 3, and prove their field theoretic  
properties in section 4.  In section 5 we describe examples of these  
invariants, and show how cup products,  Steenrod operations,  and  
Stiefel-Whitney classes arise.  We also show how the Cartan and Adem  
formulae follow from the field theoretic properties.

   The second half of the paper begins in section 6, where we deal  
with the geometric and analytic aspects of the moduli spaces, and  
give a more combinatorial, Morse theoretic description of the graph  
operations. Transversality, compactness issues, the resulting  
smoothness of the moduli spaces is studied in sections 6 through  
8.    The Morse theoretic description of  the graph operations is  
given in section 9, where they are shown to live on the level of the  
Morse-Smale chain complexes associated to Morse functions.  In  
particular the operations are defined by suitably counting the number  
of metric graph flows in $M$ that satisfy certain boundary conditions.
A geometric proof of a generalized gluing formula is also given.

There are three appendices to the paper.  Two
cover analytic issues such as regularity and Fredholm properties.   
The third gives a detailed description of the
generalized Pontrjagin-Thom construction needed to define the virtual  
fundamental classes of the moduli spaces.

The first author would like to acknowledge and thank the Department   
of  Mathematics and Statistics at Melbourne University for its  
hospitality during a visit in 2004 where much of the work in this  
paper was carried out.  The second author would like to acknowledge the support of the Australian Academy of Sciences.

\section{Categories of graphs, and the moduli space of metric-Morse  
structures on a graph.}
In this section we describe a category of graphs that will be used to  
define our moduli space of graph flows.  As we will show, the  
geometric realization of this category will
consist of graphs equipped with appropriate metrics.   The idea of  
this category was inspired by the work of Culler-Vogtmann \cite 
{cullervogtmann}, and the interpretation of this work due to Igusa  
\cite{igusa} and Godin \cite{godin}.
\med
\begin{definition}\label{categ} Define $\cc_{b, p+q}$ to be the  
category of oriented graphs of first Betti number $b$, with $p+q$  
leaves.  More specifically, the objects of $\cc_{b, p+q}$ are finite  
graphs (one dimensional CW-complexes) $\G$,  with the following  
properties:
\begin{enumerate}
\item Each edge of the graph $\G$ has an orientation.
\item  $\G$ has  $p+q$ univalent vertices, or ``leaves".   $p$ of  
these are vertices of edges whose orientation points away from the  
vertex
(toward the body of the graph).  These are called  ``incoming"   
leaves.  The remaining $q$ leaves are on edges whose orientation points
toward the vertex (away from the body of the graph).  These are  
called ``outgoing" leaves.
\item $\G$ comes equipped with a  ``basepoint", which is a  
nonunivalent vertex.
\end{enumerate}
For set theoretic reasons we also assume that the objects in this  
category (the graphs) are subspaces  of a fixed infinite dimensional  
Euclidean space, $\br^\infty$.

A morphism between objects $\phi : \G_1 \to \G_2$ is combinatorial  
map of graphs (cellular map) that satisfies:
\begin{enumerate}
\item $\phi$ preserves the orientations of each edge.
\item The inverse image of each vertex is a tree (i.e a contractible  
subgraph).
\item The inverse image of each open edge is an open edge.
\item $\phi$ preserves the basepoints.
\end{enumerate}
\end{definition}

\med
\begin{figure}[ht]
   \centering
   \includegraphics[height=4cm]{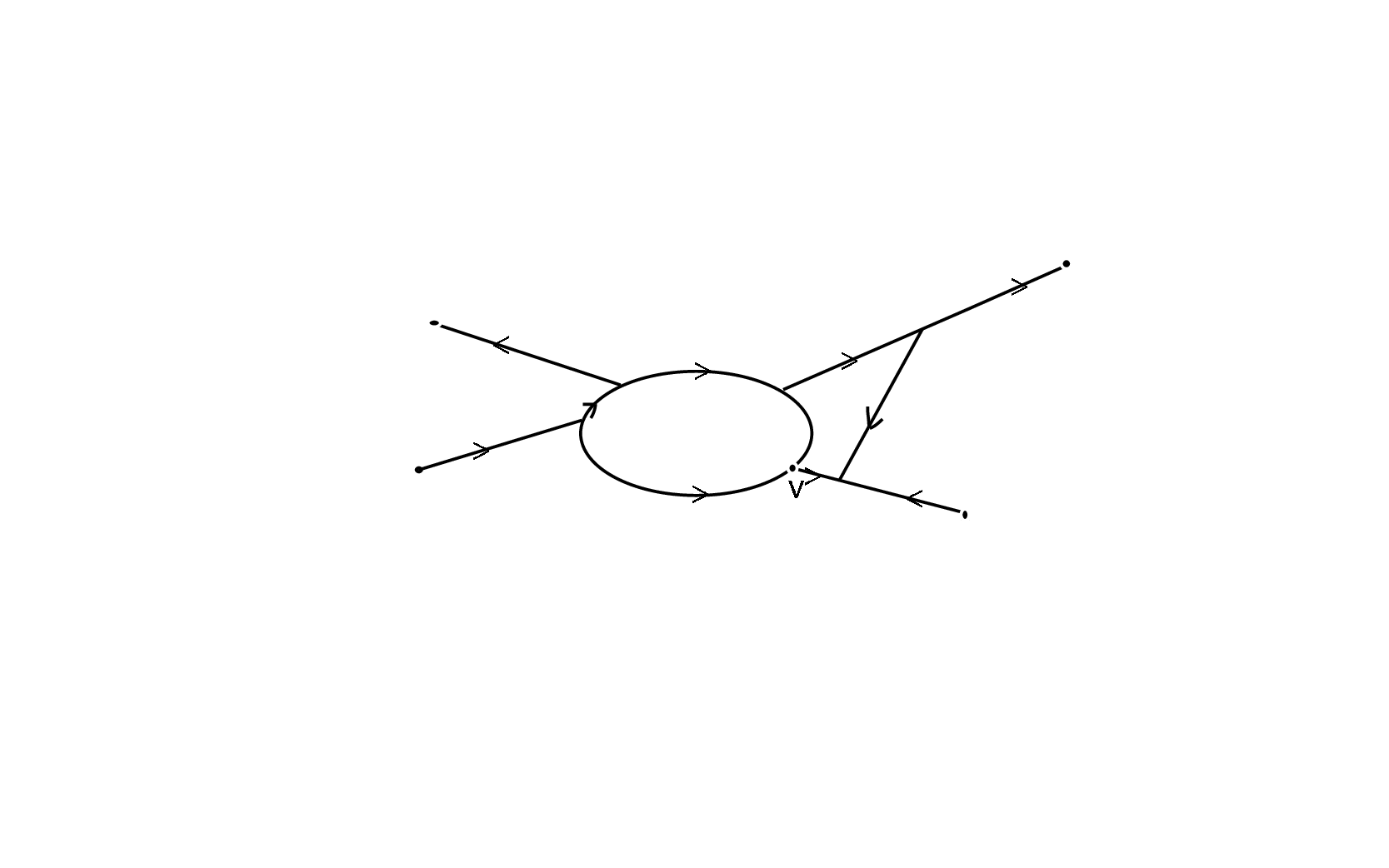}
   \caption{An object $\Gamma$ in $\cc_{2, 2+2}$}
   \label{fig:figone}
\end{figure}

We observe that by the  definition of $\cc_{b, p+q}$,  a   morphism    
$\phi : \G_1 \to \G_2$ is a basepoint preserving cellular map which  
is a homotopy equivalence.
Given a graph $\Gamma \in \cc_{b, p+q}$, we define the automorphism  
group $Aut (\G)$ to be the group of invertible morphisms from
$\Gamma$ to itself in this category.    $Aut (\G)$ is a finite group,  
as it is a subgroup of the group of permutations of the the edges.

\med
We now fix a graph $\G$ (an object in $ \cc_{b, p+q}$), and we  
describe the category of  ``graphs over $\G$", $\ccg$.  As we will  
see below, a point in the geometric realization of this category will  
be viewed as a metric on a generalized subdivision of $\G$.

\med
\begin{definition} Define $\ccg$ to be the category whose  objects  
are morphisms in $\cc_{b, p+q}$ with target $\G$:  $\phi : \G_0 \to \G$.
A morphism from $\phi_0 : \G_0 \to \G$ to $\phi_1 : \G_1 \to \G$ is a  
morphism $\psi : \G_0 \to \G_1$ in $ \cc_{b, p+q}$  with the property
that $ \phi_0  = \phi_1 \circ \psi  : \Gamma_0 \to \Gamma_1 \to \G$.
\end{definition}

\med
Notice that the identity map $id : \G \to \G$ is a terminal object in  
$\ccg$.  That is, every object $\phi : \G_0 \to \G$ has a unique
morphism to $id : \G \to \G$.  This implies that the geometric  
realization of the category, $|\ccg|$ is contractible.     But notice  
that the category $\ccg$ has a free right action of the automorphism  
group, $Aut (\G)$,  given on the objects by composition:
\begin{align}
Objects \, (\ccg) \times Aut (\G) & \to Objects \, (\ccg) \notag \\
(\phi : \G_0 \to \G ) \cdot g &\to g \circ \phi : \G_0 \xrightarrow 
{\phi} \G \xrightarrow{g} \G
\end{align}
This induces a free action on the geometric realization $\ccg$.  We  
therefore have the following:

\begin{proposition}\label{baut}
The orbit space is homotopy equivalent to the classifying space,
$$|\ccg|/ \ag  \simeq B\ag. $$
\end{proposition}

\med
We now consider the geometric realization of the category $|\ccg| 
$.     Following an idea of Culler-Vogtmann \cite{cullervogtmann} and  
Igusa \cite{igusa}, we interpret a point in this space as defining a  
metric on a generalized subdivision of the graph $\G$.

Recall that
$$
|\ccg| = \bigcup_k \Delta^k \times
   \{\G_k \xr{\psi_k} \G_{k-1}
   \xr{\psi_{k-1}}\G_{k-2} \to \cdots \xr{\psi_1} \Gamma_0 \xr{\phi} 
\G \} / \sim
$$
where the identifications  come from the face and degeneracy operations.

\begin{figure}[ht]
   \centering
   \includegraphics[height=12cm]{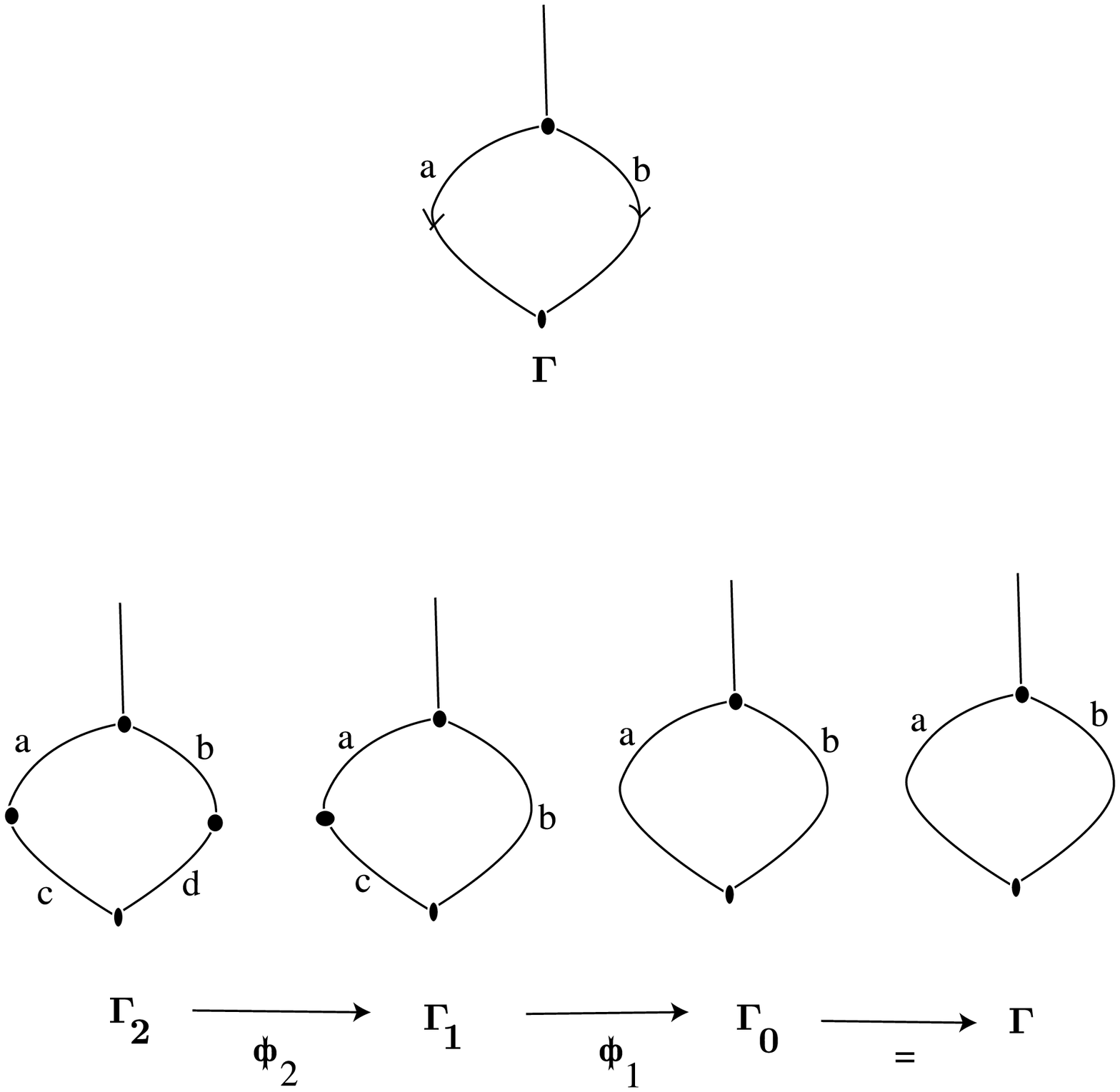}
   \caption{A $2$-simplex in $|\ccg|$.}
   \label{fig:figtwo}
\end{figure}

Let $(\vec{t}, \vec{\psi})$ be a point in $|\ccg|$, where
$\vec{t} = (t_0, t_1, \cdots , t_k)$ is a vector of positive numbers  
whose sum equals one, and  $\vec{\psi}$ is a sequence of $k$-composable
morphisms in $\ccg$.   Recall that a morphism $\phi_i : \G_i \to \G_ 
{i-1}$ can only collapse trees, or perhaps compose such a collapse  
with an automorphism.  So given a composition of morphisms,
$$
\vec{\psi} : \G_k \to \cdots \to \G_0 \to \G
$$
we may think of $\G_k$ is a  generalized   subdivision of $\G$, in  
the sense that $\G$ is obtained from $\G_k$ by collapsing various edges.

We use the coordinates $\vec{t}$ of the simplex $\Delta^k$ to define  
a metric on $\G_k$ as follows.  For each edge $E $ of $\G_k$, define  
$k+1$ numbers, $\lambda_0(E), \cdots , \lambda_k (E)$ given by
$$
\lambda_i(E) = \begin{cases} 0 \quad \text{if $E$ is collapsed by $ 
\vec{\psi}$ in $\G_i$, and}, \\
1 \quad \text{if $E$ is not collapsed in $\G_i$} \end{cases}
$$

We then define the length of the edge $E$ to be
\begin{equation}\label{length}
\ell (E) = \sum_{i=0}^k t_i \lambda_i (E).
\end{equation}

\begin{figure}[ht]
   \centering
   \includegraphics[height=12cm]{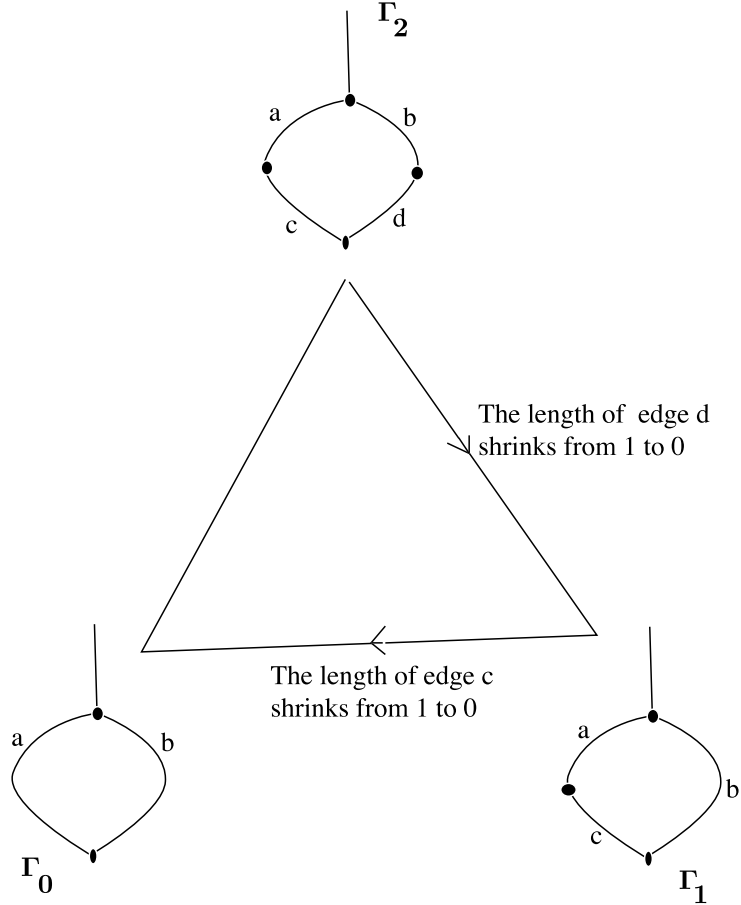}
   \caption{A $2$-simplex of metrics.}
   \label{fig:figthree}
\end{figure}

Notice also that
the orientation
on the edges and the metric deterimine parameterizations (isometries)  
of standard intervals to the edges of the graph $\G_k$ over $\G$,
\begin{equation}\label{param}
\theta_E : [0, \ell (E)] \xr{\cong} E
\end{equation}

Thus a point $ (\vec{t}, \vec{\psi}) \in|\ccg|$  determines a metric  
on a graph  $\G_k$ living over $\G$, as well as a parameterization of  
its edges.   In some sense this may be viewed as the analogue in our  
theory, of the moduli space of Riemann surfaces in Gromov-Witten  
theory.  In that theory, one
studies maps from a Riemann surface (an element of moduli space) to a  
symplectic manifold, which
satisfy the Cauchy-Riemann equations (or some perturbation of them)   
with respect to a choice of a compatible  almost complex structure on  
the symplectic manifold.  In our case, we want to study maps from an  
element of our moduli space,
i.e a graph living over $\G$, equipped with a metric and  
parameterization of the edges, to a target
manifold $M$,  that satisfies certain \sl ordinary \rm differential  
equations. These differential equations will be the gradient flow  
equations of smooth functions on $M$.  To define these, we need to  
impose more structure on our graphs, given by a labeling of the edges  
of the graph by distinct smooth functions on the manifold.  We call  
such a structure a \sl Morse labeling \rm of a graph.  We define this  
precisely as follows.

\med
Let $V$ be a real  vector space.  Let $F(V, k)$ be the configuration  
space of $k$ distinct ordered points in $V$.
That is, $  F(V, k) = \{(v_1, \cdots , v_k) \in V^k \, \text{such  
that} \, v_i \neq v_j \, \text{if} \, i \neq j \}.$
Recall that if $V$ is infinite dimensional, $F(V, k)$ is contractible.

\med
Throughout the rest of this section we let $M$ be a fixed closed,  
Riemannian manifold.

\med
\begin{definition}  An $M$-\sl Morse labeling \rm  of a graph $\G \in  
\cc_{b, p+q}$ is a pair
$(\phi_0:\G_0 \to \G, \,  c)$, where $\phi_0:\G_0 \to \G$ is an  
object of $\ccg$, and $c \in F(C^\infty (M), e(\G_0))$, where $C^ 
\infty (M)$ is the vector space of smooth, real valued functions on $M 
$, and $e(\G_0)$ is the number of edges of $\G_0$.   We think of the  
vector of functions making up the configuration $c$ as labeling the  
edges of $\G_0$.
\end{definition}

\med
Fixing our manifold $M$  and graph $\G$,  our goal now is to define  
the  moduli space of metrics and Morse structures (abbreviated  
``structures")  on $\G$, $  \cm_\G  $.  We do this as follows.

Consider the functor
$$
\mu : \ccg \to \, Spaces
$$
which assigns to a graph over $\G$,    $\phi_0 : \G_0 \to \G$, the  
configuration space $F(C^\infty (M), e(\G_0))$.    Given a morphism
$\psi : \G_1 \to \G_0$, which collapses certain edges  and perhaps  
permutes others, there is an obvious induced map,
$$
\mu (\psi) : F(C^\infty (M),  e(\G_1)) \to F(C^\infty (M), e(\G_0)).
$$
This map projects off of the coordinates corresponding to edges  
collapsed by $\psi$, and permutes coordinates corresponding to  the  
permutation of edges induced by $\psi$.

We can now do a homotopy theoretic construction, called the homotopy  
colimit (see for example \cite{bousfieldkan}).
\begin{definition}\label{structure} We define the space of metric  
structures and Morse labelings on $G$, $ \cs_\G $,  to be the  
homotopy colimit,
$$  \cs_\G = hocolim \, (\mu : \ccg \to \, Spaces).$$
\end{definition}

The homotopy colimit construction is a simplicial space whose $k$  
simplices consist of pairs,
$(\vec{f}, \vec{\psi})$, where $\vec{\psi} : \Gamma_k \to \Gamma_ 
{k-1} \to \cdots \Gamma_0 \to \Gamma$ is a $k$-tuple of composable  
morphisms in $\ccg$, and $\vec{f} \in \mu (\Gamma_k)$.  That is, $\vec 
{f}$ is an  $M$ - Morse labeling  of the edges of $\G_k$.
So we can think of a point $\sigma \in \cs_\G $ as  defining a metric  
on a graph over $\G$, together with an $M$-  Morse labeling of its  
edges.

We now make the following observation.

\begin{lemma}  The space of metric-Morse structures  $\cs_\G $ is  
contractible with a free $\ag$ action.
\end{lemma}

\begin{proof}   The contractibility follows from  standard facts  
about the homotopy colimit construction,  considering the fact that  
both $|\ccg|$ and $F(C^\infty (M), m)$ are contractible.  The free  
action of $\ag$ on $|\ccg|$ extends to an action on $\cs_\G $, since $ 
\ag$ acts by permuting the edges of $\G$, and therefore permutes the  
labels accordingly.  \end{proof}

\med
We now define our moduli space of structures.

\med
\begin{definition}  The moduli space of metric structures and $M$- 
Morse labelings on $G$, $ \cm_\G$,  is defined to be the quotient,
$$
\cm_\G  =\cs_\G / \ag.
$$
\end{definition}

\med
We therefore have the following.

\begin{corollary} The moduli space is a classifying space of the  
automorphism group,
$$
\cm_\G \simeq B\ag.
$$
\end{corollary}

\section{The moduli space of  metric-graph   flows in a manifold}

Let $M$ be a fixed, smooth, closed $n$-manifold with a   Riemannian  
metric.  Let $\G \in \cc_{b, p+q}$ be a graph.  In this section
we define  the moduli space of $\G$-flows in $M$, $  \cm_\G (M)$,   
and study its topology.  This will be an infinite
dimensional space built from the moduli space of metric-Morse  
structures, $ \cm_\G$, which in turn has an infinite dimensional
homotopy type, since $ \cm_\G \simeq B\ag$, and $\ag$ is a finite  
group.  However, given a homology class $\alpha \in H_k(\ag)$, we  
show how to define a ``virtual fundamental class",
$$
[ \cm^\alpha_\G(M)] \in H_q( \cm_\G(M))
$$
where $q = k +\chi (\G)n$, where $\chi (\G)  = 1-b$ is the Euler  
characteristic.   The smooth structures on these moduli spaces will  
be studied in later sections.  But even without knowledge of this  
structure,    these virtual fundamental classes  will be  constructed  
using generalized Pontrjagin-Thom constructions  similar to those   
defined in \cite{cohenumkehr}.    These constructions allow us to  
define invariants in the next section, which we will identify with  
classical  cohomology operations
in   section 4.
\med
Let $\sigma \in \cs_\G$ be a metric-Morse structure.  Then $\sigma =  
( \vec{t}, \vec{\psi}, c)$, where $ \vec{t}\in \Delta^k$,   $\vec 
{\psi} : \G_k \to \cdots \to \G_0 \to \G$ is a $k$-simplex in the  
nerve of $\ccg$, that is a $k$-tuple of composable morphisms, and $c$  
is a Morse labeling of the edges of $\G_k$.
\begin{definition}  A metric-$\G$-flow in $M$, is a pair $(\sigma,  
\gamma)$, where $\sigma = ( \vec{t}, \vec{\psi}, c) \in   \cs_\G$ is  
a metric-Morse structure on $\G$,
and $\gamma : \G_k \to M$ is a continuous map, smooth on the edges,  
satisfying the following property.    Given any edge $E$ of $\G_k, $   
let    $\gamma_E : [0, \ell (E)] \to M $ be the composition
$$
\gamma_E : [0, \ell (E)] \xr{\theta_E} E \subset \G_k \xr{\gamma} M,
$$   where $\theta_E$ is the parameterization of the edge $E$   
defined in (\ref{param}).   Then $\gamma_E$ is required to  satisfy  
the differential equation
$$
\frac{d\gamma_E}{dt} (s) + \nabla f_E(\gamma_E(s)) = 0.
$$
Here the collection of  labeling functions $\{f_E: M \to \br  \, :  
\,  E \, \text{is an edge of} \,    \, \G \} $ is the configuration   
$c \in F(C^\infty (M), e(\G))$   determined by the structure $\sigma$.
\end{definition}
We define the ``structure space of metric-graph flows", $\tccgm$, to  
be the space
\begin{equation}\label{graphflow}
\tccgm = \{(\sigma, \gamma) \, \text{a  metric-$\G$-flow in $M$}\},
\end{equation}
and the moduli space of graph flows to be the orbit space,
$$
\ccgm = \tccgm/\ag.
$$ Here the automorphism group  $\ag$ acts on $\tccgm$ by acting on  
the structure $\sigma$ as described above.

\med
We have not yet defined the topology on these spaces of flows.  To do  
that we first consider  the case when the graph $\G$ is a tree.  That  
is, $\G$ is contractible, so $b_1(\G) = 0$.

\med
\begin{proposition}\label{tree} Let $\G$ be a tree. Then there is  an  
$\ag$-equivariant bijective correspondence
\begin{align}
\Psi : \tccgm &\xr{\cong} \cs_\G \times M \notag \\
(\sigma, \gamma) &\to  \sigma \times \gamma (v) \notag
\end{align}
where $v$ is the fixed vertex of the graph $\G_k$ over $\G$  
determined by the structure $\sigma$.  On the right hand side, $\ag$  
acts on $\cs_\G $ as described above, and acts trivially on $M$.
\end{proposition}

\begin{proof}  This follows from the existence and uniqueness theorem  
for solutions of ODE's on compact
manifolds.  The point is that the values of $\gamma$ on the edges  
emanating from $v$ are completely determined
by $\gamma (v) \in M$, since one has a unique flow line through that  
point for any of the functions labeling these edges.
The value of $\gamma$ on these edges determines the value of $\gamma$  
on coincident edges (i.e edges that share a vertex)
for the same reason.  The  fact that $\Psi$ is a bijection now  
follows.  The $\ag$-equivariance of $\Psi$ is immediate.
\end{proof}

\med
We now topologize $\tccgm$ so that $\Psi : \tccgm \to \cs (\G) \times  
M$ is a homeomorphism.  We then have the following description of the  
moduli space of
graph flows, when $\G$ is a tree:

\med
\begin{corollary}\label{treemoduli}  Let $\G$ be a tree. Then $\Psi$  
induces   a homeomorphism,
$$
\Psi : \ccgm \xr{\cong} \cs_\G/\ag  \times M
$$
which has the homotopy type of $B\ag  \times M$.
\end{corollary}

For general connected graphs $\G$, we analyze the topology of $\ccgm$  
in the following way.
Let $\sigma \in \cs_\G$.  A \sl tree  flow \rm of  $\Gamma$ with  
respect to the structure $\sigma$ is a collection $\gamma = \{\gamma_T 
\}$ where $\gamma_T : T \to M$ is a graph flow on a maximal subtree  
$T \subset \G_k$.  The collection ranges over all maximal subtrees $T  
\subset \G_k$, and is subject only to the condition that
the values at the basepoint are the same:
$$
\gamma_{T_1} (v) = \gamma_{T_2}(v)
$$
for any two maximal trees $T_1, T_2 \subset \G_k$.  (Here $v \in T  
\subset \G$ is the fixed point vertex.)

We define
\begin{equation}
\tilde \cm_{tree}(\G, M) = \{(\sigma, \gamma) \, : \,  \sigma \in \cs_ 
\G, \, \text{and} \,
\gamma = \{\gamma_T\} \, \text{is a tree flow of } \, \Gamma \, \text 
{with respect to} \, \sigma \}
\end{equation}
and $$
\cm_{tree}(\G, M) = \tilde\cm_{tree}(\G,M)/\ag.$$

Notice that the proof of proposition \ref {tree} also proves the  
following.

\begin{theorem}     For any graph  $\G \in \cc_{b, p+q}$ there is an $ 
\ag$ -equivariant bijective correspondence,
\begin{align}
\Psi : \tilde \cm_{tree}(\G, M) &\xr{\cong} \cs_\G    \times M   
\notag \\
(\sigma, \gamma)  &\la   \sigma \times \gamma (v). \notag
\end{align}
\end{theorem}

We therefore again topologize $\tilde \cm_{tree}(\G,M)$ so that $\Psi 
$ is an equivariant homeomorphism.  Then $$\cm_{tree}(\G,M) \cong    
\cs_\G/\ag \times M \simeq B\ag \times M. $$

\med
Consider the inclusion,   $\tilde \rho : \tccgm \hk \tilde \cm_{tree} 
(\G, M)$ defined to be the map that sends a graph flow $\gamma$ to  
the tree flow obtained by restricting $\gamma$ to each maximal tree.   
We then give $\tccgm$ the subspace topology, which makes $\rho$ an   
equivariant embedding.  This defines an embedding
$\rho : \ccgm \hk \cm_{tree}(\G, M)$.

\med
We use this embedding to define virtual fundamental classes of $\ccgm 
$. Recall  that the space $\ccgm$ is infinite dimensional
because the moduli space $\cs_\G/\ag  \simeq B\ag$ is infinite  
dimensional.  We can ``cut down" this moduli space by considering
an embedding of a  compact manifold of structures, $  \tilde N  
\subset \cs_\G$. We let
$N = \tilde N/\ag \subset \cs_\G/\ag \simeq B\ag$.   We can then  
define the space
$\ccngm \subset \ccgm$ to be the subspace $\ccngm = \{  (\sigma,  
\gamma)  \in     \tccgm$ such that $\sigma \in \tilde N)/\ag \}$.   
Then the embedding $\rho : \ccgm \hk \cm_{tree}(\G, M) \cong  \cs_\G/ 
\ag \times M$ defines an embedding
\begin{equation}\label{subman}
\rho_N : \ccngm \hk N \times M.
\end{equation}

To motivate our  construction  of the virtual fundamental classes,   
suppose we know
   that $\ccngm$ is  a smooth closed submanifold of $N \times M$  of  
codimension $k$.  Then the image of its fundamental class $[\ccngm]  
\in H_*(\ccngm)$  in $H_*(\ccgm)$ would
be the image under the ``umkehr map",
$$
H_*(N \times M) \xr{(\rho_N)_! } H_{*-k}(\ccngm)) \to H_{*-k}(\ccgm)
$$
of the product of the  fundamental classes $[N] \times [M]$.  The  
umkehr map $(\rho_N)_! :  H_*(N \times M) \to H_{*-k}(\ccngm) $ is  
Poincare dual
to the restriction map in cohomology, $\rho_N^* : H^*(N \times M) \to  
H^*(\ccngm), $  induced by the embedding $\rho_N : \ccngm \hk N  
\times M$.  In particular the fundamental class $[\ccngm] \in H_{*-k} 
(\ccgm)$ only depends on the homology class represented by the  
manifold $[N] \in H_*(\cs_\G/\ag) \cong H_*(B\ag)$.

To define our ``virtual fundamental class", we avoid the question of  
whether $\ccngm$ can be given a smooth structure (we address this  
question in a later section), by directly defining  the umkehr map
\begin{equation}\label{shriek}
\rho_! : H_*(B\ag \times M) = H_*(\cs_\G/\ag \times M) \to H_{*-bn} 
(\ccgm)
\end{equation}
where $b = b_1(\G)$, and $n = dim \, M$.   Once we have this map,  
then given $\alpha \in H_q(B\ag)$, the  virtual fundamental class $ 
[\ccalgm]$ is defined by
\begin{equation}\label{virtual}
[\ccalgm] = \rho_!(\alpha \times [M]) \in H_{q-(b-1)n}(\ccgm).
\end{equation}

\med
The rest of this section will be devoted to defining the umkehr map $  
\rho_!$.
The existence of this map follows from a construction that is used to  
give a proof of  a   general existence theorem for umkehr maps by the  
first author and J. Klein in \cite{cohenumkehr}.
This construction is based on the existence of   ``Pontrjagin-Thom  
collapse maps".
We recall that given a smooth embedding of compact manifolds, $e : N  
\hk M$ of codimension $k$,  the umkehr map $e_! : H_*(N) \to H_{*-k}* 
(M)$ can be computed via the Pontrjagin-Thom collapse map,
$$
\tau_e : M \to M/M-\eta_e
$$
where $N \subset \eta_e$ is a tubular neighborhood.  This quotient  
space is the one point compactification of the tubular neighborhood,  
which is homeomorphic to the Thom space of the normal bundle, $N^ 
{\nu_e}$.    The umkehr map is then given by the composition,
$$
\begin{CD}
e_! : H_*(M) \xr{(\tau_e)_*}H_*(N^{\nu_e}) @>\cap u >\cong > H_{*-k}(N)
\end{CD}$$ where the last map is the cap product with the Thom class,  
yielding the Thom isomorphism.

\med

To apply this construction in our setting, we need to produce an open  
neighborhood
$\eta_\eps$ of the embedding  $\rho : \ccgm \hk \cm_{tree}(\G, M)  
\cong \cs_\G/\ag  \times M,$  that is homeomorphic to the total space  
of an appropriate normal bundle,
$\nu_\rho$.  We now define these objects.

Let $T \subset \G$ be a maximal tree.  We define a map
$p_T : \tilde \cm_{tree}(\G, M) \to M^{2b}$ as follows.
Since $T $ is a maximal tree, the complement $\G - T$ consists of $b  
= b_1(\G)$ open edges, $e^T_1, \cdots,  e^T_b$.  Now let    $\phi :  
\Gamma_0 \to \G$ be an object  in $\cc_\G$.  Since the inverse image  
under $\phi $ of an   edge is an   edge, then $\phi^{-1}(e^T_i) =  
e^T_i(\G_0)$ is an edge, and the tree $T(\G_0) = \phi^{-1}(T) \subset  
\G_0$
has complement $\G_0 - T(\G_0)$ given by the $b$ open edges $  e^T_i 
(\G_0), i = 1, \cdots , b$.
The edges $ e^T_i(\G_0)$ are oriented, so they have source and target  
vertices,
$s_i^T(\G_0)$, and $t_i^T(\G_0)$.

Now let  $(\sigma, \gamma)$ be a point in  $\tilde \cm_{tree}(\G, M) 
$.  So    $ \sigma = (\vec{t}, \vec \psi, c) \in \cs_\G$, and  $  
\gamma = \{\gamma_{T_j} : T_j \to M\}$, where
the $T_j$'s are the maximal trees in $\G_k$, and $\gamma_{T_j}$  is
   a graph flow on the tree  $  T_j$ with respect to the structure $ 
\sigma$.    Let $T_1 = T(\G_k) = \phi_k^{-1}(T) \subset \G_k$.

Consider the graph flow $\gamma_{T_1} : T_1 \to M$, and let $x_i $ be  
the image of the source vertex,
\begin{equation}\label{xeye}
x_i = \gamma_{T_1}(s^T_i(\G_k)) \in M.
\end{equation}
Now consider the  image of the target vertex,  $\gamma_{T_1}(t^T_i 
(\G_k)) \in M.$
The existence and uniqueness theorem for solutions of ODE's says  
there is a unique map $\alpha_i : e^T_i(\G_k)\to M$ which is graph  
flow with respect to the structure $\sigma$,  satisfying the initial  
condition,   $\alpha_i (t^T_i(\G_k)   =  \gamma_{T_1}(t^T_i(\G_k))  
\in M. $ We then define $y_i \in M$ to be the image of the source  
vertex under the map $\alpha_i$:
$$
y_i = \alpha_i(s^T_i(\G_k)) \in M.
$$
Notice that   the tree flow $\gamma$  is induced from a flow on the  
full graph $\G$ if and only if $x_i = y_i$ for all $i = 1, \cdots ,b 
$.  Said another way, we have defined a map
\begin{align}\label{peetee}
p_T : \tilde \cm_{tree}(\G,M) &\to (M^2)^b \\
(\sigma, \gamma) &\to (x_1, y_1), \cdots, (x_b, y_b) \notag
\end{align}
where  the following diagram is a pullback square:
\begin{equation}\label{pullback}
\begin{CD}
\tccgm  @>\rho > \hk >  \tilde \cm_{tree}(\G,M)  \\
@Vp_TVV    @VVp_T V \\
M^b @>\hk>\Delta^b> (M^2)^b.
\end{CD}
\end{equation}
Here $\Delta : M \to M^2$ is the diagonal.
\med
We now define our tubular neighborhood and normal bundle.  Give $M$ a  
Riemannian metric.
\begin{definition} 1.  For $\epsilon > 0$,  let $\eta_\eps \subset  
\tilde \cm_{tree}(\G,M)$
be the open set containing $\rho (\tccgm)$ defined  to be the inverse  
image of the $\eps$-neighborhood of the diagonal,
$$
\eta_\eps = \{(\sigma, \gamma) \in \tccgm \, : d(p_T(\sigma, \gamma),  
\Delta (M)) < \eps \quad \text{for every maximal tree $T \subset \G$} \}
$$ where $d $ is the Riemannian distance in $M\times M$.
2.  Let $\nu(\rho) \to \tccgm$ be the vector bundle defined as follows.
Let $p : \tccgm \to M$ be the map $(\sigma, \gamma) \to \gamma (v)$.   
This is the right hand factor of the embedding $\rho : \tccgm \hk  
\tilde \cm_{tree}(\G,M) \cong   \cs_\G \times M$.  Define
$$
\nu(\rho) = p^*(\bigoplus_b TM)
$$
to be  the pullback of the Whitney sum of $b$-copies of the tangent  
bundle.
\end{definition}
We notice that $\eta_e$ is an $\ag$-invariant open subspace of $  
\tilde \cm_{tree}(\G,M)$
and therefore defines an open neigborhood which, by abuse of notation  
we also call $\eta_\eps$ of the embedding of quotient spaces, $\rho :  
\ccgm \hk  \cm_{tree}(\G,M)$.
Similarly, $\nu(\rho)$ is an invariant bundle over $\tccgm$,and  
therefore defines
a bundle $\nu (\rho) = p^*(\bigoplus_b TM)$ over $\ccgm$.
   The following theorem will allow us to  define a Pontrjagin-Thom  
collapse map, which as observed above, will allow us to define the  
umkehr map $\rho_!$. This is a tubular neighborhood theorem for the  
embedding $\rho : \ccgm \hk  \cm_{tree}(\G,M).$  Its proof is rather  
technical, so we leave it to the appendix.

   \med
   \begin{theorem}\label{tube} For $\eps > 0$ sufficiently small,  
there is a homeomorphism $\Theta : \eta_\eps \xr{\cong} \nu (\rho)$  
taking $\ccgm$ to the zero section.
   \end{theorem}

   \med
     The homeomorphism $\Theta$ then defines a homeomorphism of the  
quotient space to the Thom space,
     $$
     \Theta :  \cm_{tree}(\G,M)/ (\cm_{tree}(\G,M) - \eta_\eps )  \,  
\la  \ccgm^{\nu (\rho)}
     $$
     and so we have a Pontrjagin-Thom collapse map,
     \begin{equation}\label{pontthom}
     \tau_\rho : \cs_\G/\ag \times M  \cong \cm_{tree}(\G,M) \xr 
{project} \cm_{tree}(\G,M)/ (\cm_{tree}(\G,M) - \eta_\eps ) \xr 
{\Theta} \ccgm^{\nu (\rho)}.
     \end{equation}

     Assuming $M$ is oriented, this defines an umkehr map,
     \begin{align}\label{umkehr}
     \rho_! :H_*(B\ag \times M) \cong H_*(\cs_\G/\ag \times M )   &  
\xr{\tau_\rho}   H_*(\ccgm^{\nu (\rho)}) \\
     &\xr{Thom \, iso} H_{*-b\cdot n}(\ccgm).  \notag
     \end{align}

\med
We are now ready to define virtual fundamental classes of these  
moduli spaces.

     \begin{definition}
     Let  $\alpha \in H_q(B\ag ; k)$, where $k$ is a coefficient field.
Define the virtual fundamental class, $[\cm^\alpha (\G,M)] \in H_{q +  
\chi (\G) n}(\ccgm; k)$ to be the image of $\alpha \otimes [M]$ under  
the umkehr map
$$
   \rho_! : H_*(B\ag ;k ) \otimes H_*(M; k) \to H_{*-bn}(\ccgm; k).
   $$
   \end{definition}

   Notice that since $1-b$ is the Euler characteristic  $\chi (\G)$,   
we have that the virtual fundamental class associated to a homology  
class $\alpha$ of degree $q$ lies
   in degree, $q+\chi (\G)n$,
   $$
   [\cm^\alpha (\G,M)] \in H_{q +\chi(\G)n}(\ccgm; k).
   $$

   These virtual fundamental classes, and more generally the umkehr  
map $\rho_!$,  will allow us to define   cohomology operations  
yielding the Morse Field Theory described in the introduction.  We  
define and study these operations in the next section.

\section{Graph operations}
In this section we describe Gromov-Witten type operations induced by  
our moduli spaces of graphs and their virtual vector bundles.  We  
actually describe two types of operations induced by a graph
$\G$, the first, $q^0_\G$, is equivariant with respect to the \sl  
bordered \rm automorphism group, $\aog$, which consists of those  
automorphisms $g \in \ag$ that fix the marked points (univalent  
vertices).  These operations are the directly analogous to Gromov- 
Witten operations.  We then show how these
operations can be extended to operations $q_\G$ that are equivariant  
with respect to the full automorphism group.

\med
Let $\G$ be an object in $\cc_{b, p+q}$, and $M$ a closed, $n$- 
dimensional manifold. In what follows we consider homology and  
cohomology with coefficients in an arbitrary but fixed field $k$.
We begin by defining the operations,

\begin{equation}
q^0_\G : H_*(B\aog) \otimes H_*(M)^{\otimes p} \la H_*(M)^{\otimes q}
\end{equation}
which raises total dimension by $\chi (\G)n -np$ where $\chi(\G)$ is  
the Euler characteristic of the graph $\G$  ($\chi (\G) = 1-b$),  $b$  
is the first Betti number of $\G$, and $p$ and $q$ are the number of  
incoming and outgoing marked points of $\G$ respectively.

Let $\cm_0(\G, M) = \tccgm/\aog \simeq B\aog$.  Consider the  
evaluation maps $ev_{in} : \tccgm \to M^p$ and $ev_{out} : \tccgm \to  
M^q$   that evaluate  a graph flow on the incoming and outgoing  
marked points, respectively. Since automorphisms in $\aog$ preserve  
these marked points, they descend to give maps $ev_{in} : \cm_0(\G,M)  
\to M^p$ and $ev_{out}: \cm_0(\G,M) \to M^q$.  Let $ev$ be the  
product map,
$$
ev = ev_{in} \times ev_{out} : \cm_0(\G, M) \la M^p \times M^q.
$$

Let $\alpha \in H_r(B\aog) = H_r(\cm_0(\G, M))$. As we did in the  
last section (\ref{virtual}) we can define a virtual fundamental class
$$
[\cm_0^\alpha(\G, M)] = \rho_!(\alpha \times [M]) \in H_{r+n-bn}(\cm_0 
(\G, M)) = H_{r+\chi(\G)n}(\cm_0(\G, M)).
$$
Consider the Gromov-Witten type invariant,
\begin{align}
\bar q^0_\G (\alpha) : H^*(M)^{\otimes p} \otimes H^*(M)^{\otimes q} & 
\to k  \notag \\
x \otimes y &\to \langle ev^*(x \otimes y), [\cm_0^\alpha(\G, M)]  
\rangle. \notag
\end{align}
Notice that $\bar q^0_\G (\alpha)$ can only be nonzero if the total  
dimension of $x \otimes y$ is $r+\chi(\G)n$.
We may  think of $\bar q^0_\G(\alpha) $ as an element of homology,
$$
\bar q^0_\G(\alpha) = ev_*( [\cm_0^\alpha(\G, M)]) \in H_*(M^{ p})  
\otimes H_*(M^{ q})
$$
of total dimension $r+\chi(\G)n$.  By applying Poincare duality to  
the left hand tensor factor, this defines a class
$$
   q^0_\G(\alpha) \in  H^{np -*}(M^p) \otimes H_*(M^q) \cong Hom (H_* 
(M)^{\otimes p}; H_*(M)^{\otimes q})
   $$
which raises total dimension by $r+\chi(\G)n-np$.  We have therefore  
defined an operation
\begin{equation}\label{q0}
q^0_\G : H_*(B\aog) \otimes H_*(M)^{\otimes p} \la H_*(M)^{\otimes q}
\end{equation}
which raises total dimension by $\chi (\G)n -np$.

\med
We now describe an extension of the operation $q^0_\G$ to an operator  
on equivariant homology,
$$
q_\G : H_*^{\ag}( M^p) \la H_{*+\chi (\G)n-np}^{\ag} (M^q).
$$
Here $\ag$ acts on $M^p$ via the permutation action determined by the  
homomorphism $\ag \to \Sigma_p$ that sends an automorphism to the  
induced permutation of the $p$-incoming marked points.
The $\ag$ action on $M^q$ is defined similarly.
The sense in which the operation $q_\G$ will extend $q^0_G$, is the  
following. Since an element $g \in \aog$  lies in  the kernel of  
the    homomorphism  $\ag \to \Sigma_p$   its action on   $M^p$ is  
trivial.  Therefore
the inclusion $\aog \subset \ag$ induces a map of homotopy orbit spaces,
$$
B\aog \times M^p \to E\ag \times_{\ag}M^p
$$
and therefore an induced map in homology, $H_*(B\aog) \otimes H_*(M)^ 
{\otimes p} \to H_*^{\ag}(M^p)$.
The compatibility of the operators $q^0_\G$ and $q_\G$ is that the  
following diagram commutes:
\begin{equation}\label{commutes}
\begin{CD}
H_*(B\aog) \otimes H_*(M)^{\otimes p}  & @>q^0_\G >>  & H_*(M^q) \\
@VVV &&    @VVV \\
H_*^{\ag}(M^p)  @>>q_\G> H_*^{\ag}(M^q) @>>> H_*^{\Sigma_q}(M^q).
\end{CD}
\end{equation}
We now define the graph operation $q_\G$.  As above, consider the  
evaluation map
$$
ev_{in} : \tccgm \to M^p,
$$ which evaluates a graph flow on the $p$ incoming marked points.
This map is $\ag$ equivariant, where as above, $\ag$ acts on $M^p$ by  
permuting the coordinates
according to the homomorphism $\ag \to \Sigma_p$.  Taking homotopy  
orbit spaces, we get a map
$$
ev_{in} : \tccgm/\ag = \ccgm \to E\ag \times_{\ag}M^p.
$$
We similarly have a map $ev_{out}: \ccgm \to E\ag  \times_{\ag}M^q.$
Notice that up to homotopy, the map $ev_{in}$ factors as the  
composition,
\begin{equation}\label{compose}
ev_{in} : \ccgm \xr{\rho} \cm_{tree}(\G, M) \cong   \cs_\G (M)/\ag  
\times M \simeq B\ag \times M \xr{\Delta^p} E\ag \times_{\ag} M^p.
\end{equation}
Here $\Delta^p: M \to M^p$ is the $p$-fold diagonal, which  maps $M$  
to the  fixed points of the $\ag$-action on $M^p$.  Therefore by  
applying homotopy orbit spaces, we have an induced map $\Delta^p : B 
\ag \times M \to  E\ag \times_{\ag} M^p$.
$\Delta^p$ is a codimension $n(p-1)$ embedding, and so there is a  
Pontrjagin-Thom map to the Thom space of the normal bundle,
$$
\tau_{\Delta^p} : E\ag \times_{\ag} M^p \to (B\ag \times M)^{\nu  
(\Delta^p)}.
$$
As described in the previous section, such a map induces an umkehr  
map in homology,
$$
(\Delta^p)_! : H_*(E\ag \times_{\ag} M^p) \to H_{*-n(p-1)}(B\ag  
\times M).
$$
Because of the factoring of $ev_{in}$ in (\ref{compose}), we can then  
define the umkehr map $(ev_{in})_!$ as the composition of umkehr maps,
\begin{align}
(ev_{in})_! : H_*(E\ag \times_{\ag} M^p) \xr{(\Delta^p)_!} H_{*-n 
(p-1)}(B\ag \times M) &\xr{\rho_!} H_{*-n(p-1)-bn}(\ccgm) \notag \\
&=H_{*+\chi(\G)n-np}(\ccgm). \notag
\end{align}

We now define the operation $ q_\G$ as follows.

\begin{definition}\label{qgamma}
Define
$$
q_{\G} : H_*^{\ag}( M^p) \la H_{*+\chi (\G)n-np}^{\ag}(M^q)
$$
to be the composition
$$
q_{\G} : H_*(E\ag \times_{\ag} M^p) \xr{(ev_{in})_!} H_{*+\chi(\G)n- 
np}(\ccgm) \xr{ev_{out}}  H_{*+\chi(\G)n-np}(E\ag \times_{\ag} M^q).
$$
\end{definition}

We now observe the following property relating the operations $q_{\G} 
$ and $q^0_{\G}$.

\med
\begin{proposition}   The operation $q_{\G}$ extends $q^0_{\G}$ in  
the sense that it makes diagram (\ref{commutes}) commute.
\end{proposition}

\begin{proof}  Let $\alpha \in H_*(B\aog)$,  $\beta \in H_*(M^p)$,  
and $x \in H^*(M^q)$
be in the image of $H^*_{\Sigma_q}(M^q) \to H^*(M^q).$  Then by  
definition,
\begin{equation}\label{first}
\langle x \, , \, q^0_\G (\alpha \otimes \beta)\rangle = \langle ev^*_ 
{out}(x)\cup ev^*_{in}(D\beta)\, , \, \rho_!(\alpha \otimes [M])\rangle,
\end{equation}
where $D : H_*(M^p) \to H^{np-*}(M^p)$ is Poincare duality.  On the  
other hand,
by the definition of $q_\G$,
$$
\langle x \, , \, q_\G (\alpha \otimes \beta)\rangle = \langle ev^*_ 
{out}(x) \, , \, \rho_!(\alpha \otimes \Delta^p_!(\beta) \rangle.
$$
But by the commutativity of the diagram
$$
\begin{CD}
H_*(M^p) @>\Delta^p_! >> H_*(M) \\
@VD VV   @AA\cap [M] A \\
H^{np-*}(M^p) @>>(\Delta^p)^* > H^{np-*}(M)
\end{CD}
$$
this quantity is equal to
\begin{equation}\label{quantity}
\langle ev^*_{out}(x) \, , \, \rho_!(\alpha \otimes  (\Delta^p)^*(D 
\beta)\cap [M])\rangle.
\end{equation}
Now the dual umkehr map,
$$
\rho^! : H^*(\cm_0(\G, M)) \to H^{*+bn}(B\aog \times M)
$$
is a map of $H^*(  \cm_\G (M))$-modules.   This implies that
$ \rho_!(\alpha \otimes  (\Delta^p)^*(D\beta)\cap [M]) = \rho_! 
(\alpha \otimes [M])\cap \rho^*((\Delta^p)^*(D\beta)).$  Thus   
quantity  (\ref{quantity}) is equal to
\begin{equation}\label{quantity2}
\langle ev^*_{out}(x) \, , \, \rho_!(\alpha \otimes [M])\cap \rho^* 
((\Delta^p)^*(D\beta))   \rangle.
\end{equation}
Now by the definition of the evaluation map $ev_{in}$, the following  
diagram commutes:
$$
\begin{CD}
H^*(M^p) @>(\Delta^p)^* >> H^*(M)   @>\hk >> H^*(M \times B\aog) \\
@V ev_{in} VV  && @VV \rho^*V \\
H^*(\cm_0(\G,M) & @>> = > & H^*(\cm_0(\G,M).
\end{CD}
$$
So quantity (\ref{quantity2}) is equal to
$$
\langle ev^*_{out}(x) \, , \, \rho_!(\alpha \otimes [M])\cap ev^*_ 
{in} (D\beta) = \langle ev^*_{out}(x)\cup ev^*_{in}(D\beta)\, , \,  
\rho_!(\alpha \otimes [M])\rangle,
$$
which is the same as the quantity in equation (\ref{first}).
\end{proof}

\med
We end this section with the observation that given any group  
homomorphism
$\theta:  G \to \ag$, the above constructions and arguments using the  
moduli space
$\cm^G_\G(M) = EG \times_\theta \tilde   \cm_\G (M)$, allow us to    
construct    an umkehr map,
$$
(ev_{in})_! : H^G_*(M^p) \to H_{*+ \chi(\G)n -np}(\cm^G_\G (M)),
$$
which in turn   allows the definition of  an operation defined  on $G 
$-equivariant homology,
\begin{equation}\label{subgroup}
q^{G}_\G : H_*^{G}(M^p) \to H_*^G(M^q)
\end{equation} that is natural with respect to homomorphisms between  
groups living over $\ag$.      That is if $\theta_1 : G_1 \to \ag$  
and $\theta_2 : G_2 \to \ag$ are homomorphisms and  $f: G_1 \to G_2$  
is a group homomorphism such that $\theta_2 \circ f = \theta_1$,  
then  then the following diagram commutes:
\begin{equation}\label{homomorphism}
\begin{CD}
H_*^{G_1}(M^p)  @>q^{G_1}_\G >>   H_*^{G_1}(M^q) \\
@Vf_* VV  @VVf_*V \\
H_*^{G_2}(M^p)  @>>q^{G_2}_\G >   H_*^{G_2}(M^q)
\end{CD}
\end{equation}
\section{Field theoretic properties of the graph operations}
In this section we describe the two properties,  invariance and  
gluing, that imply the assignment to a graph $\G$ the operation $q_\G 
$ defines a field theory.
We begin with the invariance property.  Roughly this says that a  
morphism
$$ \phi : \G_1 \to \G_2$$in $\cc_{b, p+q}$ takes the operation $q_ 
{\G_1}$ to $q_{\G_2}$.   We state this more precisely as follows.
Let $  G_1 <  \ag_1$ and $  G_2  < \ag_2$ be   subgroups.

\begin{definition} We say that a morphism $\phi : \G_1 \to \G_2$ is  
$G_1$-$G_2$ equivariant, if for every $g_1 \in G_1$ there exists a  
unique $g_2 \in G_2$ such that the following composite morphisms are  
equal:
$$
\phi \circ g_1 = g_2 \circ \phi : \G_1 \to \G_2.
$$
\end{definition}
In this setting, $\phi$ determines a homomorphism,
\begin{align}
\phi_* : G_1 &\to G_2 \notag \\
g_1 &\to g_2. \notag
\end{align}
Furthermore, one easily checks that the homomorphism  $\phi_* : \G_1  
\to \G_2$  lives over the identity in the symmetric groups. That is,  
the following diagram commutes:
$$
\begin{CD}
G_1 @>>> \Sigma_p \times \Sigma_q \\
@V\phi_* VV  @VV =  V \\
G_2 @>>> \Sigma_p \times \Sigma_q
\end{CD}
$$
where the two horizontal maps assign to an automorphism the induced  
permutation of the incoming and outgoing leaves.    The commutativity  
of this diagram then says that $\phi$ induces maps of equivariant  
homology,
$$
\phi_* : H_*^{G_1}(M^p) \to H_*^{G_2}(M^p) \quad \text{and} \quad  
\phi_* : H_*^{G_1}(M^q) \to H_*^{G_2}(M^q).
$$

\med
\begin{theorem}(Invariance)\label{invariance}
Let $\phi : \G_1 \to \G_2$ be a morphism in $\cc_{g, p+q}$,   $G_1 <  
\ag_1$ and $G_2 < \ag_2$ be such that $\phi$ is $G_1$-$G_2$  
equivariant.  Then the following diagram commutes:
$$
\begin{CD}
H_*^{G_1}(M^p) @>q^{G_1}_{\G_1} >>  H_*^{G_1}(M^q)  \\
@V\phi_* VV   @VV\phi_* V \\
   H_*^{G_2}(M^p) @>q^{G_2}_{\G_2} >>  H_*^{G_2}(M^q)
   \end{CD}
   $$
   \end{theorem}

\med
\begin{proof}  The proof of this theorem is immediate from the  
definitions, using
the naturality of the Pontrjagin-Thom collapse maps (and thus umkehr  
maps).
We leave the details of this argument to the reader.  \end{proof}
\med
We now discuss a gluing relation held by these operations.  In the  
nonequivariant
setting, a   gluing relation is proved using analytic techniques  
below.  Here
we describe and prove a gluing relation in the general equivariant  
setting.
\med
Let $\G_1$ be a graph with $p$ incoming marked points and $q$ outgoing
marked points.  Let $\G_2$ be a graph with $q$ incoming and $r$  
outgoing marked points.  Say $\G_1 \in \cc_{b_1, p+q}$, and $\G_2 \in  
\cc_{b_2, q+r}$.  By identifying the
the $q$ outgoing leaves (univalent vertices)  of $\G_1$ to the $q$  
incoming leaves of $\G_2$,  this defines
a ``glued"  graph, $\G_1 \#\G_2 \in \cc_{b_1+b_2+q-1, p+r}.$

\begin{figure}[ht]
   \centering
   \includegraphics[height=10cm]{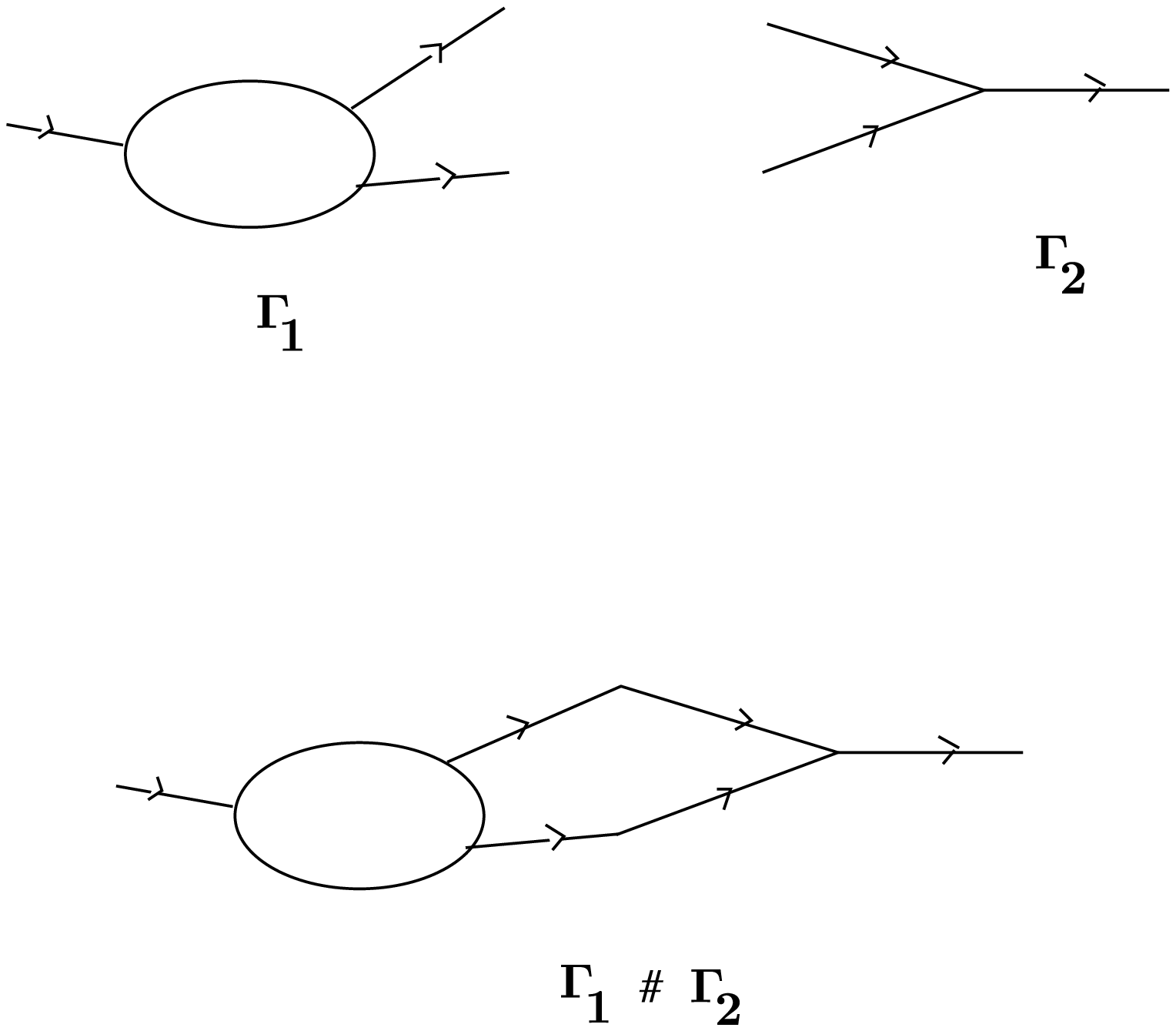}
   \caption{ $\G_1 \# \G_2$}
   \label{fig:figfour}
\end{figure}
   Let $\G_1$ and $\G_2$ be as above.  Consider the homomorphisms
$$\rho_{out} : \ago \to \Sigma_q \quad \rho_{in} : \agt \to \Sigma_q$$
   defined by the induced permutations of the outgoing and incoming  
leaves, respectively.   Let $\ago \times_{\Sigma_q} \agt$ be the  
fiber product of these homomorphisms.  That is,
$$
\ago \times_{\Sigma_q} \agt \subset \ago \times \agt
$$
is the subgroup consisting of those $(g_1, g_2)$ with $\rho_{out} 
(g_1) = \rho_{in}(g_1)$.
Let
$$
p_1 : \ago \times_{\Sigma_q} \agt \to \ago \quad \text{and} \quad  
p_2 : \ago \times_{\Sigma_q} \agt \to \agt
$$
be the projection maps.  There is also an obvious inclusion as a  
subgroup of the automorphism group of the glued graph,
$$
\iota : \ago \times_{\Sigma_q} \agt  \hk Aut(\G_1 \# \G_2)
$$
which realizes  $ \ago \times_{\Sigma_q} \agt$ as the subgroup of $   
Aut(\G_1 \# \G_2)$ consisting of automorphisms that preserve the  
subgraphs, $\G_1$ and $\G_2$.
Similarly, for any pair of homomorphisms, $\theta_1 : G_1 \to \ago$  
and $\theta_2 : G_2 \to \agt$, we have an induced homomorphism
$$
\theta_1 \times \theta_2 :  G_1 \times_{\Sigma_q} G_2 \to  \ago  
\times_{\Sigma_q} \agt  \hk Aut(\G_1 \# \G_2).
$$
We then have the following gluing theorem.

\med
\begin{theorem}\label{gluing}  Let $\G_1$, $\G_2$, $\theta_1 : G_1  
\to \ago$, and $\theta_2 : G_2 \to \agt$ be as above.  Then the   
composition of the graph operations
$$
q^{G_1 \times_{\Sigma_q} G_2}_{\G_2} \circ q^{G_1 \times_{\Sigma_q}  
G_2}_{\G_1} : H_*^{G_1 \times_{\Sigma_q} G_2} (M^p) \to H_*^{G_1  
\times_{\Sigma_q} G_2}(M^q) \to H_*^{G_1 \times_{\Sigma_q} G_2}(M^r)
$$
is equal to the graph operation for the glued graph,
$$
q^{G_1 \times_{\Sigma_q} G_2}_{\G_1 \#\G_2} : H_*^{G_1 \times_ 
{\Sigma_q} G_2} (M^p) \la H_*^{G_1 \times_{\Sigma_q} G_2} (M^r).
$$
\end{theorem}

\med
\begin{proof}
For the sake of ease of notation, we leave off the superscript $ G_1 
\times_{\Sigma_q} G_2 $ in the following description of moduli spaces  
and graph operations.  We wish to prove that
$q_{\G_1 \# \G_2} = q_{\G_2} \circ q_{\G_1}.$

Consider the restriction maps,
$$
\begin{CD}
\cm_{\G_1}(M) @<r_1<< \cm_{\G_1 \#\G_2}(M) @>r_2 >> \cm_{\G_2}( M).
\end{CD}
$$ given by restricting a graph flow on $\G_1 \# \G_2$ to $\G_1$ or $ 
\G_2$,  respectively.  Notice that the following is a pullback square  
of fibrations,
$$
\begin{CD}
\cm_{\G_1 \#\G_2}(M) @>r_2 >>  \cm_{\G_2}( M) \\
   @Vr_1 VV   @VV ev^2_{in} V \\
    \cm_{\G_1}(M) @>>ev^1_{out} > E(G_1\times_{\Sigma_q} G_2) \times_ 
{G_1\times_{\Sigma_q} G_2 } M^q.
    \end{CD}
    $$
    Here the superscripts of the evaluation maps are meant to  
represent  the  graph moduli space on which they are defined.  By the  
naturality of the Pontrjagin-Thom collapse maps used to define the  
umkehr maps (see  in the proof of theorem \ref{tube} as well as the  
more general setup described in \cite{cohenumkehr}), we have the  
following relation:
    \begin{equation}\label{natumkehr}
   r_2 \circ (r_1)_! = (ev_{in}^2)_! \circ ev^1_{out} : H_*^{G_1 
\times_{\Sigma_q}G_2}(M^q)
   \to H_*( \cm_{\G_1 \#\G_2}(M)).
   \end{equation}

   Notice furthermore that we have   commutative diagrams
   $$
   \begin{CD}
\cm_{\G_1 \#\G_2}(M)  @>ev^{1,2}_{in}>>  E(G_1\times_{\Sigma_q}G_2)  
\times_{G_1\times_{\Sigma_q}G_2} M^p     \\
   @Vr_1VV     @V= VV    \\
\cm_{\G_1}(M) @>>ev^{1 }_{in}> E(G_1\times_{\Sigma_q}G_2) \times_{G_1 
\times_{\Sigma_q}G_2} M^p
\end{CD}
$$
and
$$
   \begin{CD}
\cm_{\G_1 \#\G_2}(M)  @>ev^{1,2}_{out}>>  E(G_1\times_{\Sigma_q}G_2)  
\times_{G_1\times_{\Sigma_q}G_2} M^r     \\
   @Vr_2VV     @V= VV    \\
   \cm_{\G_2}( M) @>>ev^{2 }_{out}> E(G_1\times_{\Sigma_q}G_2) \times_ 
{G_1\times_{\Sigma_q}G_2} M^r
\end{CD}
$$

The first of these diagrams  implies, by the naturality of the  
Pontrjagin-Thom collapse maps, that
\begin{equation}\label{nat2}
(ev^{1,2}_{in})_! = (r_1)_! \circ (ev^1_{in})_! : H_*^{G_1\times_ 
{\Sigma_q}G_2}(M^p)) \to
H_*( \cm_{\G_1 \#\G_2}(M)).
\end{equation}

These naturality properties allow us to calculate:
\begin{align}
q_{\G_1\#\G_2} &= ev^{1,2}_{out}\circ (ev^{1,2}_{in})_!, \quad \text 
{by definition} \notag \\
&= ev^{1,2}_{out} \circ (r_1)_! \circ (ev^1_{in})_!, \quad \text{by}  
\, (\ref{nat2}) \notag \\
&= ev^{ 2}_{out} \circ r_2 \circ (r_1)_! \circ (ev^1_{in})_!, \quad  
\text{by the commutativity of the second diagram above}, \notag \\
&= ev^{ 2}_{out} \circ (ev^2_{in})_!\circ ev^1_{out} \circ (ev^1_{in}) 
_!, \quad \text{by} \, (\ref{natumkehr}) \notag \\
&= q_{\G_2} \circ q_{\G_1}, \quad \text{by definition}. \notag
\end{align}
This completes the proof of this theorem.
\end{proof}

\pagebreak

\med
\section{Examples}
In this section we give some examples of the equivariant operations  
$q_\G$.

\subsection{  \rm  The ``Y"-graph and the Steenrod squares. }

Let ${\G_1}$ be the graph:

\begin{figure}[ht]\label{G1}
   \centering
   \includegraphics[height=3cm]{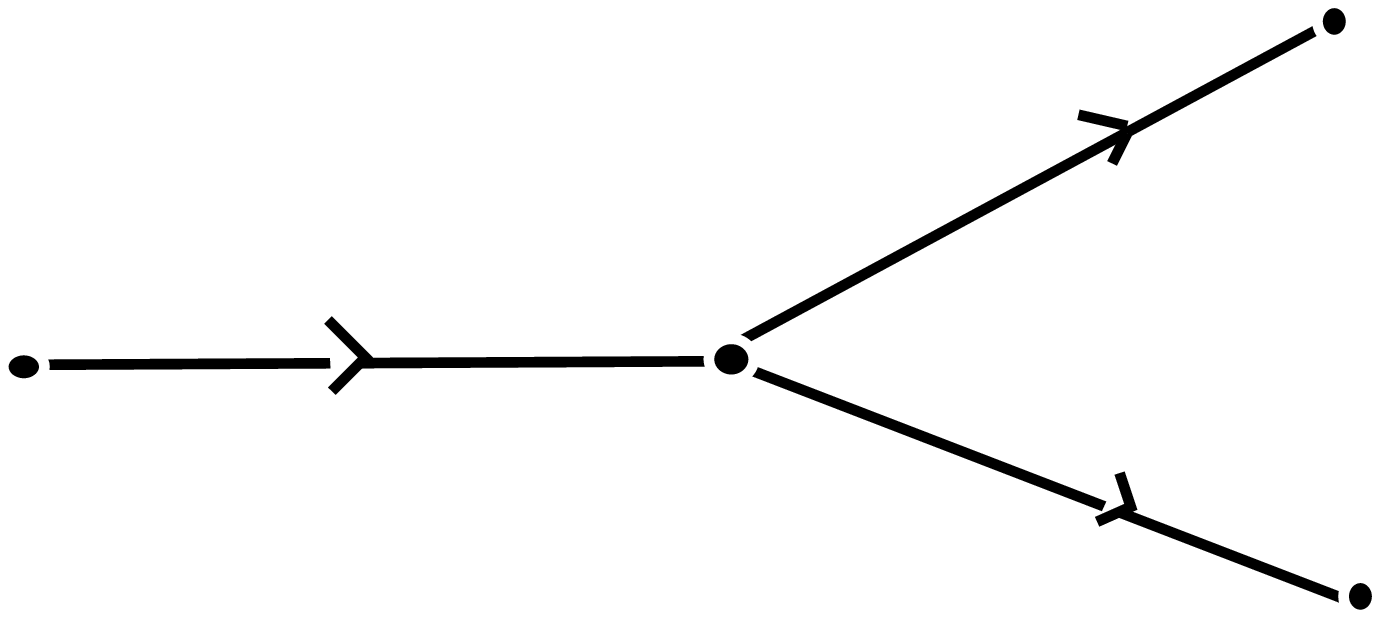}
   \caption{$\G_1$}
\end{figure}
This graph  is a tree with  one incoming and two outgoing leaves.   
The automorphism group  is the group of order 2: $\ag = \bz/2$.  The  
operation $q_{\G_1}$ is therefore a homomorphism,
$$
q_{\G_1} = ev_{out} \circ (ev_{in})_! : H_*(B\bz/2) \otimes H_*(M)  
\to H_{* }( \cm_{\G_1}(M)) \to H_*^{\bz/2}(M \times M).
$$
Since ${\G_1}$ is a tree, $\cm_{\G_1}(M) \simeq B\bz/2  \times M$,  
and clearly $ev_{in} :\cm_{\G_1}(M) \to B(\bz/2) \times M$ is  
homotopic to the identity.  This means $(ev_{in})_!$ is the identity  
homomorphism, and so $q_{\G_1} = ev_{out}$.  But as identified earlier,
$ev_{out} : \cm_{\G_1}(M) \simeq B(\bz/2) \times M \to E\bz/2 \times_ 
{\bz/2} M \times M$ is homotopic to the equivariant diagonal map.  Thus
$$
q_{\G_1} : H_*(B(\bz/2)) \otimes H_*(M) \to H^{\bz/2}_*(M \times M)
$$ is the equivariant diagonal.

Consider the dual map in cohomology with $\bz/2$-coefficients:
$$
(q_{\G_1})^* : H_{\bz/2}^*(M \times M) \la H^*(B\bz/2) \otimes H^*(M).
$$
This is Steenrod's equivariant cup product map \cite{steeneps}.   
Indeed if we considered the nonequivariant operation (associated to  
the homomorphism $\{id\} \hk \ag = \bz/2$),
then the operation
$$
(q^{id}_{\G_1})^* : H^*(M) \otimes H^*(M) \to H^*(M)
$$
is the cup product homomorphism.  In the $\bz/2$-equivariant setting,  
recall that Steenrod defined the Steenrod squaring operations $Sq^j$  
in terms of the equivariant cup product map   in the following way.  
Let  $\alpha \in H^q(M ; \bz/2)$. So $\alpha \otimes \alpha$  
represents a well defined class in $H^{\bz/2}(M \times M; \bz/2)$.  Then
\begin{equation}\label{steenrod}
(q_{\G_1})^*( \alpha \otimes \alpha)  = \sum_{j=0^2q} a^j \otimes Sq^ 
{2q-j}(\alpha).
\end{equation}
Here $a \in H^1(B\bz/2; \bz/2) = H^1(\br \bp^\infty ; \bz/2) = \bz/2$  
is the generator.

\subsection{The Cartan and Adem formulas}

   We now describe how the Cartan and Adem formulas for the Steenrod  
squares follow from the field theoretic properties (invariance and  
gluing) of the graph operations.  Consider the following graph, $\G_2$:

   \begin{figure}[ht]\label{G2}
   \centering
   \includegraphics[height=6cm]{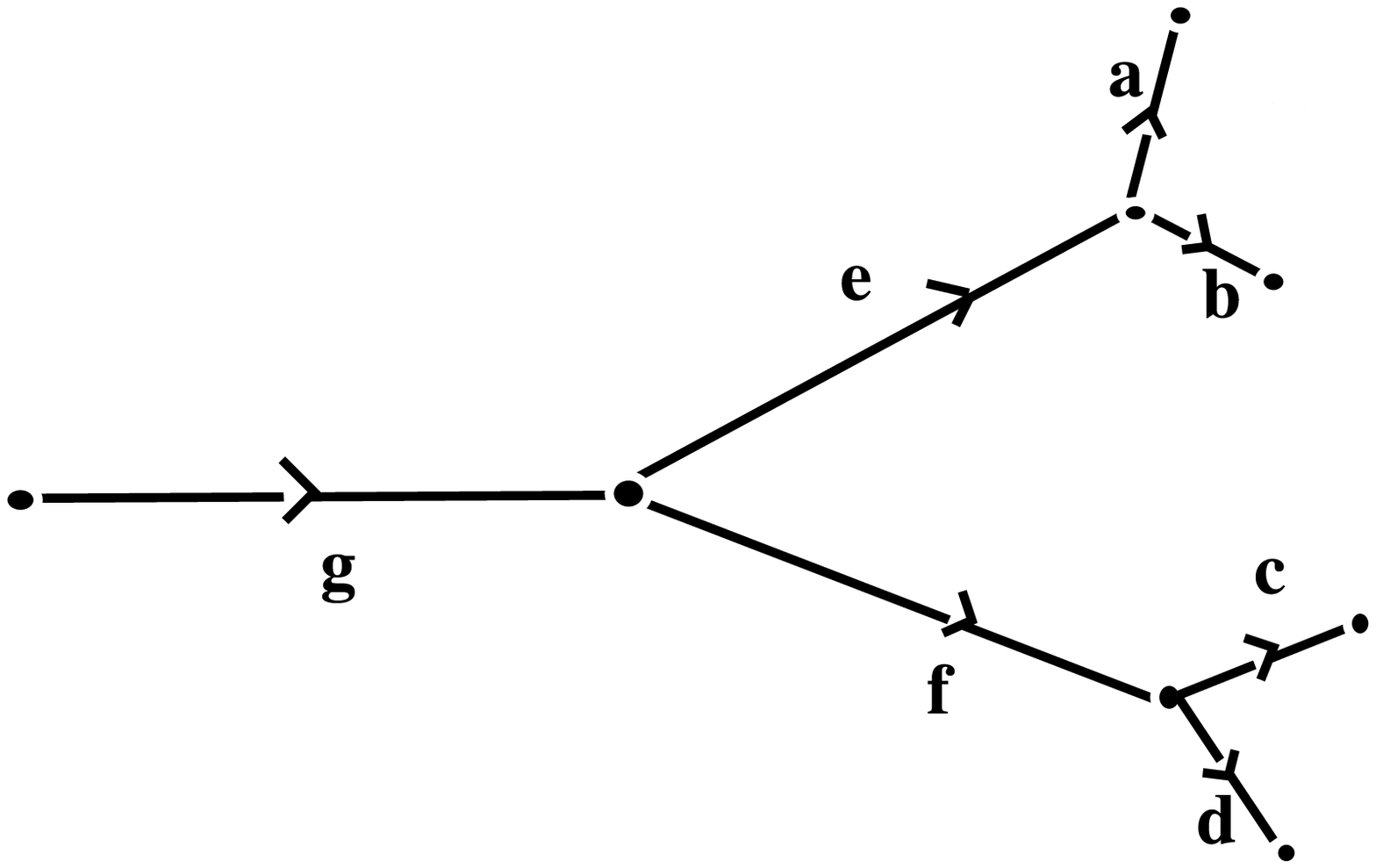}
   \caption{$\G_2$}
\end{figure}

Notice that the automorphism group, $\ag_2 \cong \Sigma_2\int\Sigma_2 
$, the wreath product of the symmetric group with itself.  It sits in  
a short exact sequence, $1 \to \Sigma_2 \times \Sigma_2 \to \Sigma_2 
\int\Sigma_2\to \Sigma_2 \to 1$.  We will view this group as   a  
subgroup of the symmetric group,
$\Sigma_2\int\Sigma_2 \hk \Sigma_4$.  Consider the subgroup $\tau : 
\bz/2 \hk \Sigma_2\int\Sigma_2$ defined by the permutation, $ 
(a,b,c,d) \to  (b,a,d,c)$.  We consider the graph operation in  
cohomology with $\bz/2$-coefficients:
$$
q^\tau_{\G_2} : H^*_\tau (M^4) \to H^*(B(\bz/2)) \otimes H^*(M),
$$
where $H^*_\tau$ is the $\bz/2$-equivariant cohomology determined by  
the embedding $\tau$.

Notice that $\G_2$ is the graph obtained by gluing two copies of the  
Y-graph $\G_1$, each having a single incoming leaf, to the two  
outgoing leaves of a third Y-graph $\G_1$.

   \begin{figure}[ht]
   \centering
   \includegraphics[height=10cm]{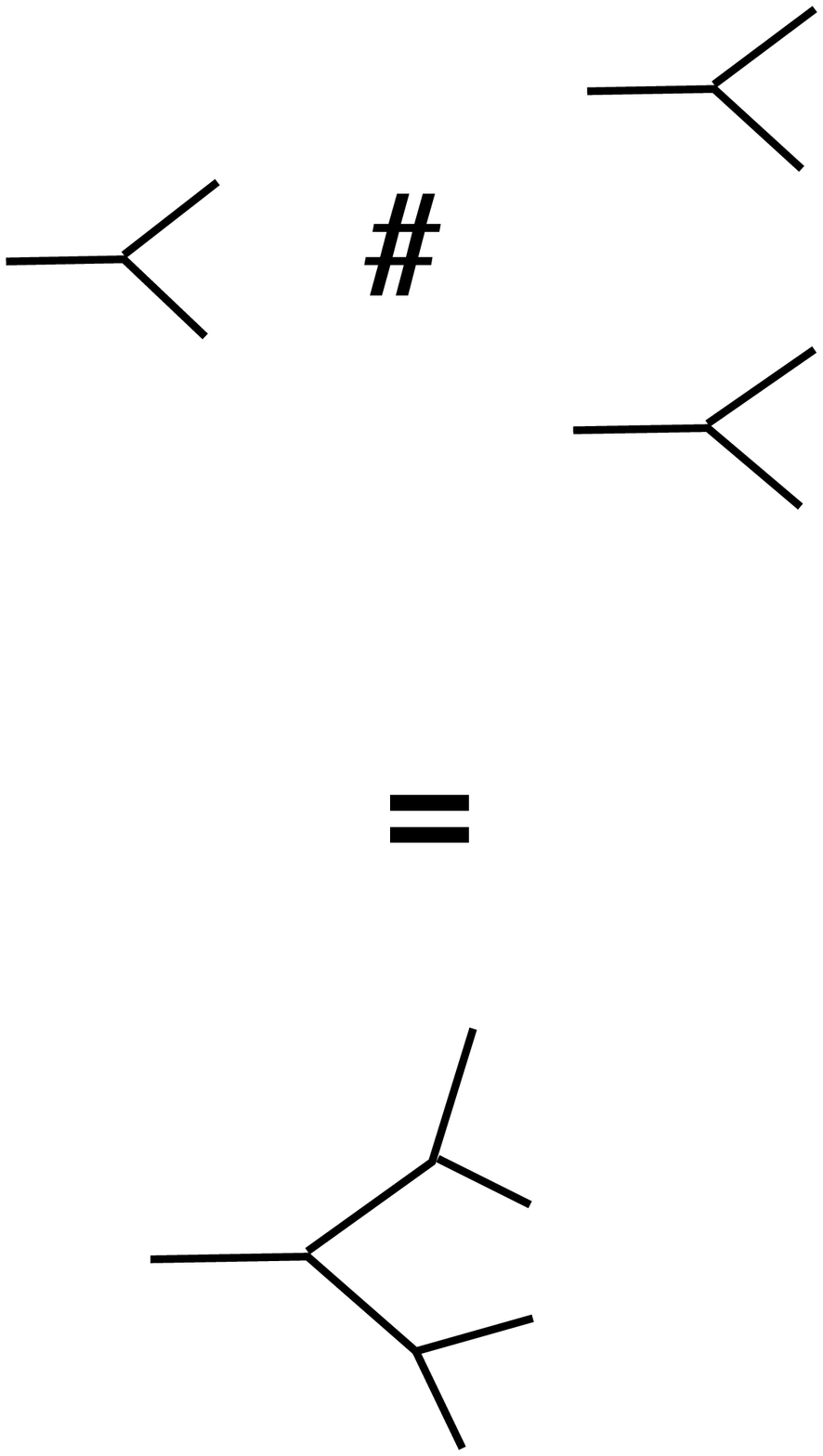}
\caption{ }
\end{figure}

By the gluing formula (theorem (\ref{gluing})) and the description of  
$q_{\G_1}$ above in terms of the (equivariant) cup product, then if $ 
\alpha \in H^q(M)$, $\beta \in H^r(M)$,  then
\begin{equation}\label{equate}
q^\tau_{\G_2}(\alpha \otimes \alpha \otimes \beta \otimes \beta) =  
\sum_{i+s+t=q+r}a^i \otimes Sq^s(\alpha)\cup Sq^t(\beta) \in H^*(B 
(\bz/2)) \otimes H^*(M).
\end{equation}

We now use the invariance property (theorem(\ref{invariance})) to  
understand this operation in another way.  Let $\G_3$ be the  
following graph:

\begin{figure}[ht]
   \centering
   \includegraphics[height=4cm]{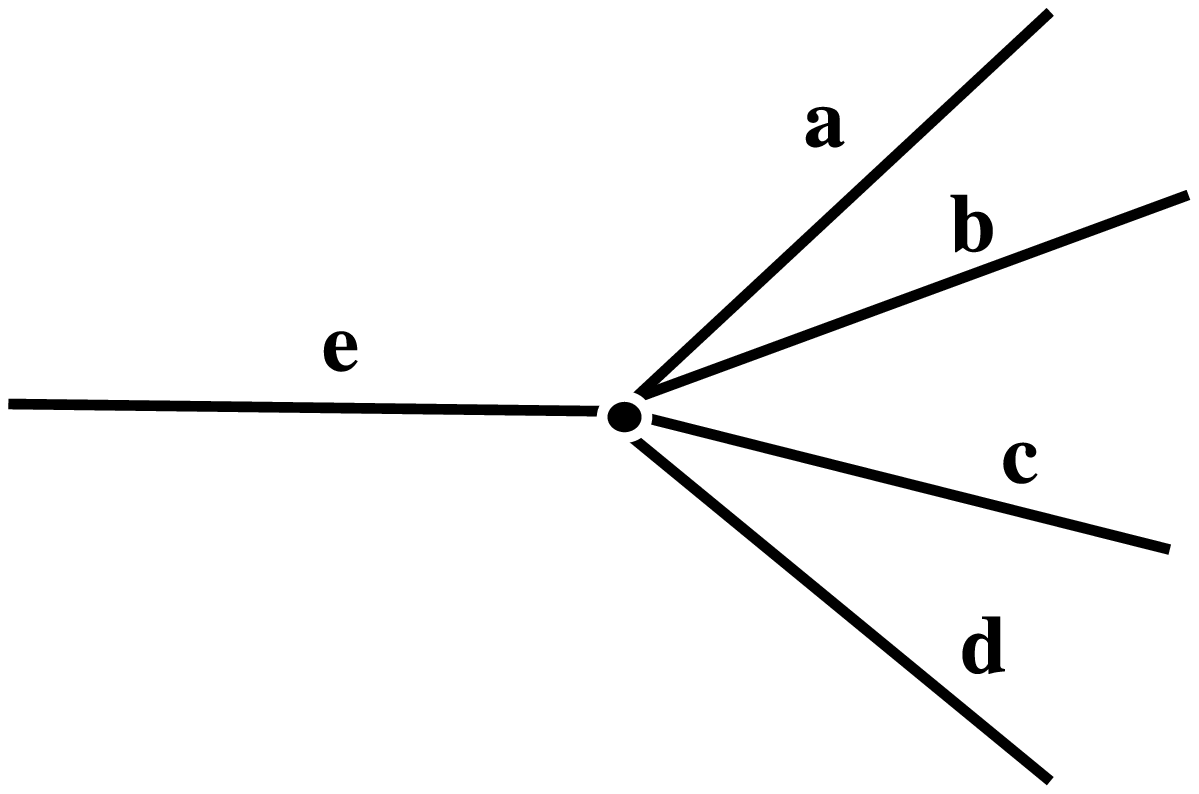}
\caption{$\G_3$}
\end{figure}

Here $Aut (\G_3) = \Sigma_4$, the symmetric group.  Consider the  
morphism,
$$
\theta : \G_2 \to \G_3
$$
obtained by collapsing the edges $e$ and $f$ in figure  6 and then  
permuting the two internal outgoing leaves.  That is, on the level of  
edges,
$$
\theta :g \to e,   \quad  a \to a,  \quad b \to c, \quad c \to b,  
\quad d \to d.
$$
$\theta$ sends the involution $\tau$ on $\G_2$ to the involution $ 
\sigma$ on $\G_3$ defined by the inclusion $\sigma : \bz/2 \hk  
\Sigma_4$, given by the permutation,
$(a, b, c, d) \to (c, d, a, b)$.  So by the invariance property, the  
following diagram commutes:
\begin{equation}\label{compat}
\begin{CD}
H^*_{\tau}(M^4) @>q^\tau_{\G_2} >> H^*(B(\bz/2)) \otimes H^*(M) \\
@V\theta V \cong V   @VV = V \\
H^*_{\sigma}(M^4) @>>q^\sigma_{\G_3} > H^*(B(\bz/2)) \otimes H^*(M) .
\end{CD}
\end{equation}

Now $\theta (\alpha \otimes \alpha \otimes \beta \otimes \beta) =  
\alpha \otimes \beta \otimes \alpha \otimes \beta  \in H^*_\sigma  
(M^4).$  Thus we know from the invariance property and formula (\ref 
{equate}), that
\begin{equation}\label{equate2}
q^\sigma_{\G_3} =  \sum_{i+s+t=q+r}a^i \otimes Sq^s(\alpha)\cup Sq^t 
(\beta).
\end{equation}

On the other hand, consider the morphism
\begin{equation}\label{phi}
\phi : \G_2 \to \G_3
\end{equation}
that also collapses edges $e$ and $f$ but maps edges $a$, $b$, $c$,  
and $d$, to $a$, $b$, $c$, and $d$ respectively.  Since the image of $ 
\sigma : \bz/2 \hk \Sigma_4$ lies in $\Sigma_2 \int \Sigma_2$, the  
invariance property implies
$$
q^\sigma_{\G_2} = q^\sigma_{\G_3} : H^*_\sigma(M^4) \to H^*(B(\bz/2))  
\otimes H^*(M).
$$  But by using figure 7 the gluing formula (theorem \ref{gluing})   
implies that
\begin{align}
q^\sigma_{\G_2} (\alpha \otimes \beta \otimes \alpha \otimes \beta)  
&= q_{\G_1}(\alpha \beta \otimes \alpha \beta) \notag \\
&= \sum_i a^i \otimes Sq^{q+r-i}(\alpha\beta), \quad \text{by (\ref 
{steenrod})}.
\end{align}
Comparing this to formula (\ref{equate2}) yields the Cartan formula,
$$
Sq^m(\alpha\beta) = \sum_{u+v=m}Sq^u(\alpha)Sq^v(\beta).
$$

\med
For the Adem relations, the graph operations don't give us new  
calculational techniques, but they do supply an interesting  
perspective on what calculations are necessary.
Namely, the Adem relations are relations involving iterates of  
Steenrod squaring operations.  From the graph point of view, the  
gluing formula tells us that these
operations  come from considering the graph $\G_2$ given in figure  
6.  As pointed out above, the automorphism group of $\G_2$ is the  
wreath product, $Aut(\G_2) = \Sigma_2\int \Sigma_2$.  In cohomology,  
the graph operation is a homomorphism,
$$
q^*_{\G_2} : H_{\Sigma_2\int \Sigma_2}^*(M^4) \to H^*(B(\Sigma_2\int  
\Sigma_2)) \otimes H^*(M),
$$
and the relevant calculation is $q^*_{\G_2}(\alpha^{\otimes 4})$ for $ 
\alpha \in H^*(M)$.
Now consider the morphism $\phi : \G_2 \to \G_3$ described above.  As  
remarked above, $\ag_4 = \Sigma_4$.  Moreover, in the language
of theorem (\ref{invariance}), $\phi$ is  $\Sigma_2\int \Sigma_2 -  
\Sigma_4$ equivariant.
Therefore by the invariance property, the following diagram commutes:

$$
\begin{CD}
H_{\Sigma_2\int \Sigma_2}^*(M^4) @>q^*_{\G_2}>>  H^*(B(\Sigma_2\int  
\Sigma_2)) \otimes H^*(M) \\
@A\phi AA   @A\iota \otimes 1 AA \\
H_{\Sigma_4}^*(M^4) @>q^*_{\G_3}>>  H^*(B(\Sigma_4)) \otimes H^*(M)
\end{CD}
$$
where $\iota : \Sigma_2\int \Sigma_2 \hk \Sigma_4$ is the inclusion  
as a subgroup.
But since $\alpha^{\otimes 4}$ lies in the image of $\phi : H_ 
{\Sigma_4}^*(M^4) \to H_{\Sigma_2\int \Sigma_2}^*(M^4)$,
we have that $q_{\G_2}^*(\alpha^{\otimes 4})$ is the image of $q_ 
{\G_3}^*(\alpha^{\otimes 4})$ under the map
$$
\iota^*  \otimes 1 : H^*(B(\Sigma_4)) \otimes H^*(M) \to H^*(B 
(\Sigma_2\int \Sigma_2)) \otimes H^*(M).
$$
Now any approach to the Adem relations involves computing the  
relative cohomologies
of $ \iota : \Sigma_2\int \Sigma_2   \hk \Sigma_4$, and in  
particular, the relative equivariant cohomologies of the permutation  
action on $M^4$.  However from this perspective, the reasons these  
calculations are forced upon us, are the gluing and invariance  
properties of the graph operations.
\subsection{Stiefel-Whitney classes}
Consider the following graph, $\G_4$:
\begin{figure}[ht]
   \centering
   \includegraphics[height=3cm]{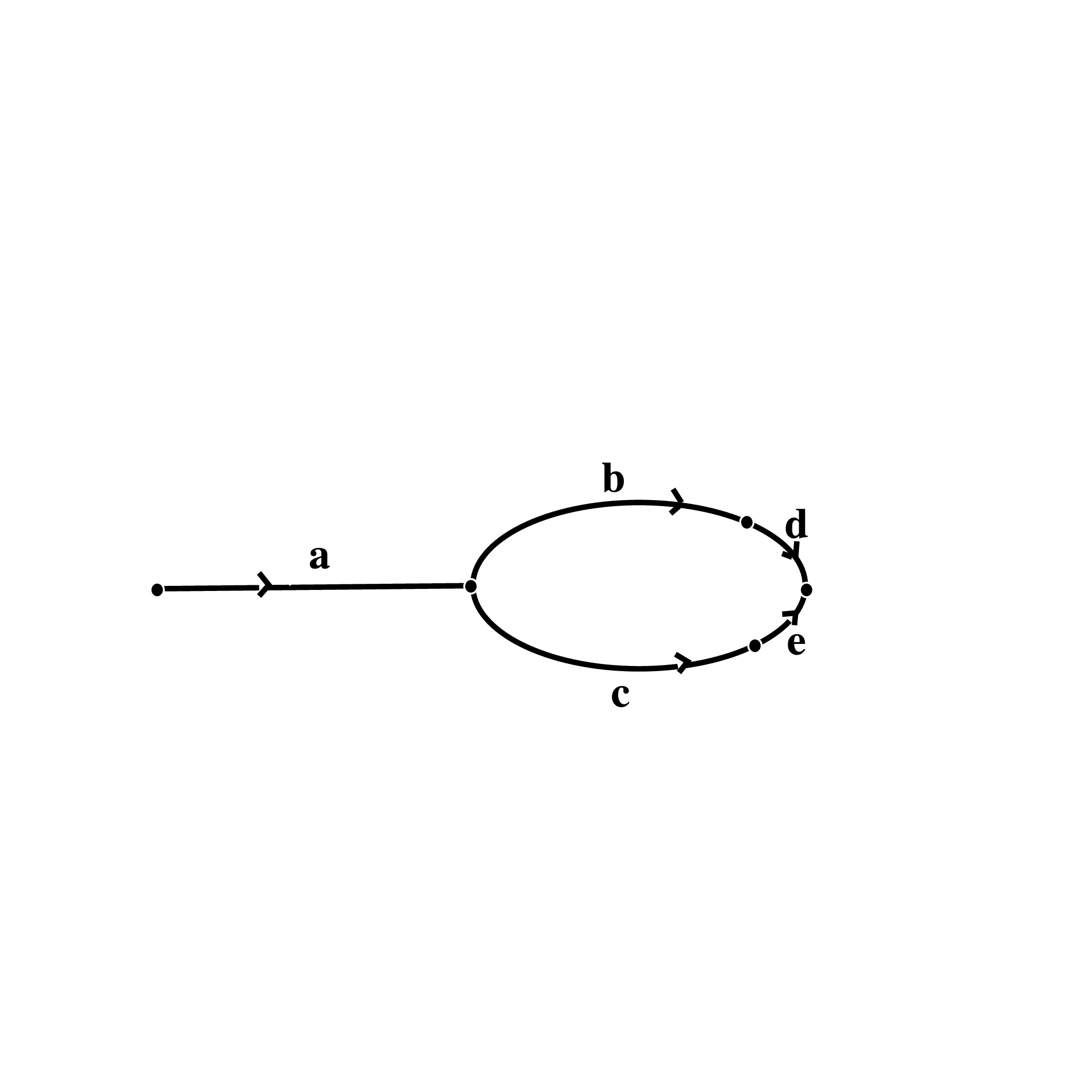}
\caption{$\G_4$}
   \label{fig:fig9}
\end{figure}
In this case the automorphism group $Aut(\G_4) \cong \bz/2$.  Also,  
since there is just one incoming leaf, the operation $q_{\G_4}$ taken  
with $\bz/2$-coefficients is a map,
$$
H_*(B(\bz/2)) \otimes H_*(M) \to \bz/2.
$$
Or, equivalently,  $q_{\G_4} \in H^*(B(\bz/2)) \otimes H^*(M)$.
The following identifies this graph operation.
\begin{theorem}
\begin{align}
q_{\G_4} &= \sum_{i=0}^n a^i \otimes w_{n-i}(M) \notag \\
&\in  H^*(B(\bz/2)) \otimes H^*(M)
\end{align}
where, $w_j(M) \in H^j(M)$ is the $j^{th}$-Stiefel-Whitney class of  
the tangent bundle of $M$, and as above, $a \in H^1(B(\bz/2))$ is the  
generator.
\end{theorem}
\begin{proof}  Let $T \subset \G_4$ be the tree obtained by removing  
the edges $d$ and $e$ in figure 9 above. $T$ has the same  
automorphism group, $Aut(T) = \bz/2$.  By restricting a $\G_4$-graph  
flow  to  $T$, one
obtains an embedding,
$$
\begin{CD}
\cm_{\G_4}(M) @>\rho >\hk > \cm_T(M) \cong (\cs_T(M)/\bz/2)  \times M  
\simeq B\bz/2 \times M.
\end{CD}
$$
   By definition (\ref{qgamma}) the operation $q_{\G_4}$ is given by  
the image of the umkehr map in cohomology,
$$
q_{\G_4}  = \rho^! (1) \in H^n(B(\bz/2) \times M).
$$
To understand this class, notice that the tree $T$ has one incoming  
and two outgoing leaves. Evaluating a graph flow on $T$ at the two  
outgoing leaves defines a map
$$
ev_{out} : \cm_T(M)  \to E(\bz/2) \times_{\bz/2} M \times M
$$
which is homotopic to the equivariant diagonal, $ \Delta : B(\bz/2)  
\times M \to E(\bz/2) \times_{\bz/2} M \times M$.  Furthermore, from  
(\ref{pullback}), the following diagram is a homotopy cartesian square:
$$
\begin{CD}
   \cm_{\G_4}(M) @>\rho >\hk > \cm_T(M) \cong (\cs(T,M)/\bz/2)   
\times M @>\simeq >> B\bz/2 \times M \\
   @V\delta VV &&@VV\Delta V \\
   B(\bz/2) \times M &@>>\Delta > & E(\bz/2) \times_{\bz/2} M \times M
   \end{CD}
   $$
   By the naturality of the Pontrjagin-Thom collapse map and the  
resulting umkehr
   map in cohomology, this homotopy cartesian square implies that
   $$
   \rho^! \circ \delta^* = \Delta^* \circ \Delta^! : H^*(B(\bz/2))  
\otimes M \to H^*(B(\bz/2)) \otimes M.
   $$
   So
   $$
    q_{\G_4}  = \rho^! (1) =  \rho^! \circ \delta^*(1) = \Delta^*  
\circ \Delta^! (1).
    $$
    But by standard properties of umkehr maps, $\Delta^* \circ  
\Delta^! (1)$ is the mod $2$ Euler class of the normal bundle of  
the    equivariant diagonal embedding,
    $\Delta : B\bz/2 \times M \hk E(\bz/2) \times_{\bz/2} M \times M. 
$   Since the normal bundle of the (nonequivariant) diagonal $ 
\Delta : M \to M \times M$ is the tangent bundle, $p: TM \to M$,  the  
normal bundle of the  equivariant diagonal is the equivariant tangent  
bundle,
    $$
    E(\bz/2) \times_{\bz/2} TM \xr{1 \times p} B(\bz/2) \times M,
    $$
    where $\bz/2$ acts fiberwise on $TM$ by multiplication by $-1$.
    The mod $2$ Euler class  is the $n^{th}$-Stiefel-Whitney class of  
this bundle, which
    is given by the sum, $\sum_{i=0}^n a^i \otimes w_{n-i}(TM)$.

    This completes the proof of this theorem.
   \end{proof}

   \subsection{Miscellaneous}
   We conclude this section with a few miscellaneous remarks about  
examples.  Here we work nonequivariantly (i.e we take $q^{1}_\G$  
where $1 \in G$ is the trivial subgroup).
\begin{itemize}
\item  Consider the operation $q^1_\G$ with field coefficients.  Then  
rather than a homomorphism $q^1_\G : H_*(M)^{\otimes p} \to H_*(M)^ 
{\otimes q}$, we may think
of $q^1_\G$ as living in the tensor product, $\bigotimes_p H^*(M)  
\otimes \bigotimes_q H_*(M).$  It is shown in section 9 below  
(corollary \ref{change}),  that if one  changes  the orientation of  
an edge connected to a univalent vertex, one changes the invariant by  
Poincare duality on that factor.
\item The previous remark shows that  when one lets $\tilde \G_1$  
be   the Y-graph as in figure 5,  except that the orientations of all  
three edges are reversed, then the operation
$$
q^1_{\tilde \G_1} : H_*(M) \otimes H_*(M) \la H_*(M)
$$
is the intersection pairing.
\item Consider the graph below with two incoming univalent vertices.
\begin{figure}[ht]
   \centering
   \includegraphics[height=2cm]{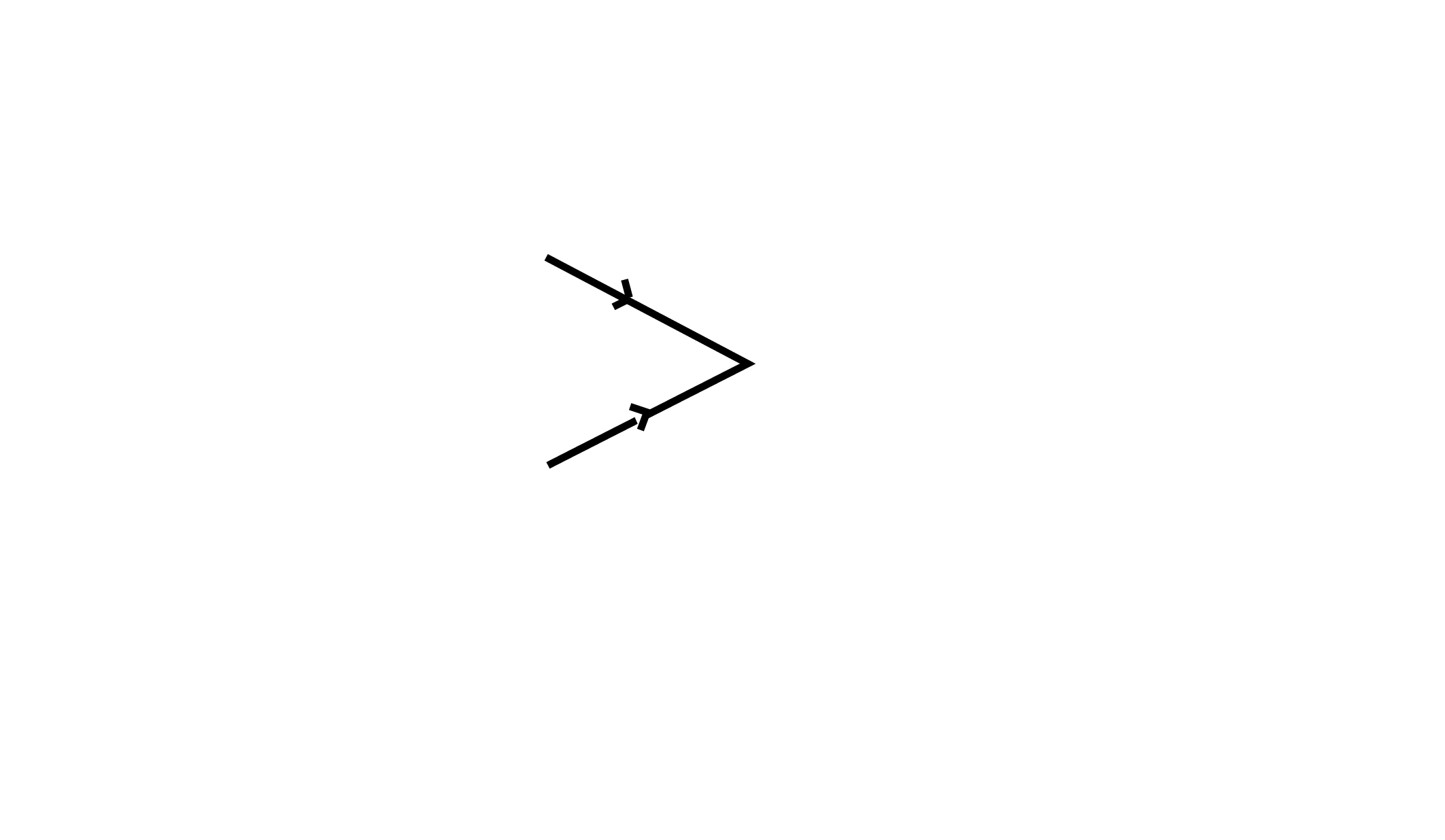}
\caption{$\G_0$}
   \end{figure}
Then the operation $q^1_{\G_0} :H_*(M) \otimes H_*(M) \to k$ is the  
nondegenerate intersection pairing.  Thus the Frobenius algebra  
structure of $H_*(M)$ is encoded in the the Morse field theory  
structure.

\end{itemize}
\section{Transversality.}
\label{sec:trans}

We now give a differential topological construction of the graph  
invariants $q_\G$
defined in section 3.  Throughout this section, and the rest of the  
paper, we will only be using  the automorphism group $\aog  $  
introduced in section 3, that consists of those automorphisms that  
preserve the univalent vertices (leaves).   Since we will be using  
this group exclusively throughout the remainder of the paper, we ease  
notation by simply writing $\ag$ for $\aog$, $\cm_\G(M)$ for $\tilde 
\cm_{\G}(M)/\aog$.

Giving this alternative definition of the graph operations   involves  
studying the smoothness properties of the moduli spaces.  This is the  
main goal of this section.   Our plan for this section is the following.

\med

We will consider the ``graph flow map"
\[ \begin{array}{ccccl}
\Phi&:&\cp_{\g}(M)&\to&\cp_{\g}(TM)\\
&&\gamma&\mapsto&\displaystyle\frac{d\gamma_E}{dt}+\nabla
f_E(\gamma(t))
\end{array}\]
for each edge $E$ of the metric graph $\g_k$,
where $(\g_k,f_E)\in\cs_{\g}$ and where $\cp_{\g}(M)$ (which will be  
defined carefully below) is a space consisting of pairs $(\g_k,  
\gamma)$, where $\g_k$ is a graph over $\G$ (i.e an object in $\cc_\G 
$), and
$\gamma:\g_k\to M$ is a map.
    The map $\Phi$ is a section of the vector bundle
$\cp_{\g}(TM)$ over $\cp_{\g}(M)$ with fibres given by sections of
$\gamma^*(TM)$.

The universal moduli space of metric-graph flows  can therefore be  
thought of as the quotient by $\ag$
   of the zero set of the section $\Phi$
\[\cm_{\g}(M)=\tilde\cm_{\g}(M)/\ag,\ \ \
\tilde\cm_{\g}(M)=\Phi^{-1}(0)\subset\cp_{\g}(M).\]
We will show that it  is a smooth, orientable manifold by an  
application of the implicit function theorem.    Furthermore, we will  
show that the projection map
\[    \begin{array}{c}
     \cgm(M)\\
     \ \ \downarrow\pi \\
     \cgm
     \end{array}  \]
     is smooth, and has virtual codimension (i.e the dimension of $ 
\cgm$ minus the dimension of $\cgm (M)$) equal to
   $-{\rm dim\ }M\cdot\chi(\g)$.   Thus for any submanifold $N\subset 
\cm_{\g}$ transverse to the map
   $\pi$,   the space \[\cm^N_{\g}(M)=\pi^{-1}(N)\] is a
smooth manifold of dimension ${\rm dim\ }M\cdot\chi(\g)+{\rm dim\
}N$.

The evaluation map $ev_v:\cm^N_{\g}(M)\to M$ of a graph flow at a     
univalent
vertex $v\in\g$ allows one to cut down the moduli space further.
Given a Morse function $f$ on $M$ associate to an outgoing
univalent   vertex $v\in\g$ a critical point $a_v$ of $f$ with stable
manifold $\cw^s(a_v)\subset M$.  Then we will see that  $N\subset\cm_ 
{\g}$ can be chosen
transverse to the map  $\pi : \cgm(M) \to \cm_\G $, and so that
$ev_v(\cm^N_{\g}(M))$ intersects $\cw^s(a_v)$ transversely in $M$.   
This will imply that
\[\cm^N_{\g}(M;a_v)=\cm^N_{\g}(M)\cap ev_v^{-1}(\cw^s(a_v))\]
is a smooth manifold of dimension ${\rm dim\ }M\cdot\chi(\g)+{\rm dim\
}N-{\rm index}(a_v).$ By repeated application of this on a
collection of critical points $\vec{a}=\{a_v\}$ of $f$ labeled by the
univalent   vertices of $\g$, one can choose $N$ to get a smooth  
manifold
$\cm^N_{\g}(M;\vec{a})$ of dimension
\[{\rm dim\ }\cm^N_{\g}(M;\vec{a})={\rm dim\ }\cm^N_{\g}(M)
     -\sum_{v \,\, incoming}\left( {\rm dim\ }M-{\rm index}(a_v)\right)
     - \sum_{v \,  \, outgoing} {\rm index}(a_v)\]
where the two sums are taken over the set of incoming and outgoing  
univalent vertices, respectively.

In section 7 we will   prove that the zero dimensional moduli spaces $ 
\cm^N_{\g}(M;\vec{a})$ are compact
and hence one can count the number of points in the moduli space to
get invariants of the manifold. (We will study more general  
compactness issues in section 8.)  These invariants take their values in
formal sums of critical points of $f$ and can be interpreted as
homology classes in the Morse chain complex of $f$.  This will lead   
to a
differential topology construction of the invariants $q_\G$ which we   
do in section 9.

The simple purpose of this section is to prove that transversality can
be arranged.  This is a generalisation of the fact that Morse-Smale
functions exist.  We will use the Sard-Smale theorem so we must first
put a Banach manifold structure on the universal moduli space
$\cm_{\g}(M)$.

\subsection{Mapping spaces.}

To study the flow map
\[\Phi:\cp_{\g}(M)\to\cp_{\g}(TM)\]
we put a Banach manifold structure on the spaces of maps $\cp_{\g}(M)$
and $\cp_{\g}(TM)$, then linearise $\Phi$, and prove regularity.  When
the moduli spaces are finite-dimensional an index calculation gives
the dimension.

Continuous maps from a graph $\g$ to a compact manifold $M$ are best
understood when one equips $\g$ and $M$ with metrics.  More precisely,
equip $M$ with a smooth Riemannian metric and take an oriented metric
graph $\g_k \to \G$ homotopy equivalent to $\g$. (Strictly speaking  
we are considering a point in the geometric realization of $|\cc_\G|  
$, and interpreting it as a metric-graph over $\G$ as discussed in  
section one.) The mapping space
$\cp_{\g}(M)$ consists of continuous maps with square integrable
derivative of all metric graphs $\g_k$ homotopy equivalent to $\g$ as
follows.
\begin{definition}
     For an oriented metric graph $\g_k$ define $\cp_{\g_k}(M)$ to be  
the
     subset of continuous maps $\g_k\to M$ with square integrable
     derivative
     \[\cp_{\g_k}(M)=\left\{\gamma:\g_k\to M\left|\ \gamma {\rm\
     continuous},\
     \int_{\g_k}\left|\frac{d\gamma}{dt}\right|^2dt<\infty
     \right\}.\right.  \]
\end{definition}
Put the $W^{1,2}$ metric on $\cp_{\g_k}(M)$ to give it a Banach
manifold structure, i.e. take {\em continuous} sections $s$ of
$\cv=\gamma^*TM$ satisfying
\[\|s\|^2=\int_{\g_k}\left(\left|\frac{ds}{dt}\right|^2+|s|^2\right) 
dt<\infty.\]
Note that the Sobolev embedding theorem
\[W^{1,2}(E,\cv)\subset C^0(E,\cv)\]
on the interior of edges shows that the requirement of continuity on $s$
can be stated more weakly as continuity at vertices.

We wish to take the union of $\cp_{\g_k}(M)$ over all
$\sigma=(\g_k,\vec{f})\in\cs_{\g}$.  Before doing this, we need a
Banach manifold structure on $\cs_{\g}$.  This is achieved by building
up $\cs_{\g}$ from finite dimensional manifolds, so that the Banach
manifold structure is simply obtained by taking finite objects in
$\cs_{\g}$.

In the construction of $\cs_{\g}$, take $V\subset C^{\infty}(M)$ to be
an $N$-dimensional vector space and allow only $\g_k$ for $k$ less
than $N$.  That  is, $g_k$ is a point on the $k$-skeleton of $|\cc_ 
\G|/\ag$, for $k <N$.   (We have unnecessarily chosen the bound on $k 
$ to coincide
with the dimension of $V$.)  Then $\cs_{\g}$ is built up out of the
union of $\cs^{(N)}_{\g}\subset\cs_{\g}$.
\begin{definition}
\begin{itemize}
     \item Define
     \[\cp_{\g}(M)=\bigcup_{(\g_k,\vec{f})\in\cs_{\g}}\{\g_k,\vec{f}\}
     \times \cp_{\g_k}(M)\]
     and equip it with the topology induced from $\cs_{\g}$ and
     $\cp_{\g_k}(M)$.

    \item  Define the vector bundle $\cp_{\g}(TM) \to \cp_{\g}(M)$ so  
that
for each $\gamma\in\cp_{\g}(M)$ the fibre over $\gamma$
consists of $L^2$ maps.
\begin{eqnarray*}
     \left(\cp_{\g}(TM)\right)_{\gamma}&=&\left\{(\gamma,\xi):\g_k\to  
TM\left|\
     \int_{\g_k}|\xi|^2dt<\infty\right\}.\right.
\end{eqnarray*}
\end{itemize}
\end{definition}

\subsection{Surjectivity}
In  the previous section we proved that the   space $\cp_\G(M)$
is a manifold.  We next show that $0$ is a regular
value of $\Phi$.

\begin{theorem}   \label{th:reg}
        $  \Phi:\cp_{\g}(M)\to\cp_{\g}(TM)$ intersects the zero  
section transversally.
\end{theorem}
\begin{proof}
     The tangent space at a point $(\g_k,\vec{f},\gamma)\in\cp_{\g}(M)$
     is given by
     \[ T_{(\g_k,\vec{f},\gamma)}\cp_{\g}(M)=T_{(\g_k,\vec{f})}
     \cs_{\g}\oplus W^{1,2}(\g_k,\cv)\]
     where $\cv=\gamma^*(TM)$ is a vector bundle over $\g$.  (If $M$ is
     orientable then $\cv$ is trivial so $W^{1,2}(\br^n)$ suffices.)   
The
     linearisation of $\Phi$ decomposes into $D\Phi=(I,D_1+D_{\g_k})$  
where
     $I$ is the identity on the $T\cs_{\g}$ part and
     \[ D_1+D_{\g_k}:T_{(\g_k,\vec{f})}\cs_{\g}\oplus W^{1,2}(\g_k,\cv)
     \to L^2(\g_k,\cv)\ .\]
     We must show that for all points of the universal moduli space
     $(\g_k,\vec{f},\gamma)\in\cm_{\g}(M)$,
     $D\Phi_{(\g_k,\vec{f},\gamma)}$ is surjective and has a right  
inverse.

     A tangent vector in $T_{(\g_k,\vec{f},\gamma)}\cp_{\g}(M)$ is
     given by a triple $(\lambda,\vec{h},s)$ where
     $\lambda=\{\lambda_E\}$ is the infinitesimal change in the length
     of $E$, $\vec{h}=\{h_E\}$ is the infinitesimal change in the
     smooth function labeling $E$, and $s=\{s_E\}$ is a section of the
     vector bundle $\cv=\gamma^*TM$ over $\g$.
     \begin{eqnarray*}
     \gamma_E&\mapsto&\gamma_E+\epsilon s_E\\
     l_E&\mapsto&l_E+\epsilon\lambda_E\\
     f_E&\mapsto&f_E+\epsilon h_E.
     \end{eqnarray*}
     To linearise
     $\Phi(\g_k,\vec{f},\gamma)=\{\dot{\gamma_E}+\nabla f_E(\gamma_E)\}$
     assume for the moment that $l_E>0$ and reparametrise $E$
     by $\tau\in[0,1]$ so $t=\tau l_E$.  Then
     $\ell_Ed\gamma/dt=d\gamma/d\tau$ and $\ell_E\left(d\gamma_E/dt+ 
\nabla
     f_E(\gamma_E)\right)= d\gamma_E/d\tau_E+\ell_E\nabla f_E 
(\gamma_E).$
     \[    \ell_E\left(\Phi+\epsilon\Delta\Phi\right)=
     \frac{d}{d\tau_E}(\gamma_E+\epsilon s_E)+(\ell_E+\epsilon\lambda_E)
     \nabla(f_E+\epsilon h_E)\ \ \ \ \ \ \ \ \ \ \ \ \ \ \ \ \ \ \ \  
\ \ \
     \ \ \ \ \ \ \ \ \ \ \ \ \ \ \ \ \ \ \ \ \ \]
     \[=\frac{d}{d\tau_E}\gamma_E+\ell_E\nabla f_E(\gamma_E)
     +\epsilon
     (\frac{d}{d\tau_E}s_E+\ell_E\nabla\nabla f_E\cdot s_E+\lambda_E
     \nabla f_E+\ell_E\nabla h_E)\]
     \[=\ell_E\left(\frac{d}{dt_E}\gamma_E+\nabla f_E(\gamma_E)+\epsilon
     (\frac{d}{dt_E}s_E+\nabla\nabla f_E\cdot s_E+\frac{\lambda_E} 
{\ell_E}
     \nabla f_E+\nabla h_E)\right).\]
     hence for $\ell_E>0$
     \[ D_1(\lambda,\vec{h})+D_{\g_k}s=(\frac{\lambda_E}{\ell_E} 
\nabla f_E
     +\nabla h_E)+(\dot{s}_E+\nabla\nabla f_E\cdot s_E).\]
     For the case $\ell_E=0$ we first need to understand more about the
     cokernel of $D_{\g_k}$.  Since $L^2(\g,\cv)$ is a Hilbert space
     the cokernel of $D_{\g_k}$ can be identified with the orthogonal
     complement of its image.  Thus
     \[\coker D_{\g_k}=\{ r\in L^2(\g,\cv)|\langle r,D_{\g_k}
     \phi\rangle=0\ {\rm for\ all}\ \phi\in C^{\infty}_0(\g)\}\]
     which gives good local behaviour of an element $r$ of the
     cokernel on the interior of an edge and at a vertex.
     \begin{lemma}   \label{th:coker}
     On the interior of any edge $E\subset\g_k$, an element
     $r\in\coker D_{\g_k}$ is smooth and satisfies
     \begin{equation} \label{eq:coker}
         \dot{r}_E-(\nabla\nabla f_E)^T\cdot r_E=0.
     \end{equation}
     At a vertex $v\in\g_k$, $r$ is free to be discontinuous up to
     the codimension $1$ condition
     \begin{equation}  \label{eq:cokerv}
         \sum_{E\ni v}(-1)^{v(E)}r_E(v)=0
     \end{equation}
     where $v(E)=0$ (or 1) when $E$ is incoming (outgoing).
     \end{lemma}
     \begin{proof}
     The first part of the lemma is standard so we defer that to
an appendix.  To prove (\ref{eq:cokerv}) consider $\phi\in
     C^{\infty}_0(\g)$ whose support lies in a neighbourhood of the
     vertex $v\in\g$.  Then
     \begin{eqnarray*}
         0&=&\int_{\g}\langle r,\dot{\phi}+A\phi\rangle dt\\
         &=&\sum_{E\ni v}(-1)^{v(E)}r_E(v)\phi(v)-\int_{\g}\langle
         \dot{r},\phi\rangle dt+\int_{\g}\langle A^Tr,\phi\rangle
         dt\\
         &=&\sum_{E\ni v}(-1)^{v(E)}r_E(v)\phi(v).
     \end{eqnarray*}
     \end{proof}
     When $\ell_E=0$ we can make sense of $D_1(\lambda,\vec{h})$ weakly
     in $L^2$ as follows.  Along an edge $E$, $\nabla f_E\in\ker D_ 
{\g_k}$
     since it gives an infinitesimal change in parametrisation.  Thus,
     for any $r_E\in\coker D_{\g_k}$, \[\frac{d}{dt}\langle r_E(t), 
\nabla
     f_E(t)\rangle= \langle \dot{r}_E-(\nabla\nabla f_E)^T\cdot r_E
     ,\nabla f_E(t)\rangle+ \langle r_E(t),D_{\g_k}\nabla
     f_E(t)\rangle=0\] so \[\langle r_E(t),\nabla
     f_E(t)\rangle_2=\ell_E \langle r_E(0),\nabla f_E(0)\rangle\]
     since $r_E(0)$ makes sense (as a one-sided limit.)  Therefore,
     \begin{equation} \label{eq:D1}
     \langle r_E,D_1(\lambda,0)\rangle_2=
     \frac{\lambda_E}{\ell_E}\langle r_E,\nabla f_E\rangle_2=
     \lambda_E \langle r_E(0),\nabla f_E(0)\rangle
     \end{equation}
     which makes sense when $\ell_E=0$.

     When building up $\cs_{\g}$ from finite-dimensional $V\subset
     C^{\infty}(M)$, choose $V$ to be generated by $\{
     f_1,\ldots,f_N\}$ such that $\{\nabla f_1,\ldots,\nabla f_N\}$
     span $T_xM$ at every point $x\in M$.  At a vertex $v\in\g_k$,
     identify $\g_k$ with $\g_{k+d}\to\g_k$ a metric graph with $d=\dim
     M$ extra edges that contract to $v$ and associate to these edges
     smooth functions with gradients spanning $T_{\gamma(v)}M$.  Then
     an infinitesimal increase in the length of one of the zero length
     edges $E'$ at $v$ gives any direction $\nabla f_{E'}(0)$ and thus
     from (\ref{eq:D1}) $\lambda_{E'} \langle r_{E'}(0),\nabla
     f_{E'}(0)\rangle=0$ implies $r_E(0)=r_{E'}(0)=0$ and since $r_E$
     satisfies an ODE this implies that $r_E(t)\equiv 0$ so $D\Phi$ is
     onto.

     It is proven in an appendix that $D_{\g_k}$ is Fredholm and it is
     a standard fact that this implies $D\Phi$ has a right inverse.
\end{proof}

\begin{theorem}
     The universal moduli space of graph flows $\cm_{\g}(M)$ is a
     smooth Banach manifold. The projection map
     \[ \pi:\cm_{\g}(M)\to\cm_{\g}\]
     has virtual codimension $-\dim M\cdot\chi(\g)$. Its cover
     $ \tilde\cm_{\g}(M) $  inherits a natural coorientation from an
     orientation of $M$.
\end{theorem}
\begin{proof}
     Since   $\Phi$ intersects the zero section transversally,  from  
the implicit function
     theorem it follows that $\tilde\cm_{\g}(M)=\Phi^{-1}(0)$ is a  
manifold
     and since $\ag$ acts freely on $\tilde\cm_{\g}(M)$ so too is
     $\cm_{\g}(M)=\tilde\cm_{\g}(M)/\ag$.

     The projection map
     \[ \pi:\tilde\cm_{\g}(M)\to\cs_{\g}\]
     is a Fredholm map between Banach manifolds since it is clear that
     $\ker D\pi=\ker D_{\g_k}$ and with a little more thought one can
     see that $D_1$ induces an isomorphism between $\coker D\pi$ and
     $\coker D_{\g_k}$.  Since $D_{\g_k}$ is Fredholm (see
     the appendix), $D\pi$ is Fredholm with index equal to
     ${\rm index\ }D_{\g_k}$.

     \med
     \begin{lemma}  \label{th:index} The index of the operator $D_ 
{\G_k}$ is given by
     $${\rm index\ }D_{\g_k}={\rm dim\ }M\cdot\chi(\g).$$
     \end{lemma}

     \med
     \begin{proof}
     The index remains unchanged on the continuous families of
     operators
     \[ D_{\g_k}(\lambda)=\frac{d}{dt}+\lambda\nabla\nabla f,\
     \lambda\in[0,1]\]
     so we may replace $D_{\g_k}$ by $\frac{d}{dt}$ which is
     differentiation on $\gamma^*(TM)$, an $\br^d$ bundle over
     $\g$.  The operator $\frac{d}{dt}$ is well-defined even if
     $\gamma^*(TM)$ is non-trivial since we choose trivialisations
     so that $\gamma^*(TM)$ is the sum of a trivial bundle and a
     non-trivial line bundle with transition function
     multiplication by $-1$.  Thus the transition function commutes
     with $\frac{d}{dt}$.  By Lemma~\ref{th:coker} $\frac{d}{dt}$
     is self-adjoint up to boundary terms and elements $r$ of the
     cokernel also satisfy $\frac{d}{dt}r=0$.

     If the bundle $\gamma^*(TM)$ is trivial then
     $\ker\frac{d}{dt}$ consists of constant sections so it has
     dimension $\dim M$.  Elements of the cokernel are constant
     along edges but may be discontinuous at vertices, satisfying a
     codimension one condition there.  The cokernel is spanned by
     sections constant around a cycle in $\g$ representing a
     generator of $H_1(\g)$ and zero outside the cycle.  Thus it
     has dimension $\dim M\cdot b_1(\g)$ and ${\rm index\
     }D_{\g_k}={\rm dim\ }M\cdot\chi(\g).$

     If the bundle $\gamma^*(TM)$ is non-trivial then
     $\ker\frac{d}{dt}$ has dimension $\dim M-1$ since it is
     trivial on the non-trivial sub-line bundle of $\gamma^*(TM)$.
     This is because there has to be a cycle in $\g$ on which there
     is an odd number of transition functions given by
     multiplication by $-1$.  But then since any element of the
     kernel is a constant $c$, we must have $c=-c=0$.  To see that
     the cokernel has dimension ${\rm dim\ }M\cdot b_1(\g)-1$ it is
     enough to consider the non-trivial sub-line bundle of
     $\gamma^*(TM)$ since the argument above takes care of the
     trivial sub-bundle.  Note that $H_1(\g)$ can be generated by
     $b_1(\g)$ cycles such that $\gamma^*(TM)$ is trivial around
     all but one cycle and it has an odd number of transition
     functions given by multiplication by $-1$ along one cycle.
     (To see this, take $b_1(\g)$ cycles in $\g$ that generate
     $H_1(\g)$.  If the bundle is non-trivial on two cycles
     $\alpha$ and $\beta$ in the generating set, then replace
     $\alpha$ and $\beta$ by $\alpha+\beta$ and $\beta$.  Continue
     this until the bundle is non-trivial on only one generator.)
     Again the cokernel is spanned by sections constant around a
     cycle in the generating set of $H_1(\g)$, hence it is zero on
     the non-trivial cycle and constant on the other $b_1(\g)-1$
     cycles.  As before ${\rm index\ }D_{\g_k}={\rm dim\
     }M\cdot\chi(\g).$
     \end{proof}
     The normal bundle of $ \pi : \tilde\cm_{\g}(M) \to \cs_\G$ at
     $\sigma=(\g_k,\vec{f})\in\cs_{\g}$ is canonically isomorphic to
     $(\coker D\phi)^*=(\coker D\g_k)^*$.  A coorientation of
     $\tilde\cm_{\g}(M)$ in $\cs_{\g}$ which is an orientation of its
     normal bundle is thus a section of the line bundle $\wedge^{\rm
     max} (\coker D\Phi)^*$ over $\tilde\cm_{\g}(M)$.  This line bundle
     coincides with the determinant line bundle \[\det
     D\Phi=\wedge^{\rm max}\ker D\Phi\otimes\wedge^{\rm max} (\coker
     D\Phi)^*\] since $D\Phi$ is surjective.  The determinant line
     bundle extends to a locally trivial line bundle over all of
     $\cp_{\g}(M)$, and since $\cp_{\g}(M)$ is contractible the
     determinant line bundle is globally trivial .  Thus
     $\pi(\tilde\cm_{\g}(M))$ is coorientable.  A coorientation is {\em
     canonically} determined from an orientation on $M$ via the
     evaluation map.
\end{proof}
The following is a corollary of what we have just proved.  It is the  
main result of this section.  It gives the smoothness of the finite  
dimensional moduli spaces.

\med
\begin{theorem}   \label{th:sman}
     For a generic submanifold with boundary $N\subset\cm_{\g}$, the
     moduli space $\cm^N_{\g}(M)$ is a manifold with boundary, and has
     dimension
     \begin{equation*}
     {\rm dim\ }\cm^N_{\g}(M)={\rm dim\ }M\cdot\chi(\g)+{\rm
     dim\ }N.
     \end{equation*}
\end{theorem}

\med
\begin{proof}
     A strong version of the Sard-Smale theorem guarantees that    
that any submanifold of $\cs_{\g}$ can be perturbed to
     an arbitrarily close submanifold $N$ that is transverse to
     $\pi: \cm_{\g}(M)\to \cm_\G$, with its boundary transverse to
     $\pi(\cm_{\g}(M)) $.  Hence $\cm^N_{\g}(M)=\pi^{-1}(N)$ is a
     manifold with boundary.  The dimension formula follows immediately
     from the codimension of $\pi(\cm_{\g}(M))$ in $\cm_{\g}$.  (In the
     case   $b_1(\g)=0$,    this means that $\cm_{\g}(M)$ maps
     onto $\cs_{\g}$ with fibre of dimension ${\rm dim\ }M$.)
\end{proof}
{\em Remark.} The moduli space $\cm_{g,n}(M,\beta)$ of stable maps
from genus $g$ curves with $n$  marked   points to a variety $M$ (or
symplectic manifold) $\gamma:\Sigma\to M$ with image representing
$\beta=[\gamma(\Sigma)]\in H_2(M)$ uses the entire parameter space
$N=\cm_{g,n}$ and has (real) dimension
\[ {\rm dim\ }\cm^N_{g,n}(M)=\frac{1}{2}{\rm dim\
}M\cdot\chi(\Sigma)+{\rm dim\ }N+2\langle
c_1(M),[\gamma(\Sigma)]\rangle.\]
The analogue of Lemma~\ref{th:index} is the Riemann-Roch formula given
in terms of {\em complex} dimensions
\[ {\rm dim\ }H^0(\gamma^*(TM))-{\rm dim\ }H^1(\gamma^*(TM))
=\frac{1}{2}{\rm dim\
}\gamma^*(TM)\cdot\chi(\Sigma)+{\rm degree\ }\gamma^*(TM).\]
Both Riemann-Roch and Lemma~\ref{th:index} are index theorems relating
the index of a differential operator to topological information.
The topological term $\langle c_1(M),[\gamma(\Sigma)]\rangle$
specifies different components of the moduli space and is detected in
the dimension formula.  Similarly, connected components of the moduli
space of graph flows have constant $\langle
w_1(M),[\gamma(\g)]\rangle$.  One might expect different dimensions
for different connected components, however the term $\langle
w_1(M),[\gamma(\g)]\rangle$ is not detected in the dimension formula,
although curiously it does appear in the calculation of the dimension.

\section{Zero dimensional moduli spaces and counting.}
\label{sec:zerod}

Given a Morse function $f$ on $M$, the {\em Morse complex} of $f$ is a
chain complex generated by the critical points of $f$, with boundary
maps obtained from counting gradient flows.  Using this
description of the homology of $M$, we will show how the graph  
operations defined earlier can be defined geometrically on  the chain  
level as maps between formal  linear combinations of
critical points of $f$.  The graph moduli spaces are ideal for
defining such chain level maps.

The stable and unstable manifolds of critical points of $f$ represent
homology and cohomology classes on $M$ and they intersect the image of
the moduli space of graph flows under the evaluation map.  This will  
allow us to give
a realisation of the umkehr maps defined in section 3 from tensor  
products of the homology and cohomology of $M$ to
   to the homology of the  moduli space of
graph flows.

The intersection of the image of the evaluation map with stable and
unstable manifolds will  be interpreted in terms of a moduli space of
graph flows for a {\em non-compact} graph.  Non-compact edges will  
map to
gradient flows of the Morse function $f$. As we will see,  the Morse  
condition---that
the critical points of $f$ are non-degenerate---arises because the
gradient flows live on a non-compact graph.  Until now the degeneracy
of critical points of a smooth function on $M$ has been of no concern
to the construction of the moduli spaces because only compact graphs
have been used.

For the remainder of the paper we work with the non-compact graph
$\gn$, obtained from $\g$ by adding, for each  univalent   vertex $v 
\in\g$,
a non-compact edge $E_v$ oriented incoming or outgoing according to
whether $v$ is incoming or outgoing.      A graph flow is a  
continuous map $\gamma:\gn\to M$ which is
the previously defined graph flow on $\g$, and on non-compact edges it
is the gradient flow of the Morse function $f$.

To define the moduli space of graph flows we now specify a collection
of critical points $\vec{a}=\{ a_v\}$ of the Morse function $f$ and
require that the gradient flow of the non-compact edge $E_v$, of
$\gn$ converges to the critical point $a_v$.  The graph flow map is
defined on appropriate path spaces (defined below) and is given by:
\[ \Phi(\gamma)=\displaystyle\frac{d\gamma_E}{dt}+\nabla f_E(\gamma 
(t))\]
where $E\subset\gn_k$ varies over all edges of $\g_k$ and non-compact
edges of $\gn_k$.   In this notation $f_E$ is the restriction of $f$  
to the edge  $E$.

The universal moduli space of graph flows of $\gn$ is notated by
\[    \begin{array}{c}
     \tilde\cm_\G(M; \vec{a})=\Phi^{-1}(0)\\
     \ \ \downarrow\pi \\
     \cs_{\g}
     \end{array}  \]
where $\g$ encodes $\gn$ through its oriented  univalent   vertices.   
Notice that
the  space of structures  remains $\cs_{\g}$, (i.e it has not changed  
even though we are now working with the enlarged graph $\tilde \G$),  
because  there is a fixed function $f$ labeling all the noncompact  
edges.

Most of the results for compact graphs generalise to these particular  
non-compact
graphs.  The following theorem encapsulates these generalisations.

\begin{theorem}   \label{th:moduli}
     For a generic submanifold with boundary $N\subset\cm_{\g}$, the
     moduli space of graph flows $\cm^N_{\g}(M;\vec{a})=\pi^{-1}(N)/ 
\ag \subset \cm_{\g}(M;\vec{a}) $ is a
     manifold with boundary, of dimension
     \begin{equation*}
     {\rm dim\ }\cm^N_{\g}(M;\vec{a})={\rm dim\ }\cm^N_{\g}(M)
     -\sum_{v \, \, incoming}\left( {\rm dim\ }M-{\rm index}(a_v)\right)
     - \sum_{v \, \, outgoing} {\rm index}(a_v).
     \end{equation*}
\end{theorem}

{\em Remarks.}  1. The moduli space $\cm^N_{\g}(M;\vec{a})$ is in  
general
not compact.

2. There is a canonical orientation on $\tilde{\cm}^N_{\g}(M;\vec{a}) 
$ induced by an orientation on $M$.

\begin{proof}
     To prove the theorem we must define the Banach manifold structure
     on the mapping spaces, construct the universal moduli space of
     graph flows, prove that the projection $\pi$ to the  structure
     space is Fredholm, calculate its index and prove regularity.
      Except for the proof that the operator is Fredholm, these results
     require only small adjustments to the compact graph case.

     For an oriented metric graph $\g_k$ with  univalent   vertices  
attach a
     half-line $E_v$ to each  univalent   vertex $v$ to get the non- 
compact
     graph $\gn_k$.  The non-compact edge is oriented according to the
     orientation of $v$ and this is realised in the parametrisation of
     incoming $E_v$ by $t\in(-\infty,0]$ and outgoing $E_v$ by
     $t\in[0,\infty)$.
     \begin{definition}
     Define $\cp_{\g_k}(M;\vec{a})$ to be the subset of continuous
     maps from $\gn_k\to M$, that converge on non-compact edges to
     $a_i$, with square integrable derivative
     \[\cp_{\g_k}(M;\vec{a})=\left\{\gamma:\g_k\to M\left|\ \gamma
     {\rm\ continuous},\ \lim_{t\to(-)\infty}\gamma_{E_v}(t)=a_v,\
     \int_{\gn_k}\left|\frac{d\gamma}{dt}\right|^2dt<\infty
     \right\}\right.  \] and for any section $s$ of the vector
     bundle $\cv=\gamma^*TM$ over $\gn_k$ define its norm using the
     $W^{1,2}$ metric
     \[\|s\|^2=\int_{\gn_k}\left(\left|\frac{ds}{dt}\right|^2+|s|^2 
\right)dt.\]
     \end{definition}
     This Banach manifold contains the solutions to the graph flow
     equation since on a non-compact (outgoing, say) edge $E$
     associated to the critical point $a_v$,
     \[\int_E\left|\frac{d\gamma}{dt}\right|^2dt=
     -\int_0^{\infty}\left\langle\frac{d\gamma}{dt}, \nabla
     f_E\right\rangle dt=-\int_0^{\infty}\frac{df_E}{dt}dt=
     f_E(\gamma(v))-f_E(a_v)<\infty.\] Then $\cp_{\g}(M;\vec{a})$ is
     defined as a union of $\cp_{\g_k}(M;\vec{a})$ in the same way that
     $\cp_{\g}(M)$ is defined.  Similarly define
     $\cp_{\g}(TM;\vec{a})$.  The graph flow map
     $\Phi:\cp_{\g}(M;\vec{a})\to\cp_{\g}(TM;\vec{a})$ defines the
     universal moduli space as its zero set:
     \[\cm_{\g}(M;\vec{a})=\Phi^{-1}(0)\subset\cp_{\g}(M;\vec{a}).\]

     Regularity of $\Phi$ at 0 requires the following minor  
adjustments to the
     compact case.  The proof of Theorem~\ref{th:reg} shows that any
     element $r$ of the cokernel of $D\Phi$ must vanish on $\g_k$, the
     compact part of $\gn_k$.  But on a non-compact edge $E_v$, $r$ is
     determined, via the codimension 1 condition (\ref{eq:cokerv}) at
     $v$, by its values on compact edges containing $v$ and hence it
     vanishes on $E_v$ and so vanishes everywhere on $\gn_k$.  Thus
     $D\Phi$ is onto.  It has a right inverse since $D_{\gn_k}$ is  
Fredholm  (see the appendix).

     There is no change to the proof of a canonical coorientation on $ 
\cm_{\g}(M;\vec{a})$.  As before, the virtual codimension of $\pi$  
follows from an index calculation.
     The projection $\pi$ to the parameter space is
     Fredholm, since $D_{\gn_k}$ is Fredholm, and
     ${\rm index\ }D\pi={\rm index\ }D_{\gn_k}$.
     \begin{lemma}
     \[{\rm index\ }D_{\gn_k}={\rm dim\ }M\cdot\chi(\g)
     -\sum_{v>0}\left( {\rm dim\ }M-{\rm index}(a_v)\right)
     - \sum_{v<0} {\rm index}(a_v).\]
     \end{lemma}
     \begin{proof}
     Choose trivialisations of $\cv=\gamma^*TM$ over $\gn_k$ so
     that transition functions are simply multiplication by $\pm
     1$.  With respect to these local trivialisations
     $D_{\gn_k}=\frac{d}{dt}+A(t)$ is well-defined since the
     operator commutes with multiplication by $\pm 1$.  The index
     remains unchanged under continuous deformations of $D_{\gn_k}$
     although we cannot deform $A(t)$ to zero as in the compact
     case, because at infinity $A(t)$ looks like the Hessian of the
     Morse function $f$ at each critical point and hence it is
     invertible there.  However, we may deform $A(t)$ so that it is
     diagonal, zero on $\g_k\subset\gn_k$ and constant outside of a
     compact subset of $\gn_k$ that contains $\g_k$.

     For $A={\rm diag}(\lambda_1(t),...,\lambda_d(t))$ we can
     explicitly solve the system for the kernel:
     \[\dot{s}_i=-\lambda_i(t)s(t),\ i=1,...,d\ .\] Since
     $\lambda_i(t)=\lambda_i^v$ is constant near infinity along
     $E_v\subset\gn$ then $s(t)\sim e^{-\lambda_it}$ near infinity.
     Thus, $s\in W^{1,2}(\br^+,\br)$ only when $\lambda_i^v<0$
     (respectively, $\lambda_i^v>0$) when $E_v$ is incoming
     (respectively, outgoing).  If the $i$th eigenvalue does not
     satisfy this condition for a single $v$ then the solution must
     vanish on $E_v$ and hence by continuity on all of $\gn$.  We
     see then that the dimension of the kernel is given by the
     number of $\lambda_i(t)$ with $\lambda_i^v<0$ for all $v$
     oriented positively ($E_v$ incoming) and $\lambda_i^v>0$ for
     all $v$ oriented negatively ($E_v$ outgoing.)

     For the cokernel we use $-A^T$ so negate each $\lambda_i^v$.
     Then $r_i\in W^{1,2}(\br^+,\br)$ when $\lambda_i^v>0$
     (respectively, $\lambda_i^v<0$) when $E_v$ is incoming
     (respectively, outgoing).  It is no longer true that if $r_i$
     vanishes along one edge then it vanishes on all of $\gn$.  For
     each $i$ we get a contribution to the cokernel from each
     incoming (outgoing) edge $E_v$ with $\lambda_i^v>0(<0)$.

     In order to calculate
     \[{\rm index}\ D_{\gn_k}=\dim\ker D_{\gn_k}-\dim\coker D_{\gn_k}\]
     change the index of a critical point and observe the change in
     ${\rm index}\ D_{\gn_k}$.  For an incoming edge $E_v$ change
     ${\rm index}(a_v)$ to ${\rm index}(a_v)-1$, so take
     $\lambda_i^v<0$ and send it to $-\lambda_i^v$.  Either
     $\lambda_i^v$ contributes to the kernel (it cannot contribute
     to the cokernel) then $-\lambda_i^v$ cannot contribute to the
     cokernel and we lose $1$ from ${\rm index\ }D_{\gn_k}$, or
     $\lambda_i^v$ does not contribute to the kernel in which case
     $-\lambda_i^j$ does contribute to the cokernel and we again
     lose $1$ from ${\rm index\ }D_{\gn_k}$.  A similar argument
     shows that on an outgoing edge $E_v$, the change ${\rm
     index}(a_v)\mapsto {\rm index}(a_v)+1$ affects the change
     ${\rm index\ }D_{\gn_k}\mapsto{\rm index\ }D_{\gn_k}-1$.  Thus
     $$ {\rm index}\ D_{\gn_k}=\sum_{v>0}{\rm index}(a_v)
     - \sum_{v<0} {\rm index}(a_v)+\ {\rm constant}.$$
     To determined the constant, suppose that ${\rm
     index}(a_v)={\rm dim\ }M$ (i.e. $\lambda_i^v<0$ for all $i$)
     for each incoming $E_v$ and ${\rm index}(a_v)=0$ (i.e.
     $\lambda_i^v>0$ for all $i$) for each outgoing $E_v$.  Then
     the non-compact edges make no contribution to the cokernel and
     there is no obstruction to the kernel.  Hence the index is the
     same as that for the compact graph, i.e. \[{\rm index\
     }D_{\gn_k}={\rm index\ }D_{\g_k}= {\rm dim\ }M\cdot\chi(\g)\]
     and the constant agrees with the statement of the lemma.  (In
     terms of the graph flow, we have just seen that when incoming
     and outgoing edges converge respectively to maxima and minima
     of $f$, locally it is as if there is no critical point
     restriction.)
     \end{proof}
     This completes the proof of the theorem.
\end{proof}

To the collection $\vec{a}$ of $l$ critical points of $f$ associate
the product of stable and unstable manifolds $W(\vec{a})\subset M^l$
\[ W(\vec{a})=\prod_{v>0}\cw^u(a_v)\times\prod_{v<0}\cw^s(a_v).\] Now  
consider the evaluation map on the univalent vertices,
$ev:\cm^N_{\g}(M) \to M^{p+q}$(we are assuming $p$ incoming leaves  
and $q$ outgoing leaves).   It is clear that
\[\cm^N_{\g}(M;\vec{a})=\cm^N_{\g}(M) \cap ev^{-1}(W(\vec{a})).\] In the
introduction to Section~{\ref{sec:trans}, we claimed that $N$ can be
chosen so that $ev(\cm^N_{\g}(M))$ intersects $W(\vec{a})$  
transversally.
The proof of this does not give a new proof that $\cm^N_{\g}(M;\vec 
{a})$ is a manifold since it uses the proof of that fact, although it  
is a more intuitive way of seeing the
manifold structure and its dimension and it will be used in the  
compactness arguments.

\begin{lemma}
     If $\gamma\in\cm^N_{\g}(M;\vec{a})$ is a regular point of the
   flow map $\Phi$,  then $ev(\cm^N_{\g}(M))$ intersects
     $W(\vec{a})$ transversally in $M^{p+q}$.
\end{lemma}
\begin{proof}
     If $ev(\cm^N_{\g}(M))$ does not intersect $W(\vec{a})$  
transversally
     at $\vec{x}\in M^{p+q}$ then there is a vector $\vec{\xi}\in
     T_{\vec{x}}M^{p+q}$ orthogonal to the tangent spaces of
     $ev(\cm^N_{\g}(M))$ and $W(\vec{a})$.  Take a non-zero component of
     $\vec{\xi}$ in one factor $M$ of $M^{p+q}$, corresponding to the
      univalent   vertex $v\in\g$.  Along the non-compact edge $E_v$
     parametrised by $t\in[0,\infty)$ solve the equation
     $\dot{r}(t)-(\nabla\nabla f)^T\cdot r(t)=0$ with $r(0)=\xi$.
     Since $\xi$ is orthogonal to $ev(\cm^N_{\g}(M))$ at $x$, $r(t)$
     decays at infinity and lives in $L^2$.  Put $r=0$ on the rest of
     the graph $\gn_k$.

     Since $D_{\g_k}$ is surjective, $r\in{\rm im\ } D_{\g_k}$ so
     $r=D_{\gn_k}s$ for some $s\in W^{1,2}(\gn_k)$.  Now
     \[\int_{\gn_k}\langle r,r \rangle dt=\int_0^{\infty}\langle r,
     D_{\gn_k}s\rangle dt=\int_0^{\infty}d/dt\langle r,s\rangle dt
     =\langle \xi,s(0)\rangle.\]
     But $r\equiv 0$ on $\g_k$ so $D_{\g_k}s=0$ so
     $s(0)\in T_x ev(\cm^N_{\g}(M))$ and $\langle \xi,s(0)\rangle=0$  
which
     is a contradiction.
\end{proof}

The following corollary is a generalisation of the Morse-Smale
condition.
\begin{corollary}   \label{th:msmale}
Any submanifold   $N\subset\cm_{\g}(M)$ can be perturbed   so that $ev 
(\cm^N_{\g}(M))$ intersects
     $W(\vec{a})$ transversally {\em for all} collections of critical
     points $\vec{a}$ of $f$.
\end{corollary}
\begin{proof}
     For each collection of critical points $\vec{a}$ of $f$, the proof
     of Theorem~\ref{th:moduli} supplies a universal moduli space
     together with a map to the parameter space
     $\pi:\cm_{\g}(M;\vec{a})\to\cs_{\g}$.  For a given $\vec{a}$, the
     Sard-Smale theorem allows one to make an arbitrarily small
     deformation of a submanifold with boundary $N_0\subset\cm_{\g}$ to
     $N_1$ that is transverse to $\pi(\cm_{\g}(M;\vec{a}))$.  Take
     another collection $\vec{a'}$ and again apply the Sard-Smale
     theorem to choose a deformation $N_2$ of $N_1$ small enough so
     that it remains transverse to $\pi(\cm_{\g}(M;\vec{a}))$ and so
     that it is also transverse to $\pi(\cm_{\g}(M;\vec{a'}))$.  Take
     the finite list of all collections of critical points $\vec{a}$
labeled by a given set of  univalent   vertices of $\g$, and update
     $N_0, N_1, N_2,\ldots$ to get a finite sequence that finishes at
     $N\subset\cm_{\g}$ simultaneously transverse to all the spaces,
     $\pi(\cm_{\g}(M;\vec{a}))$.
\end{proof}

\section{Compactness.}  \label{sec:compact}
The graph moduli spaces are non-compact due to the non-compact edges
of the graph.  This will imply, as we will see, the non-compactness  
and gluing issues essentially
reduce to these same issues for  spaces of gradient flows of a Morse  
function.

\subsection{Piecewise graph flows.}
   We begin by recalling  the natural compactification of the space  
of gradient flow
lines converging to two fixed critical points of the Morse function
$f$.

The space of flow-lines of the Morse function $f$ from the critical
point $a$ to the critical point $b$  can be viewed as using  the  
noncompact graph $\g=\br$ which has
a one-dimensional space of translational symmetries.  The moduli space
of flows is the quotient space
\[\cm(a,b)=\cm_{\br}(M;a,b)/\br.\]

Notice that  $\cm_{\br}(M;a,b)$ is the intersection of the unstable  
manifold of $a$ with the stable manifold of $b$,   $\cm_{\br}(M;a,b)  
= W^u_a \cap W^s_b$.

Now assume that $M$ is equipped with a metric so that $f: M \to \br $  
satisfies the
Morse-Smale condition.  This says that the intersections of stable  
and unstable manifolds are all transverse.   Recall the partial  
ordering on the set of critical points in this setting,
$a\geq b$ if there is a gradient flow connecting $a$ and $b$, i.e $\cm 
(a,b) \neq \emptyset$.

Define the space of {\em piecewise flow lines} connecting critical
points $a$ and $b$ by:
\[\overline{\cm}(a,b) = \bigcup_{a=a_0 \geq a_1 \geq ... \geq a_j = b}
\cm(a,a_1) \times \cm(a_1,a_2) \times ...  \times \cm(a_{j-1},b)\]
where the union is taken over all  nonincreasing finite sequences of
critical points.  For  example, $a\geq b$  implies
$\cm(a,b)\subset\overline{\cm}(a,b)$.

Since
$f$ satisfies the Morse-Smale condition, $a > b$ implies that
$f(a) > f(b)$ and ${\rm index\ }a > {\rm index\ }b$.
The result is that $\overline{\cm}(a,b)$ is compact, which is a simple
equicontinuity argument, and that it contains $\cm(a,b)$ as an open
dense subset.  This    is often expressed as a gluing theorem since it
implies the existence of true flows arbitrarily close to piecewise
flows.  A uniqueness part of gluing further implies that
$\overline{\cm}(a,b)$ is a manifold with corners.  For our purposes
it is sufficient to consider at most 1-dimensional moduli spaces.  In  
the one dimensional
$\overline{\cm}(a,b)$ is a 1-manifold with boundary.  In
particular a
deleted neighbourhood of any  boundary component in $\overline{\cm} 
(a,b)$
is a connected, open interval.

By analogy, we define the space of piecewise graph flows by
\[\overline\cm_{\g}^N(M; \vec{a}) =
\bigcup_{\vec{b} } \cm_{\g}^N(M; \vec{b})  \times
\prod_{v\, \, incoming }\overline\cm(a_v,b_v) \times
\prod_{v \, \, outgoing }\overline\cm(b_v,a_v)\]
where the union is taken over all collections of critical points
$\vec{b}$ labeling  the  univalent   vertices of $\g$.  Notice that  
the restriction of such piecewise
graph flow $\gamma$ to a compact edge is a gradient flow of the  
function labeling that edge, and when restricted to a noncompact  
edge, it is a piecewise flow line.

\begin{proposition}
     When $N$ is compact, $\overline\cm_{\g}^N(M;\vec{a})$ is compact.
\end{proposition}
\begin{proof}
     For any $(\g_k,\vec{f})\in N$, the gradient vector fields $\nabla
     f_E$ along the edge $E\subset\gn_k$ are bounded and uniformly
     continuous, uniformly in $N$, since $M$ is compact.  (As usual, we
     express the Morse function $f$ by $f_E$ for any non-compact edge
     of $\gn_k$.)

     Hence the space of maps $\overline\cm_{\g}^N(M;\vec{a})$ is an
     equicontinuous family since the derivatives $d\gamma_E/dt=-\nabla
     f_E$ are uniformly bounded.  Let
     $\{\gamma^j\}\subset\overline\cm_{\g}^N(M;\vec{a})$ be a sequence
     of piecewise graph flows.  Take any  univalent   vertex $v\in\gn 
$, or
     any point on a non-compact edge labeled by its parameter $T<0$
     ($T>0$) for an incoming (outgoing) edge $E$ of the metric graph.
     Both of these give well-defined choices of points in any metric
     graph in $N$.  Since $M$ is compact, the sequence $\gamma^j(v)$,
     or $\gamma^j(T)$, has a convergent subsequence converging to a
     point $x\in M$.  By differentiating $\nabla f_E$ over $M$ one gets
     a uniform $C^2$ bound on the $\{\gamma^j\}$ and thus the limit of
     the subsequence satisfies the flow equation.  Thus the flow from
     the limit point $x$ is a uniform limit of the subsequence of graph
     flows.  It may be a graph flow or a gradient flow of $f$.  As we
     choose different points on non-compact edges, we get different
     gradient flows of $f$ that are also uniform limits of   a   
subsequence of graph flows.

     So the limit of a sequence of piecewise flows is locally a flow
     and to prove that it is itself a piecewise flow it remains to show
     that the limit is a continuous map from $\gn_k$ to $M$.
     Canonically parametrise $\{\gamma^j\}$ by $s=f(\gamma^j(t))$ so
     they satisfy $d\gamma^j(s)/ds+\nabla f/|\nabla f|=0$.  Again one
     gets a uniform bound on $d\gamma^j(s)/ds$ so by equicontinuity the
     limit of the subsequence is a continuous map from $\gn_k$ to $M$
     and hence a piecewise flow.
\end{proof}
{\em Remark.} In the above proof   it  is clear that the
   that non-compactness of the moduli space of graph flows     arises  
due to the non-compact edges of the graph.
    We say that a sequence {\em bubbles} along a
non-compact edge if its limit is not a smooth flow there.

\begin{corollary}  \label{th:zerod}
     For generic choice of $N$, if ${\rm dim\ }\cm_{\g}^N(M;\vec{a})=0$
     then $\cm_{\g}^N(M;\vec{a})$ is compact.
\end{corollary}
\begin{proof}
     Choose $f$ to be Morse-Smale and $N$ as in
     Corollary~\ref{th:msmale} so that all moduli spaces $\cm_{\g}^N(M;
     \vec{b})$ are manifolds of the expected dimension.  Suppose that a
     sequence of graph flows bubbles along an incoming edge $E_v$ and
     converges to a piecewise graph flow.  Since $f$ is Morse-Smale,
     $\cm(a_v,b_v)$ is non-empty only if ${\rm index\ }a_v>{\rm index\
     }b_v$.  But then \[{\rm dim\ }\cm_{\g}^N(M;\vec{b})<{\rm dim\
     }\cm_{\g}^N(M;\vec{a})=0\] so by transversality
     $\cm_{\g}^N(M;\vec{b})$ is empty, contradicting the claim that the
     sequence bubbles.  The same argument works for an outgoing edge.
     Thus no bubbling can occur and $\overline\cm_{\g}^N(M;\vec{a})=
     \cm_{\g}^N(M;\vec{a})$.
\end{proof}

\med

\begin{theorem}   \label{th:glue}
     For generic choice of $N$, if ${\rm dim\ }\cm_{\g}^N(M;\vec{a})=1$
     then $\overline\cm_{\g}^N(M;\vec{a})$ is a 1-manifold with
     boundary     $\displaystyle\bigcup_v\bigcup_{b_v}
     \cm_{\g}^N(M;\vec{b})\times\cm(a_v,b_v)$
     where for each $v$ ${\rm index\ }b_v={\rm index\ }a_v\pm 1$,
     and it contains $\cm_{\g}^N(M;\vec{a})$ as an open dense
     subset.
\end{theorem}

\med
\begin{proof}
     The same argument as in the proof of Corollary~\ref{th:zerod}
     shows that for a 1-dimensional moduli space
     $\cm_{\g}^N(M;\vec{a})$, any sequence
     $\{\gamma^j\}\subset\cm_{\g}^N(M;\vec{a})$ bubbles at most once.
     If a sequence bubbles along the incoming edge $E_v$ then its limit
     is given by the pair $(\gamma,\mu)$ satisfying
     \begin{enumerate}
     \item[(i)] $\gamma\in\cm_{\g}^N(M;\vec{b})$,
     \item[(ii)] $\mu\in\cm(a_v,b_v)$ uniquely defined up to rescaling,
     \item[(iii)] ${\rm index\ }b_v={\rm index\ }a_v-1$ so ${\rm dim\
     }\cm_{\g}^N(M;\vec{b})=0$.
     \end{enumerate}
     {\em Conversely}, to prove the theorem we need to show that any
     $(\gamma,\mu)$ satisfying (i), (ii) and (iii) is a unique end of
     $\cm_{\g}^N(M;\vec{a})$.  The same argument will apply to an
     outgoing edge.

     We follow the approach in \cite{BFKEle}.  The idea is to find a
     manifold with boundary $\cp$ and a smooth manifold $\cn$ that lie
     inside a common ambient space, such that the broken flow
     $(\gamma,\mu)$ maps to a point in both these manifolds.  If $\cp$
     and $\partial\cp$ intersect $\cn$ transversely then $(\gamma,\mu)$
     is a unique end of the 1-dimensional intersection $\cp\cap\cn$.
     More is proven in \cite{BFKEle} for higher-dimensional moduli
     spaces, where $\cp$ is a product of manifolds with boundary, so a
     manifold with corners, hence the transversal intersection inherits
     a structure of a manifold with corners.

     Put $f(b_v)=c$.  Choose $\epsilon>0$ small enough so that $c$ is
     the only critical value in $[c-\epsilon,c+\epsilon]$.  Define
     \[M^{\pm}=f^{-1}(c\pm\epsilon)\subset M\]
     and
     \[\cp\subset M^+\times M^-\]
     by pairs $(x^+,x^-)$ that flow to the same point $x\in f^{-1}(c)$
     under the forward, respectively backward, gradient flow (possibly
     flowing for infinite time.)

     Let $\vec{a}(-v)$ be $\vec{a}$ with $a_v$ removed.  Define $W_v^s$
     to be all those points of $M$ that flow under the gradient flow of
     $f$ to $ev_v(\cm_{\g}^N(M;\vec{a}(-v))$ and
     \[\cn=W^u_{a_v}\cap M^+\times W_v^s\cap M^-\subset M^+\times
     M^-.\]
     The ``stable manifold'' $W_v^s$ is a manifold of dimension $d-{\rm
     index\ }a_v+2$ for $d={\rm dim}\ M$ so $\cn$ is a $d$ dimensional
     manifold.

     It is proven in \cite{BFKEle} that $\cp$ is a $d-1$ dimensional
     manifold with boundary that intersects $\cn$ transversally inside
     the $2(d-1)$ dimensional manifold $M^+\times M^-$.  The critical
     point $b_v$ is contained inside the intersection $\cp\cap\cn$ and
     a neighbourhood of $b_v$ in $\cp\cap\cn$ is a 1-manifold $K$ with
     boundary $b_v$.

     The arguments in \cite{BFKEle} require the Morse function $f$ to
     be Morse-Smale, and we must choose either a metric on $M$ that is
     standard near critical points of $f$, or replace the gradient flow
     with a Morse-like vector field on $M$.  If we choose the latter,
     the analysis in Section~\ref{sec:zerod} does not change since it
     depends only on the fact that $\nabla\nabla f$ is invertible at
     infinity and this is still true of Morse-like vector fields.
     Thus, in our adaption of the arguments in \cite{BFKEle} we will
     require the same conditions on $f$ and replace the gradient
     vector field on external edges by a Morse-like vector field.

     Finally, we will prove that the analogues of stable and unstable
     manifolds for a graph flow intersect stable and unstable manifolds
     of $f$ transversally.  This is a slight adjustment of
     Corollary~\ref{th:msmale} which shows that the image of the moduli
     space of graph flows under the evaluation map, $ev(\cm^N_{\g})$,
     intersects the stable and unstable manifolds of $f$ transversally.

     First notice that the stable manifold for a graph flow, $W_v^s$,
     constructed from $ev_v(\cm_{\g}^N(M;\vec{a}(-v))$, is a moduli
     space of graph flows as follows.  For $\sigma=(\g_k,\vec{f})\in N$
     define $\sigma^+=(\g_k^+,\vec{f})$ on the graph $\g_k^+$ obtained
     from $\g_k$ by adding a compact edge $E$ at $v\in\g_k$ oriented
     inwards and assigning to $E$ the vector field $\nabla f$ and
     length $\ell_E$ any positive real number.  This gives a family of
     structures $N^+$ with $\dim N^+=\dim N+1$.  (Since $\g^+\to\g$ is
     a homotopy equivalence the set $N^+$ is almost a subset of
     $\cm_\g$ except that the lengths of edges do not add to 1.)

     The argument in Corollary~\ref{th:msmale} also shows that for any
     length $\ell_E$ on the extra compact edge $E\subset\g_k^+$, for
     generic choice of $N$ the image of the moduli space of graph flows
     under the evaluation map at the univalent vertex of $E$ intersects
     the unstable manifolds of $f$ transversally.  Thus, as we vary
     $\ell_E$ transversality is unchanged so $W_v^s$ intersects the
     unstable manifolds of $f$ transversally.  Note that
     $N\subset\cs_{\g}$ is chosen so that all moduli spaces
     $\cm_{\g}^N(M; \vec{b})$ are manifolds of the expected dimension
     which is {\em independent} of $\ell_E$.

     The same construction works for a negatively oriented vertex $v$
     by adding an outward pointing compact edge at $v$ to get $\g_k^-$
     and thus showing that $W_v^u$ intersects the stable manifolds of
     $f$ transversally.
\end{proof}

{\em Remark.}The main ingredient in gluing is the transversality
of the intersection of the image of the evaluation map and stable and
unstable manifolds of $f$, which follows from surjectivity of
$D\Phi_{\gn_k}$.  Gluing can be defined directly from surjectivity of
$D\Phi_{\gn_k}$.  One uses the energy functional defined on
$\cp_{\g_k}(M;\vec{a})$
\begin{eqnarray*}
     {\mathcal
     E}(\gamma)&=&\frac{1}{2}\int_{\gn_k}\left(\left|
     \frac{d\gamma}{dt}\right|^2+
     |\nabla f(\gamma)|^2\right)dt\\&=&f(\alpha)-f(\beta)+
     \frac{1}{2}\int_{\gn_k}\left|\frac{d\gamma}{dt}+\nabla
     f(\gamma)\right|^2dt
\end{eqnarray*}
where the first expression shows that $\mathcal E$ is non-negative and
the second expression shows that $\mathcal E$ is minimised by graph
flows.  A broken flow yields a path with small energy---an approximate
flow.  The implicit function theorem shows that there is a unique true
flow nearby.  Details for the case of the Morse complex can be found
in \cite{SchMor}.\\

Using the same gluing constructions as in the proof of
Theorem~\ref{th:glue} we will now show how to remove an edge
$E\subset\g$ leaving two marked vertices given by its endpoints.  (An
endpoint of $E$ must not coincide with an existing marked vertex of
$\g$.  To find such an edge it may be necessary to take an edge
$E\subset\g_1\to\g$ and consider all $\g_k\to\g_1\to\g$.)  The edge
$E$ may or may not be separating.  We denote $\g-E$ to be the graph,
or union of two graphs, with marked vertices those of $\g$ and the
endpoints of $E$, oriented according to the orientation of $E$.

Choose $N$ so that $\cm_{\g}^N(M;\vec{a})$ is a smooth
zero-dimensional moduli space for all $\vec{a}$ and so that for a
given edge $E\subset\g$ the induced structure on $\g-E$ gives a smooth
zero-dimensional moduli space for all $\vec{b}$ .  Here
$\vec{b}=(\vec{a},a^-,a^+)$ is a vector of critical points of $f$
associated to the marked vertices of $\g-E$, and $a^-=a^+$ is named
twice because it is used twice.  The critical point $a^-$ is
associated to the negatively oriented (outgoing) endpoint of $E$ and
$a^+$ is associated to the positively oriented (incoming) endpoint of
$E$.  To the edge $E$, each metric-Morse structure in $N$ should
associate the gradient vector field $\nabla f$ of the external Morse
function.

\med
\begin{theorem}  \label{th:cobord}
     The moduli spaces $\cm_{\g}^N(M;\vec{a})$ and
     $\displaystyle\bigcup_{(a^-,a^+)}
     \cm_{\g-E}^N(M;\vec{a},a^-,a^+)$ are cobordant.
\end{theorem}

\med
\begin{proof}
     Define the one-dimensional moduli space
     $\cm_{\g}^{N_E}(M;\vec{a})$ using a family $N_E$ of structures
     with $\dim N_E=\dim N+1$ as follows.  For
     $\sigma=(\g_k,\vec{f})\in N$, take the edge
     $E_k=\phi^{-1}(E)\in\g_k$ where $\phi:\g_k\to\g$ is the homotopy
     equivalence and assign to it any length
     $\ell\in[\ell_{E_k},\infty)$.  This gives a family $\tilde{N}_E$
     of structures on $\g$ with $\partial \tilde{N}_E=N\cup N|_{\g-E}$.
     As in the proof of Theorem~\ref{th:glue} the family $N_E$ is not
     contained in $\cm_\g$ since the lengths of edges do not add to 1
     so we use an enlargement of $\cs_{\g}$ to allow $E$ to have an
     arbitrarily large edge length.  Inside this space of parameters
     take an arbitrarily small deformation $N_E$ of $\tilde{N}_E$ that
     fixes the boundary so that $N_E$ is transversal to
     $\pi(\cgm(M;\vec{a}))$ for all $\vec{a}$.  Then it immediately
     follows that $\cm_{\g}^{N_E}(M;\vec{a})$ is a one-dimensional
     manifold with compact and non-compact ends.  At the compact ends
     it is a manifold with boundary $\cm_{\g}^N(M;\vec{a})$ and we will
     show that it can be compactified at the non-compact ends so that
     the 1-manifold gives the cobordance stated in the theorem.  In
     other words
     \begin{equation}  \label{eq:boundary}
     \partial\overline\cm_{\g}^{N_E}(M;\vec{a})=
     \cm_{\g}^N(M;\vec{a})\cup \displaystyle\bigcup_{(a^-,a^+)}
     \cm_{\g-E}^N(M;\vec{a},a^-,a^+).
     \end{equation}
     The same transversality argument as in the proof of
     Corollary~\ref{th:zerod} shows that any sequence
     $\{\gamma^j\}\subset\cm_{\g}^N(M;\vec{a})$ bubbles at most once
     along the edge $E$ to give a graph flow in
     $\cm_{\g-E}^N(M;\vec{a},a^-,a^+)$ for critical point $a^-=a^+$
     with index so that $\dim\cm_{\g-E}^N(M;\vec{a},a^-,a^+)=0$.  (The
     expected dimension is the same as the actual dimension.)  As
     usual, if $\g-E$ is disconnected then
     $\cm_{\g-E}^N(M;\vec{a},a^-,a^+)$ is the product of moduli spaces
     for each component of $\g-E$ and by a graph flow in
     $\cm_{\g-E}^N(M;\vec{a},a^-,a^+)$ we mean a pair of graph flows.

     The theorem will be proven if we can show that for any flow in the
     zero-dimensional moduli space $\cm_{\g-E}^N(M;\vec{a},a^-,a^+)$
     there is a unique flow nearby in $\cm_{\g}^{N_E}(M;\vec{a})$.
     This gluing result follows the proof of Theorem~\ref{th:glue}
     exactly.  Once again we construct a manifold with boundary $\cp$
     and a smooth manifold $\cn$ that lie inside a common ambient
     space, such that a broken graph flow given by a flow in
     $\cm_{\g-E}^N(M;\vec{a},a^-,a^+)$ maps to a point in both these
     manifolds.  In fact
     $\overline\cm_{\g}^{N_E}(M;\vec{a})=\cp\cap\cn$ and the
     intersection will be transverse so
     $\overline\cm_{\g}^{N_E}(M;\vec{a})$ is a 1-manifold with boundary
     and in particular any broken flow is a unique end of this
     1-manifold.

     Put $f(a^{\pm})=c$ then there is no change to the definition of
     $\cp\subset M^+\times M^-$ for $M^{\pm}=f^{-1}(c\pm\epsilon)$.  In
     the definition of $\cn$ we now use stable and unstable manifolds
     of graph flows:
     \[\cn=W_{v^-}^u\cap M^+\times W_{v^+}^s\cap M^-\subset M^+\times
     M^-\]
     where $v^{\pm}$ are the endpoints of $E$ and $W_{v^-}^u$ and
     $W_{v^+}^s$ are defined in the proof of Theorem~\ref{th:glue}.
     Arguing as in Corollary~\ref{th:msmale} it can be shown that when
     $N_E$ is chosen transversally to $\pi(\cgm(M;\vec{a}))$ the stable
     and manifolds $W_{v^-}^u$ and $W_{v^+}^s$ intersect transversally
     so the theorem follows.
\end{proof}

{\em Remark.}  In the previous two theorems, if the moduli spaces are
oriented then the orientation on the 1-manifold agrees with the
orientations on the boundary.  This is because the orientations are
canonically induced from the evaluation map, and the gluing
construction also used the evaluation map.

\section{Cohomology operations on the Morse chain complex.}

In this section we represent the homology $H_*(M)$ in terms of the
Morse complex of the Morse function $f:M\to\br$ and express the
homology operation $$q_{\g}:H_*(BAut_0(\g))\otimes H_*(M)^{\otimes
p}\to H_*(M)^{\otimes q}$$ with respect to this representation.

\med
Recall that the Morse complex of a Morse function $f$ is the chain  
complex of abelian
groups
\[ C_n\stackrel{\partial}\to C_{n-1}\stackrel{\partial}\to\ldots
\stackrel{\partial}\to C_1\stackrel{\partial}\to C_0 \]
generated by the critical points of $f$, graded by their index.  The
boundary operator $\partial$ is defined by counting points in the
moduli space of solutions to the gradient flow equation converging to
critical points of consecutive degrees.  We will give the proof that
that this does indeed define a complex, i.e.
$\partial\circ\partial=0$, since an analogous proof is used to show
that the graph moduli spaces define homological invariants.

Let $a$ and $b$ be critical points of $f$ of index $k+1$ and $k$
respectively.  It follows from the analysis in Section~\ref{sec:trans}
that $\cm(a,b)$ is a zero-dimensional oriented compact manifold.  Thus
it makes sense to count the points, with sign, in $\cm(a,b)$.  Put
$n(a,b)=\#\cm(a,b)$ and define the linear operator
\[\partial a=\Sigma n(a,b)b\]
where the sum is over all critical points $b$ of index $k$.
\begin{lemma}
\[\partial^2=0\ .\]
\end{lemma}
\begin{proof}
     By linearity \[\partial\partial a=\Sigma n(a,b)\partial b=\Sigma
     n(a,b)n(b,c)c\] where the sum is over all critical points $b$ of
     index $k$ and $c$ of index $k-1$.  We will show that for fixed $c$
     the sum $\Sigma n(a,b)n(b,c)c$ over all intermediate critical
     points $b$ vanishes.  By Theorem~\ref{th:glue} the compactified
     one-dimensional moduli space $\overline{\cm(a,c)}$ is a manifold
     with boundary.  That is the boundary points, which are piecewise
     flows, each correspond to a unique edge.  Since one-dimensional
     compact manifolds can only be a finite collection of closed
     intervals this means that the ends come in pairs.  Thus the
     contributions to $\partial^2(a)$ come in pairs.  This immediately
     gives the vanishing of each component modulo two.  Furthermore the
     orientations on the 1 dimensional moduli space and its boundary
     agree, meaning that $n(a,b)n(b,c)=-1(+1)$ if that boundary
     component is oriented negatively (positively).  This is because
     the orientations are defined canonically using the evaluation map
     and the gluing construction also uses the evaluation map.  Thus
     the boundary points are oriented oppositely so the oriented sum
     vanishes.
\end{proof}

Choose an $\aog$-invariant submanifold $\tilde N\subset\cs_{\g}$ such  
that the quotient $N = \tilde N/\aog \in \cm_\G$ is transverse to the
image of the universal moduli space.  Given the Morse-Smale
function $f$, let $C_*(M,f)$ be the associated Morse-Smale chain complex
generated by the critical points, and let $C^*(M,f)$ be the dual
cochain complex.  The cochains are negatively graded so that the  
evaluation pairing
$C^*(M,f) \otimes C_*(M,f) \to \bz$ is of degree zero.

Define  a class  $q_{\g}^N$ to be an element of the tensor product  
complex,
$$
\bigotimes_{v \, \, incoming} C^*(M,f)  \bigotimes_{v\, \, outgoing }  
C_*(M,f)
$$
in the following manner.
Consider those collections of critical points $\vec{a}$
such that ${\rm dim}\, \msa = 0$.
These spaces contain a finite number of
oriented points which can be
counted with sign (if $\msa$ is oriented---otherwise this is well  
defined mod
$2$, and we take  coefficients to
be  ${\bf Z}_2$).

\begin{definition}
$$
q_{\g}^N= \sum \#\msan [\vec{a}] \in \bigotimes_{v \, \, incoming} C^* 
(M,f)
\bigotimes_{v\, \, outgoing } C_*(M,f).
$$
\end{definition}

Theorem~\ref{th:glue} and the definition of the boundary operator in
the Morse-Smale complex yields the following.

\begin{lemma}  $$dq = 0.$$
\end{lemma}
\begin{proof}  Recall that the boundary operator on the tensor  
product of the chain complexes  is given by
\[\partial:\bigotimes_{1 \leq i \leq k} C_*(M,f)\to\bigotimes_
{1 \leq i \leq k} C_*(M,f)\]
\[(a_1,...,a_n)\mapsto\Sigma_i(a_1,...,\partial_i(a_i),...,a_k)\]
where $\partial_i$ is defined using $f_i$.  Then if we think of $q$  
as a map
\[ q:\bigotimes_{1 \leq i \leq n_1} C_*(M,f)\to
\bigotimes_{n_1+1 \leq i \leq n} C_*(M,f)\ ,\] the requirement
that $dq=0$ is equivalent to the requirement that $q$  is a chain map:
   $\partial q=q\partial$.  Choose $\vec{a}=(\vec{b},\vec{c})$
so that $\dim\ \cm^N_{\G}(M; \vec {a})=1$.  We have divided $\vec{a}$  
into critical points
$\vec{b}$ corresponding to incoming flows and $\vec{c}$ corresponding
to outgoing flows.  Notice that for $\partial\vec{b}=\Sigma\vec{b}^j 
$, then
$\dim\cm_{\g}^N( M;(\partial\vec{b},\vec{c}))=0$ so
$q(\partial\vec{b})\in\bigotimes_{v<0} C_*(M,f)$ is obtained by counting piecewise graph flows,
containing a piecewise gradient flow along an incoming edge, from $\vec{b}$ to $\vec{c}$, and it takes it values in the module generated by $\vec{c}$.  The composition $\partial q(\vec{b})\in\bigotimes_{v<0} C_*(M,f)$ is given by piecewise graph flows, containing a piecewise gradient flow along an outgoing edge, and takes its values in the same module generated by $\vec{c}$ so it makes
sense to compare $q(\partial\vec{b})$ and $\partial q(\vec{b})$.  We will show that there is a pairing between the two types of piecewise graph flows which gives $\partial q=q\partial$.

The one-dimensional manifold $\msa$ is compact with boundary so it is
a finite collection of closed intervals.  Each boundary point of an
interval corresponds to a piecewise graph flow with exactly one external
edge not a true gradient flow.  This is the key fact behind the proof.
If more than one external edge were to break then the true graph flow
inside this piecewise graph flow would lie in a moduli space of negative
dimension, thus contradicting its existence.  These boundary piecewise
graph flows are paired by the interval they bound.  

There are three types of components of the one-dimensional manifold $\msa$ and thus  
three types of pairings of piecewise flows.  The first type of component
consists of an interval whose two boundary points correspond to 
piecewise gradient flows both containing a piecewise gradient flow along an incoming edge.  The sign, or orientation, given to the piecewise flow is the product of the signs, or orientations, given to the two components of the piecewise flow.  But this is the orientation induced from the one-dimensional moduli space.  Since the two boundary components of the one-dimensional moduli are oriented oppositely - they are two ends of an oriented interval - the two piecewise graph flows contribute a total of $1-1=0$ to 
$q(\partial\vec{a})$.

The second type of component consists
of an interval whose two boundary points correspond to two piecewise gradient flows both containing a piecewise gradient flow along an incoming edge.  It behaves like the first type of component and the two piecewise graph flows contribute $1-1=0$ to $\partial q(\vec{a})$.  

The third type of component consists
of an interval whose two boundary points correspond to two piecewise gradient flows containing, respectively, a piecewise gradient flow along an incoming edge and a piecewise gradient flow along an outgoing edge.   The one-dimensional moduli space gives an oriented cobordism between the two piecewise graph flows, so they contribute, respectively, $q(\partial\vec{a})$ and to $\partial q(\vec{a})$ with the same sign.  

We pair piecewise flows arising from the third type of component and cancel pairs of piecewise flows arising from the other two types of components to get $q(\partial\vec{a})=\partial q(\vec{a})$
and the lemma is proven. \end{proof}

We shall therefore view $q_{\g}^N$ as an element of the associated  
homology,
$$
q_{\g}^N\in H^*(M)^{\otimes n_1} \otimes H_*(M)^{\otimes n_2}.
$$
In fact $q_{\g}^N$ is independent of the choice of $N\subset\cs_{\g}$
and only depends on the homology class of $N$.  We prove this in the  
following
proposition.

\begin{proposition}   \label{th:hom}
If $N_1$ and $N_2$ are homologous, then $q_{\g}^{N_1}=q_{\g}^{N_2}$.
\end{proposition}
\begin{proof}
     A cobordism between $N_1$ and $N_2$ produces a non-compact
     1-dimensional moduli space with boundary.  Its compactification has
     boundary components consisting of the moduli spaces associated to
     $N_1$ and $N_2$ and to broken flows which correspond to
     compositions with the boundary operator.  Thus the compactified
     1-dimensional moduli space defines a chain homotopy equivalence
     between the invariants so on the level of homology
     $q_{\g}^{N_1}(M)=q_{\g}^{N_2}(M)$.
\end{proof}

It is easy to see that $q_{\g}^N$ coincides with the algebraic
topology version of the invariant defined in section 2.  This is  
because of the standard relationship between umkehr maps and  
intersection theory of chains.  \\

We end this section by giving an analytic version of the gluing
construction in Section 4.

Let $\Gamma_1$ and $\Gamma_2$ be oriented graphs.  Let
$\Gamma^{i\#j}_{1,2}$ be the oriented graph obtained by gluing  
incoming edge
$i$ of $\Gamma_1$ to outgoing edge $j$ of $\Gamma_2$.

\med
\begin{proposition}  \label{th:glu}
     $$
     q(\Gamma^{i\#j}_{1,2},M) = q(\Gamma_1, M) \Diamond^{i,j} q 
(\Gamma_2,M),
     $$
     where $\Diamond^{i,j}$ denotes tensorial contraction of cohomology
     in the $i$th coordinate with homology in the $j$th coordinate.
\end{proposition}

\med
\begin{proof}
     This uses Theorem~\ref{th:cobord} repeatedly to glue together any
     number of edges between $\g_1$ and $\g_2$.  As in the proof of
     Proposition~\ref{th:hom}, the compactified 1-dimensional moduli  
spaces
     have boundary components consisting of components of the
     zero-dimensional moduli spaces and broken flows so this gives a
     chain homotopy equivalence between the invariants.
\end{proof}

\med
\begin{corollary}\label{change}
     Changing the orientation of a non-compact edge induces the
     Poincare duality isomorphism on the relevent tensor coordinate of
     the invariant $q_{\g}(M)$.
\end{corollary}

\med
\begin{proof}
     Let $\Gamma$ be a given graph with outgoing edge $E$.  Recall  
the graph with two incoming univalent vertices discussed in section 5  
above.  It is   pictured in Figure 10.  Glue  this    graph  to $ 
\Gamma$ at $E$ to get a graph we'll call  $\Gamma'$.  By
     Proposition~\ref{th:glu} $q_{\Gamma'}$ is the composition of
     $q_{\Gamma}$ with the Poincare duality isomorphism.  Contract the
     internal glued edge to a point.  By Proposition~\ref{th:hom} this
     does not change the invariant.
\end{proof}

One can use the contractible graph with one incoming vertex and one
outgoing vertex to get a  chain  homotopy between the Morse complexes
of different Morse
functions.  By gluing this graph onto the external edges of
any other graph using Proposition~\ref{th:glu}, one sees that the
cohomology operations do not depend on the choice of external Morse
function.  This also follows from the definition of the
invariants in Section 3.

\section{Appendix:  Proof of theorem \ref{tube}.}
In this section we give a proof of theorem \ref{tube}.  Let $M$ be a  
closed $n$-dimensional manifold with a fixed Riemannian metric.  We  
begin by describing an extension of  the bundle $\bigoplus_b TM \to M 
$ to an $\ag$-equivariant bundle over
$\tilde \cm_{tree}(\G,M) \cong   \cs_\G (M) \times M$.

Let $(\sigma, \gamma) \in \tilde \cm_{tree}(\G,M)$.   Let $T \subset  
\G$ be a fixed maximal subtree, and let $T_1 \subset \Gamma_k$ be the  
inverse image of $T$ under the composite morphism $\phi_k : \G_k \to  
\G_{k-1} \to \cdots \to \G_0 \to \G$ determined
by the structure $\sigma$.   Write $p_{T_1}(\sigma, \gamma) = ((x_1,  
y_1), \cdots (x_b, y_b)) \in (M^2)^b$ as in
(\ref{peetee}).  Recall that  $x_i = \gamma_{T_1}(s^T_i(\G_k)) \in M. 
$   We  define a vector bundle.
$$
\zeta \to \tilde \cm_{tree}(G,M)
$$
to have as its fiber over $(\sigma, \gamma)$ the sum of the tangent  
spaces,
$$
\zeta_{(\sigma, \gamma)} = \bigoplus_{i=1}^b T_{x_i}M
$$

It is clear that the bundle $\zeta$ is $\ag$-equivariant.  This is  
because  if  $g \in \ag$,  the action of $g$ on the element $(\sigma,  
\gamma) \in  \tilde \cm_{tree}(\G,M)$, is given by $(g\sigma, \gamma) 
$, where the structure $g\sigma$ is determined by the morphism    $g 
\phi_k : \G_k \to \G$  given by  the composition of the morphism $ 
\phi_k$ with the automorphism $g$.  $g\phi^{-1}(T)$ is the inverse  
image under $\phi_k$ of the tree  $gT \subset \G$.

\med
Now let  $
\eps >0$ be chosen so that for every point $x$, if $B_\eps(T_xM)$ is  
the ball centered at the origin of radius $\eps$, then the  
exponential map,
$$
exp : B_\eps(T_x M) \to M
$$
is a diffeomorphism onto its image.  Let $U_\eps (x)$  be this image.
Consider a point $(\sigma, \gamma ) \in \eta_\eps$.   Notice that  
each $y_i \in U_\eps(x_i)$.  Thus there is a unique $u_i \in B_\eps(T_ 
{x_i}M)$ with $exp (u_i) = y_i$.  The assignment $(\sigma, \gamma)  
\to (u_1, \cdots, u_k)$ defines an $\ag$-equivariant section
$$
\theta : \eta_\eps \to \zeta
$$
of the restriction $\zeta_{|_{\eta_\eps}} \to \eta_\eps$.
For each $i$, the curve $t \to exp(tu_i)$ in $M$ is  a path from $x_i 
$   at $t=0$,  to $y_i$ at $t=1$.  This is a gradient flow line
of the distance function $d_{x_i} : M \to \br$, defined to be the  
distance from $x_i$,    $$d_{x_i}(x) = d(x_i, x).$$
This allows us to construct a morphism $\psi_{k} : \G_{k+1} \to \G_k$  
in $\cc_\G$, as follows.  Let $\G_{k+1}$ be the graph obtained
from $\G_k$ by  replacing each vertex $s^T_i(\G_k)$ with an edge of  
length 1.  The morphism $\psi_k$ collapses each of these edges to a  
point.  If we label these new edges by the functions $ d_{x_1},  
\cdots, d_{x_k}$, we have now defined a new structure $\sigma^\prime$.
Notice that the element $(\sigma^\prime, p(\sigma, \gamma)) \in \cs_ 
\G \times M \cong \tilde \cm_{tree}(\ag, M)$ lies in the image of $ 
\tccgm \hk \tilde \cm_{tree}(\ag, M)$.  This is because the  
coordinate $y_i$ assigned to the pair $(\sigma^\prime, p(\sigma,  
\gamma))$ is the same as the $y_i$ coordinate for the pair $(\sigma,  
\gamma)$.  But the $x_i$ coordinate assigned to the pair
$(\sigma^\prime, p(\sigma, \gamma))$ is equal to $exp(u_i) = y_i$.   
Thus the projection
$$
p_T (\sigma^\prime, p(\sigma, \gamma)) \in (M^2)^b
$$
lies in the image of $\Delta^b : M^b \subset (M^2)^b$.  By the  
pullback square (\ref{pullback}), the pair $(\sigma^\prime, p(\sigma,  
\gamma))$  lies in the image of $\tccgm$.
Sending $(\sigma, \gamma)$ to $  (\sigma^\prime, p(\sigma, \gamma))$
defines an map
$$
\pi : \eta_\eps \to \tccgm.
$$
We now show that the  section $\theta : \eta_\eps \to \zeta$ defines  
an equivariant lifting $\Theta : \eta_\eps \to \nu(\iota)$ making the  
following diagram commute:

\begin{equation}\label{Theta}
\begin{CD}
\eta_\eps @>\Theta >> \nu (\iota) \\
@V=VV  @VVV \\
\eta_\eps @>>\pi > \tccgm
\end{CD}
\end{equation}
The lifting $\Theta$ is defined as follows:

Consider the unique geodesic path in the tree $T_1$ from $s^T_i(\G_k) 
$ to the fixed vertex $v \in T_1$.   Then its image under the tree  
flow $\gamma_{T_1}$ is a parameterized curve from $x_i$ to $\gamma_ 
{T_1}(v) = p(\sigma, \gamma)$ in $M$.  (Recall that $ p : \tccgm \to M 
$ maps $(\sigma, \gamma)$ to $\gamma (v)$.)   Using the Levi-Civita  
connection, we define
$w_i \in T_{p(\sigma, \gamma)}M$ to be the image of $u_i  \in T_{x_i}M 
$ under the parallel transport operator along this path:
$$
\tau_{\gamma_{T_1}}: T_{x_i}M) \xr{\cong} T_{p(\sigma, \gamma)}M.
$$

This construction defines an $\ag$-equivariant  map
\begin{align}
\Theta : \eta_\eps &\to p^*(\bigoplus_b TM) = \nu (\rho) \\
(\sigma, \gamma) &\to (w_1, \cdots , w_b) \in \bigoplus_b T_{p 
(\sigma, \gamma)}M \notag
\end{align}
making the diagram (\ref{Theta}) commute.

\med
We claim that $\Theta$ is a homeomorphism. One can see this be  
directly constructing an inverse map
$$
\Theta^{-1} : \nu (\rho) \to \eta_\eps.
$$
This is constructed as follows.  Given $(u_1, \cdots, u_b) \in  
\bigoplus_b T_{p(\sigma, \gamma)}M$ where $(\sigma, \gamma) \in \tccgm 
$, one can parallel translate along geodesic paths in $T_1(\G)k)$ to  
obtain the vector $(w_1, \cdots , w_b) \in B_\eps (T_{x_1}M \oplus   
\cdots \oplus T_{x_b}M$. By scaling these vectors, if necessary, one  
can consider the curves $t \to exp(-tu_i)$
to define a new structure $\sigma^{''}$ so that the point $(\sigma^ 
{''}, p(\sigma, \gamma))$ lives in $\eta_\eps \subset \cs (\G) \times  
M \cong \tilde \cm_{tree}(\G, M)$.  Notice that the coordinates $\{(x_i, y_i)\}$ 
associated to $(\sigma^{''}, p(\sigma, \gamma))$ are  
the points $(exp(-u_i), x_i)$ where $(x_1, \cdots, x_b) \in M^b$ is  
$p_T(\sigma, \gamma)$ as in diagram (\ref{pullback}).  The assignment  
$(u_1, \cdots, u_b) \to (\sigma^{''}, p(\sigma, \gamma)$ defines a  
map $\nu (\rho) \to \eta_{\eps}$ which is easily checked to be  
inverse to $\Theta$.

Thus $\Theta : \eta_\eps \to \nu(\rho)$ is an equivariant  
homeomorphism, and so induces a homeomorphism on orbit spaces.
This completes the proof of theorem (\ref{tube}).

\section{ Appendix: Regularity}

\begin{lemma}
     On the interior of any edge $E\subset\g_k$, an element $r\in\coker
     D_{\g_k}$ is smooth and satisfies
     \begin{equation} \label{eq:coker}
     \dot{r}_E-(\nabla\nabla f_E)^T\cdot r_E=0.
     \end{equation}
\end{lemma}
\begin{proof}
     Take an open interval $I=(t_1,t_2)\subset E$ and choose $\phi\in
     C^{\infty}_0(I)$.  Trivialise $\cv=\gamma^*TM$ over $I$ and put
     $\nabla\nabla f_E=A(t)$ with respect to this trivialisation.  Then
     $\langle r,\dot{\phi}+A\phi\rangle_2=0$.  Now
     $\phi(t)=\int^t_{t_0}\dot{\phi}(\tau)d\tau$ so \[\int_I\langle
     r(t),\dot{\phi}(t)\rangle dt+\int_I\langle A^T(t)r(t),
     \int^t_{t_0}\dot{\phi}(\tau)d\tau\rangle dt=0\] and by Fubini's
     theorem \[\int_I\langle r(\tau),\dot{\phi}(\tau)\rangle
     d\tau+\int_I\int^{t_1}_{\tau} \langle
     A^T(t)r(t),\dot{\phi}(\tau)\rangle dtd\tau=0\] Thus
     \[\int_I\langle
     r(\tau)-\int^{\tau}_{t_1}A^T(t)r(t)dt,\dot{\phi}(\tau)\rangle
     d\tau=0\ {\rm for\ all}\ \phi\in C^{\infty}_0(I)\ .\] Since
     $\dot{\phi}$ has mean zero and the set of such functions is dense
     in $L^2(I)$ we have \[ r(\tau)-\int^{\tau}_{t_1}A^T(t)r(t)dt={\rm
     constant}\ .\] This integral equation supplies us with information
     about the behaviour of $r$ in $I$.  For a start it says that $r$
     is absolutely continuous with derivative equal to the integrand
     almost everywhere.  At points of continuity of $A$ the derivative
     of $r$ is equal to the integrand.  Furthermore, regularity of $A$
     gives regularity of $r$.  This can be seen as follows.  At a point
     $\tau_0$ of continuity of $A$
     \[\left|\frac{1}{2\delta}\int^{\tau_0+\delta}_{\tau_0-\delta}A^T 
(t)r(t)dt-
     A^T(\tau_0)r(\tau_0)\right|\leq\epsilon M\] where
     $\epsilon=\sup_{(\tau_0-\delta,\tau_0+\delta)}\{ |A(t)|,|r(t)|\}$
     tends to zero as $\delta$ tends to zero since $A(t)$ and $r(t)$
     are continuous at $\tau_0$.  This shows that the derivative of $r$
     exists there and
     \begin{equation} \label{eq:rder}
     \dot{r}(\tau_0)=A^T(\tau_0)r(\tau_0)
     \end{equation}
     If $A$ is differentiable in a neighbourhood of $\tau_0$ then by
     (\ref{eq:rder})
     \[\ddot{r}(\tau)=(\dot{A}^T(\tau)+A^T(\tau)^2)r(\tau)\] in that
     neighbourhood, and so on.  Thus $r_E$ satisfies (\ref{eq:coker}).
\end{proof}
\section{Appendix: The Fredholm operator}

Let $D_{\g_k}$ be the linearisation of the graph flow equation along
the graph flow $\gamma:\g_k\to M$ of the compact metric graph $\g_k$
so $D_{\g_k}s=\dot{s}_E+\nabla\nabla f_E\cdot s_E$ for $s$ a section
of $\cv=\gamma^*TM$.

\begin{proposition}
     $D_{\g_k}:W^{1,2}(\g_k,\cv)\to L^2(\g_k,\cv)$ is Fredholm.
\end{proposition}
\begin{proof}
     Put $D_{\g_k}s=\dot{s}_E+\nabla\nabla f_E\cdot
     s_E=\dot{s}_E(t)+A(t)s_E(t).$
     \[\int_{\g_k}|\dot{s}+As|^2dt=\int_{\g_k}\left(\frac{1}{2}|
     \dot{s}+2As|^2+\frac{1}{2}|\dot{s}|^2-|As|^2\right)dt
     \geq\int_{\g_k}\left(\frac{1}{2}|\dot{s}|^2-|As|^2\right)dt\ .\]
     Thus using $|A(t)\cdot s(t)|\leq\|A(t)\|\cdot|s(t)|$ and setting
     $c_A={\rm max}_{\g_k}\|A(t)\|$, we have
     \[\int_{\g_k}|\dot{s}+As|^2dt\geq\frac{1}{2}\int_{\g}
     |\dot{s}|^2-c_A\int_{\g_k}|s|^2dt\ .\]
     Hence there is a $c>0$ satisfying
     \[\int_{\g_k}(|s|^2+|\dot{s}|^2)dt\leq c\int_{\g_k}(|s|^2+
     |\dot{s}+As|^2)dt\ .\]
     In other words,
     \begin{equation}  \label{eq:semiF}
     \|s\|_{W^{1,2}(\g_k)}\leq c(\|s\|_{L^2(\g_k)}+
     \|D_{\g_k}s\|_{L^2(\g_k)}).
     \end{equation}
     It is a rather standard consequence of (\ref{eq:semiF}) that
     $D_{\g_k}$ is semi-Fredholm, meaning that it has
     finite-dimensional kernel and closed range (see
     \cite{SalSei,SchMor} for example).  This can be seen as follows.

     The operator
     \[ K:W^{1,2}(\g_k)\stackrel{\rm cpt.}\hookrightarrow L^2(\g_k)\]
     is compact by Rellich's lemma.  Thus the image under
     $K$ of any bounded sequence in the kernel of $D_{\g_k}$ has a
     convergent subsequence which is necessarily Cauchy.  The
     inequality (\ref{eq:semiF}) then implies that the subsequence is
     Cauchy in $W^{1,2}(\g_k)$.  Thus the unit ball in the kernel of
     $D_{\g_k}$ is compact, showing that the kernel is
     finite-dimensional.

     To show that the image is closed, consider a bounded sequence $\{
     s_i\}\subset W^{1,2}(\g_k)$ such that $\{ D_{\g_k}s_i\}$ is
     Cauchy in $L^2(\g_k)$.  Choose a subsequence $\{ s_{i_j}\}$ such
     that $\{ Ks_{i_j}\}$ is Cauchy in $L^2(\g_k)$.  It follows from
     (\ref{eq:semiF}) that $\{ s_{i_j}\}$ is Cauchy thus converging to
     $s$, say.  Since $D_{\g_k}$ is continuous, $\{ D_{\g_k}s_i\}$
     converges to $D_{\g_k}s$.  In fact, the sequence $\{ s_i\}$ can
     be arranged to be bounded as follows.  By the Hahn-Banach theorem
     there exists a closed subspace $U\subset W^{1,2}(\g_k)$
     satisfying \[\ker D_{\g_k}\oplus U=W^{1,2}(\g_k)\ .\] Project
     $\{ s_i\}$ onto $\{\tilde{s}_i\}\subset U$.  This has to be
     bounded since otherwise a subsequence of
     $\{\tilde{s}_i/\|\tilde{s}_i\|\}$ converges to $s\in U$ with
     $\|s\|=1$ and $D_{\g_k}s=0$ in contradiction to the construction
     of $U$.  Thus $D_{\g_k}$ has closed range.

     To complete the proof of the proposition we must show that $\coker
     D_{\g_k}$ is finite-dimensional.  In Lemma~\ref{th:coker} it was
     shown that elements $r$ of the cokernel satisfy the differential
     equation $D_{\g_k}^*r=0$ which is much like the equation
     $D_{\g_k}s=0$, the only difference being that $r$ need not be
     continuous at the vertices of $\gn_k$.  Nevertheless, as for $\ker
     D_{\g_k}$ the unit ball in the kernel of $D_{\g_k}^*$ is compact
     and the dimension of $\coker D_{\g_k}$ is finite.  Hence
     $D_{\g_k}$ is Fredholm.
\end{proof}
{\em Remark.} Since elements of the kernel and cokernel are smooth an
alternative proof of finite-dimensionality follows from uniqueness of
solutions to ODEs.  Still, to prove Fredholmness one must show that
the image of $D_{\g_k}$ is closed and there is no smoothness here to
work with.  This is why we used standard Banach space arguments rather
than the more intuitive uniqueness of solutions to ODEs.\\

To prove that $D_{\gn_k}$ is Fredholm for non-compact $\gn_k$
requires a further argument.

\begin{proposition}
     $D_{\gn_k}:W^{1,2}(\gn_k,\cv)\to L^2(\gn_k,\cv)$ is Fredholm.
\end{proposition}
\begin{proof}
     For any $T>0$, construct the compact graph $\g_k^T$ lying between
     $\g_k\subset\g_k^T\subset\gn_k$ by cutting $\gn_k$ off at the
     parameter $T$ on outgoing edges and $-T$ on incoming edges.  The
     proof uses the following lemma.
     \begin{lemma}  \label{th:estimate}
     For large enough $T$, there exists $c=c(T)$ such that
     \begin{equation}  \label{eq:semiF}
         \|s\|_{W^{1,2}(\gn_k)}\leq c(\|s\|_{L^2(\g_k^T)}+
         \|D_{\gn_k}s\|_{L^2(\gn_k)}).
     \end{equation}
     \end{lemma}
     \begin{proof}
     Put $$D_{\gn_k}=\frac{d}{dt}+\nabla\nabla
     f_E=\frac{d}{dt}+A(t)$$
     with respect to a trivialisation of $\cv$ over $\gn_k$ with
     transition functions $\pm 1$.  Since $f$ is Morse,
     $\lim_{t\to\infty}A(t)$ is invertible.

     The estimate of $\|s\|_{W^{1,2}(\gn_k)}$ breaks into one part
     near infinity and a compact part.  Near infinity, the graph
     flow equation is just the usual gradient flow equation so we
     can use a result whose proof can be found in \cite{SchMor}.
     Given $A(t)$ with $\lim_{t\to\infty}A(t)$ invertible, there are
     constants $T>0$, $c_1(T)>0$ such that
     \[\|s\|_{W^{1,2}(\gn_k)}\leq c_1(T)\|\dot{s}+A(t)s\|_2\ \ {\rm
     for\ all}\ s\in W^{1,2}(\gn_k),\ s_|{\g_k^T}=0.\]

     For the compact part use the previous proposition applied to
     $\g_k^T$.
     \end{proof}
     To put together the part near infinity and the compact part, define
     a cut-off function $\beta\in C^{\infty}(\gn_k,[0,1])$ with the
     properties
     \[\beta|_{\g_k^T}=1,\ \ \beta(t)=0\ {\rm for}\ |t|\geq T+1,\
     {\rm and}\ \dot{\beta}(t)\neq 0\ {\rm for}\ |t|\in (T,T+1).\]
     In the following, put $\|\cdot\|_{L^2(\gn_k)}=\|\cdot\|_2$ and
     $\|\cdot\|_{W^{1,2}(\gn_k)}=\|\cdot\|_{1,2}$ and only specify the
     compact graph in the norm. Also choose $c_4$ and $c$ large enough.
     \begin{eqnarray*}
     \|s\|_{1,2}&=&
     \|\beta s+(1-\beta)s\|_{1,2} \leq\|\beta
     s\|_{1,2}+ \|(1-\beta)s\|_{1,2}\\
     &\leq& c_4(\|\beta s\|_{L^2(\g_k^T)}+\|D_{\gn_k}(\beta
     s)\|_2+ \|D_{\gn_k}((1-\beta)s)\|_2)\\
     &\leq&c_4(\|\beta s\|_{L^2(\g_k^T)}
     +2\|\dot{\beta}s\|_2+\|D_{\gn_k}s\|_2)\\
     &\leq&c(\|s\|_{L^2(\g_k^T)}+\|D_{\gn_k}s\|_2)\ .
     \end{eqnarray*}
     and the proof goes through as in the compact case.
\end{proof}

\end{document}